\PassOptionsToPackage{hyperfootnotes=false}{hyperref}
\documentclass[a4paper]{cas-sc}

\usepackage[utf8]{inputenc}
\usepackage[numbers,sort&compress]{natbib}
\usepackage{amsthm,mathtools,mathrsfs,bm}
\usepackage{microtype}
\usepackage{enumitem}
\usepackage{verbatim}
\usepackage{aliascnt}
\usepackage[capitalise,nameinlink,noabbrev]{cleveref}

\setlist[enumerate]{itemsep=2pt,topsep=4pt}
\setlist[itemize]{itemsep=2pt,topsep=4pt}
\allowdisplaybreaks
\numberwithin{equation}{section}
\hypersetup{hidelinks}
\ExplSyntaxOn
\keys_set:nn { stm / mktitle } { nologo = true }
\RenewDocumentEnvironment { Abstract } { o }
  {
    \par\noindent
    \box_use_drop:N \g_stm_key_box
    \hfill
    \begin{minipage}[t]{.65\textwidth}
      \IfNoValueF { #1 } { \def\abstractname{#1} }
      \parindent=0pt
      \noindent\abstractname\par
      \vspace{-4pt}
      \noindent\rule{\linewidth}{.2pt}\par\footnotesize
      \ignorespaces
  }
  { \par\end{minipage}\par }
\ExplSyntaxOff

\theoremstyle{plain}
\newtheorem{theorem}{Theorem}[section]

\newaliascnt{lemma}{theorem}
\newtheorem{lemma}[lemma]{Lemma}
\aliascntresetthe{lemma}

\newaliascnt{proposition}{theorem}

\aliascntresetthe{proposition}

\newaliascnt{corollary}{theorem}

\aliascntresetthe{corollary}

\theoremstyle{definition}

\newaliascnt{definition}{theorem}
\newtheorem{definition}[definition]{Definition}
\aliascntresetthe{definition}

\newaliascnt{assumption}{theorem}
\newtheorem{assumption}[assumption]{Assumption}
\aliascntresetthe{assumption}

\newaliascnt{example}{theorem}

\aliascntresetthe{example}

\newaliascnt{remark}{theorem}
\newtheorem{remark}[remark]{Remark}
\aliascntresetthe{remark}

\DeclareMathOperator*{\esssup}{ess\,sup}

\DeclareMathOperator{\sgn}{sgn}
\DeclareMathOperator{\tr}{tr}

\crefname{theorem}{Theorem}{Theorems}
\Crefname{theorem}{Theorem}{Theorems}
\crefname{lemma}{Lemma}{Lemmas}
\Crefname{lemma}{Lemma}{Lemmas}
\crefname{proposition}{Proposition}{Propositions}
\Crefname{proposition}{Proposition}{Propositions}
\crefname{corollary}{Corollary}{Corollaries}
\Crefname{corollary}{Corollary}{Corollaries}
\crefname{definition}{Definition}{Definitions}
\Crefname{definition}{Definition}{Definitions}
\crefname{assumption}{Assumption}{Assumptions}
\Crefname{assumption}{Assumption}{Assumptions}
\crefname{example}{Example}{Examples}
\Crefname{example}{Example}{Examples}
\crefname{remark}{Remark}{Remarks}
\Crefname{remark}{Remark}{Remarks}
\crefname{equation}{equation}{equations}
\Crefname{equation}{Equation}{Equations}
\crefname{section}{Section}{Sections}
\Crefname{section}{Section}{Sections}

\newcommand{\R}{\mathbb{R}}
\newcommand{\E}{\mathbb{E}}
\newcommand{\Pp}{\mathbb{P}}

\newcommand{\Bcal}{\mathcal{B}}
\newcommand{\Pcal}{\mathcal{P}}
\newcommand{\Fcal}{\mathcal{F}}
\newcommand{\Ucal}{\mathcal{U}}
\newcommand{\Aop}{\mathcal{A}}
\newcommand{\Bop}{\mathcal{B}}
\newcommand{\Cop}{\mathcal{C}}

\newcommand{\dd}{\mathop{}\!\mathrm{d}}
\newcommand{\dt}{\mathop{}\!\mathrm{d}t}
\newcommand{\ds}{\mathop{}\!\mathrm{d}s}
\newcommand{\dr}{\mathop{}\!\mathrm{d}r}
\newcommand{\de}{\mathop{}\!\mathrm{d}e}
\newcommand{\dP}{\mathop{}\!\mathrm{d}\mathbb{P}}
\newcommand{\dW}{\mathop{}\!\mathrm{d}W}

\DeclarePairedDelimiterX{\ip}[2]{\langle}{\rangle}{#1, #2}

\DeclarePairedDelimiter{\abs}{\lvert}{\rvert}
\newcommand{\1}{\mathbf{1}}

\newcommand{\tilmu}{\tilde{\mu}}
\newcommand{\assitemref}[1]{\Cref{#1}}

\begin{document}

\shorttitle{HJB integro-differential equations  with non-Lipschitz running costs}
\shortauthors{Z. Wu and Q. Zhang}
\title[mode=title]{Viscosity Solutions of Hamilton--Jacobi--Bellman
Integro-Differential Equations with Non-Lipschitz Running Costs}

\tnotemark[1]
\tnotetext[1]{This paper is supported by the National
Key R\&D Program of China (No.~2022YFA1006101)
and the National Natural Science Foundation of
China (No.~12371445).}

\author[1]{Zihao Wu}[orcid=0009-0008-3426-9661]
\ead{zhwu22@m.fudan.edu.cn}

\author[1,2]{Qi Zhang}[orcid=0000-0002-8225-1211]
\cormark[1]
\ead{qzh@fudan.edu.cn}

\affiliation[1]{
    organization={School of Mathematical Sciences, Fudan University},
    addressline={220 Handan Road},
    city={Shanghai},
    postcode={200433},
    country={China}
}
\affiliation[2]{
    organization={Laboratory of Mathematics for Nonlinear Science, Fudan University},
    addressline={220 Handan Road},
    city={Shanghai},
    postcode={200433},
    country={China}
}
\cortext[1]{Corresponding author.}

\begin{abstract}
We study a class of parabolic Hamilton--Jacobi--Bellman (HJB) 
integro-differential equations arising from stochastic recursive
control problems with jumps. These equations involve a
L\'evy-type operator and an additional nonlocal term in the running cost. The underlying backward generator satisfies
a monotonicity condition with polynomial growth in the backward
variable and is uniformly continuous, but not necessarily Lipschitz continuous, 
in the state, Brownian, jump and control variables.
We establish well-posedness, comparison and stability results
for the corresponding backward stochastic differential equations
with jumps, and derive a dynamic programming principle for the
associated recursive control problem. Using approximation
arguments and the stability of viscosity solutions, we show
that the value function provides a viscosity solution of the
HJB integro-differential equation. We further establish a comparison principle for
bounded sub- and supersolutions and uniqueness within a
prescribed growth class, under an additional
uniform bound on the jump amplitudes.
\end{abstract}


\maketitle
\hypersetup{
    pdftitle={Viscosity Solutions of Hamilton--Jacobi--Bellman
Integro-Differential Equations with Non-Lipschitz Running Costs},
    pdfauthor={Zihao Wu and Qi Zhang},
    pdfsubject={Non-Lipschitz HJB integro-differential equations and stochastic recursive control problems},
    pdfkeywords={HJB integro-differential equations; viscosity solutions; comparison principle; non-Lipschitz generators; BSDEs with jumps; stochastic recursive control}
}

\section{Introduction}\label{sec:introduction}

Hamilton--Jacobi--Bellman (HJB) integro-differential equations
arise naturally in stochastic optimal control problems involving
jump-diffusion processes. In contrast to the purely diffusive
setting, the presence of jumps leads to nonlocal operators and
requires an appropriate notion of viscosity solution. A further
nonlinear structure appears when the cost functional is defined
recursively through a backward stochastic differential equation
(BSDE) with jumps. In this case, the associated HJB equation
contains not only the usual L\'evy-type operator, but also an
additional nonlocal term in the running cost. More
precisely, we consider a class of parabolic HJB integro-differential equations of
the form
\begin{equation}\label{eq:hjb-intro}
\left\{
\begin{aligned}
\partial_t u(t,x)
+\sup_{v\in U}\Bigl\{
&\frac{1}{2}\tr\bigl( \sigma \sigma^\top(t,x,v)D^2 u (t,x) \bigr) + \langle b(t,x,v),D u (t,x) \rangle\\
&+\int_E{( u (t,x+\gamma(t,x,v,e)) - u (t,x) - \ip{D u (t,x)}{\gamma(t,x,v,e)} ) \lambda(\de)}\\
&+f\left(t,x,u(t,x),
\sigma^\top(t,x,v)Du(t,x),
\int_E{( u (t,x+\gamma(t,x,v,e)) - u (t,x) )}\ell(t,e)\lambda(\de),v\right)
\Bigr\}=0,\\
&\hspace{-10mm}(t,x)\in[0,T)\times\mathbb R^{d'},\\
u(T,x)&=h(x),
\qquad x\in\mathbb R^{d'}.
\end{aligned}
\right.
\end{equation}
Our objective is to investigate the existence, comparison and
uniqueness of viscosity solutions of this equation under
non-Lipschitz conditions.

The probabilistic structure underlying \eqref{eq:hjb-intro}
is provided by stochastic recursive control problems with
jumps. In such problems, the state evolves according to a
controlled jump-diffusion, while the cost functional is
determined by the solution of a BSDE driven by both Brownian
motion and a compensated Poisson random measure.
The value function is
obtained by optimizing the cost functional over
admissible controls. The dynamic programming principle then
provides a probabilistic
representation for the associated HJB integro-differential equation. This
demonstrates that BSDEs with jumps are a natural tool to construct viscosity solutions of nonlinear
integro-differential equations.

The relationship between BSDEs and nonlinear partial
differential equations, well-known as the nonlinear Feynman-Kac formula, originated in the work of Peng \cite{Peng1991}. In the presence of jumps,
Tang and Li \cite{TangLi1994} developed stochastic control
methods for systems with random jumps, while Barles,
Buckdahn and Pardoux \cite{BarlesBuckdahnPardoux1997}
established the connection between BSDEs with jumps and
integro-differential equations through viscosity solutions.
On the analytical side, the fundamental viscosity solution
framework of Crandall, Ishii and Lions
\cite{CrandallIshiiLions1992} has been extended to nonlocal
equations in various settings; see, for example,
Alvarez and Tourin \cite{AlvarezTourin1996}, Barles and
Imbert \cite{BarlesImbert2008}, and Jakobsen and Karlsen
\cite{JakobsenKarlsen2005}. These developments provide
essential tools for handling nonlocal operators, stability
properties and comparison principles. For stochastic
recursive control problems, the backward semigroup approach
introduced by Peng \cite{Peng1992,Peng1997} established a
dynamic programming framework for nonlinear cost
functionals. Such recursive formulations are also closely
related to stochastic differential utility and nonlinear
evaluations; see Duffie and Epstein \cite{DuffieEpstein1992}
and El Karoui, Peng and Quenez
\cite{ElKarouiPengQuenez1997}. In the jump-diffusion
setting, Li and Peng \cite{LiPeng2009} proved the dynamic
programming principle and characterized the value function
as a viscosity solution of the corresponding HJB
integro-differential equation under Lipschitz conditions
on the generator.

Another important direction concerns the relaxation of
Lipschitz assumptions on BSDE generators. Such assumptions
may be restrictive for nonlinear evaluations and may also
be stronger than those naturally imposed on the nonlinear
terms in viscosity solution theory. Various existence,
uniqueness and stability results for BSDEs with continuous,
monotone or polynomial-growth generators have been obtained
by Lepeltier and San Martin
\cite{LepeltierSanMartin1997}, Pardoux \cite{Pardoux1999},
Briand and Carmona \cite{BriandCarmona2000}, and Ma,
Fan and Song \cite{MaFanSong2013}, among others.
Uniformly continuous dependence on the Brownian integrand
has also been studied in this context; see Jia
\cite{Jia2008}. For equations driven by Brownian and
Poisson noises, Kruse and Popier \cite{KrusePopier2016}
investigated BSDEs with monotone generators in a general
filtration. Concerning stochastic recursive control
problems without jumps, Pu and Zhang \cite{PuZhang2018}
established the dynamic programming principle and the
associated viscosity solution property for generators
which are monotone with polynomial growth in $y$ and independent of $z$.
Subsequently, Zhuo, Dong and Pu \cite{ZhuoDongPu2020}
extended the analysis to generators with uniformly
continuous dependence on $z$ and
obtained the corresponding existence and uniqueness
results for viscosity solutions of HJB equations.

The preceding results provide two closely related, but
distinct, approaches to the equations under consideration.
The recursive control theory with jumps in
\cite{LiPeng2009} relies on Lipschitz conditions on the
BSDE generator, whereas the non-Lipschitz recursive control
results in \cite{PuZhang2018,ZhuoDongPu2020} are developed
in the Brownian setting. The simultaneous presence of
jumps and non-Lipschitz nonlinearities therefore requires
additional arguments, both for the probabilistic
representation and for the viscosity solution theory of
the resulting nonlocal HJB equation. In the present paper,
we consider a generator \(f\) which satisfies a
monotonicity condition with polynomial growth in \(y\) 
and is uniformly continuous with a common modulus of at most linear
growth in \(x,z,k\) and \(v\). Moreover, \(f\) is assumed to be nondecreasing
in \(k\), which is essential for the
comparison arguments involving the jump component.

The first difficulty concerns the probabilistic
representation under these weaker regularity assumptions.
The usual Lipschitz estimates for BSDEs cannot be applied
directly. 
Furthermore,
comparison in the presence of jumps requires the
dependence on the aggregated jump variable to be handled
consistently with the compensator of the Poisson random
measure. To overcome these difficulties, we employ
inf-convolution approximations together with a Girsanov
transformation involving both Brownian and Poisson
components. The resulting comparison, a priori and
stability estimates provide the continuous dependence
properties required for the recursive control problem.
These properties
allow us to establish the dynamic programming principle
through the backward stochastic semigroup method.

The second difficulty lies in the viscosity solution
theory of \eqref{eq:hjb-intro}. 
The simultaneous presence of the two nonlocal operators 
must be accounted for in the
viscosity formulation, the approximation procedure and
the comparison argument. For the existence result, we
combine truncation and inf-convolution approximations
of the generator with a stability theorem for viscosity
solutions of the associated nonlocal equations. For
comparison and uniqueness, we use the nonlocal
Jensen--Ishii framework of Barles and Imbert
\cite{BarlesImbert2008}, together with estimates for the
two nonlocal terms and the modulus of continuity of the
generator. To treat solutions with unbounded growth,
we further employ a reference-function argument based
on the method of Barles, Buckdahn and Pardoux
\cite{BarlesBuckdahnPardoux1997}.

The main result of this paper is the existence and uniqueness of viscosity
solutions of \eqref{eq:hjb-intro}. Under the standing
assumptions, we first construct a continuous viscosity
solution through the value function of the associated
stochastic recursive control problem with jumps. The
construction rests on the well-posedness and stability
results for the corresponding BSDEs, the regularity of
the value function, and the dynamic programming
principle. By approximating the non-Lipschitz generator
with Lipschitz generators and passing to the limit through
the stability of viscosity solutions, we obtain the
desired viscosity characterization. We then establish
a comparison principle for continuous viscosity
subsolutions bounded from above and supersolutions
bounded from below, provided that the terminal
ordering holds. Finally, under an additional uniform bound
on the jump amplitudes, we extend the uniqueness
result to a prescribed class of continuous functions
with unbounded growth. Consequently, in this setting,
the probabilistically constructed value function is the
unique viscosity solution within that growth class.

The remainder of the paper is organized as follows.
In \Cref{sec:preliminaries}, we introduce the controlled
forward-backward system with jumps, state the standing
assumptions and establish the necessary well-posedness,
comparison and stability results for BSDEs with
non-Lipschitz generators. In \Cref{sec:dpp}, we prove
the determinism and continuity of the value function
and establish the dynamic programming principle.
In \Cref{sec:hjb}, we formulate the HJB
integro-differential equation, discuss the viscosity
solution formulations and their stability, and prove
the viscosity solution property of the value function.
Finally, in \Cref{sec:uniqueness}, we establish the
comparison principle and prove uniqueness in the
prescribed growth class.

\section{Preliminaries}\label{sec:preliminaries}

Let \((\Omega,\Fcal,\Pp)\) be a complete probability space
    carrying a \(d\)-dimensional standard Brownian motion \(W\) and a Poisson random measure \(\mu\) on \(\R_{+}\times E\) with compensator \(\lambda(\de) \dt\),
    where \(E=\R^{d_e}\setminus\{0\}\), for a fixed integer \(d_e\ge1\), is equipped with its Borel \(\sigma\)-algebra \(\Bcal(E)\), \(W\) and \(\mu\) are mutually independent, and \(\lambda\) is a \(\sigma\)-finite L\'evy measure on \((E,\Bcal(E))\) satisfying \(\int_E{(1\wedge |e|^2)\lambda(\de)}<\infty\).
    Denote by \(\tilmu(\dr,\de)=\mu(\dr,\de)-\lambda(\de) \dr\) the compensated Poisson random measure.

    For a given \(T>0\), let \(\mathbb{F}=\{ \Fcal_t \}_{0\le t\le T} \) be the filtration generated by \(W\) and \(\mu\), augmented 
by \(\mathcal{N}\), the set of all \(\Pp\)-null sets.
Let \(\Pcal\) denote the predictable \(\sigma\)-algebra on \(\Omega \times [0,T]\) associated with \(\mathbb{F}\).

For \(p\ge 2\), \(0\leq t\leq T\) and a Hilbert space \(H\), we introduce the following spaces, which will be used throughout this paper:
{\small\begin{align*}
    &L^{p}_{\Fcal_t}(\Omega ; H)=\bigl\{ \xi : H\text{-valued}\ \Fcal_t\text{-measurable}\ \text{random}\ \text{variable},\ \E[ |\xi|^{p} ] <\infty \bigr\}; \\
    &\mathcal{S}^{p}(t,T;H)=\Bigl\{ (\psi_s)_{s\in [t,T]} : H\text{-valued}\ \mathbb{F} \text{-adapted}\ \text{c\`adl\`ag}\ \text{process},\ \E\Bigl[ \sup_{s\in [t,T]}|\psi_s|^{p} \Bigr] <\infty \Bigr\}; \\
    &\mathcal{L}^{p}(t,T;H)=\Biggl\{ (\psi_s)_{s\in [t,T]} : H\text{-valued}\ \mathbb{F} \text{-progressively}\ \text{measurable}\ \text{process},\ \E\Biggl[ \biggl( \int_t^T{|\psi_s|^2\ds} \biggr) ^{p/2} \Biggr] <\infty \Biggr\}; \\
    &\mathcal{L}^{p}_{\lambda}(t,T;H)=\Biggl\{ K : \Omega \times [t,T]\times E\rightarrow H, \ \Pcal\otimes \Bcal(E)\text{-measurable},\ \E\Biggl[ \biggl( \int_t^T{\int_E{|K_s( e ) |^2\lambda(\de) \ds}} \biggr) ^{p/2} \Biggr] <\infty \Biggr\}; \\
    &L^2(E;H)=\biggl\{ K : E\rightarrow H,\ \Bcal(E)\text{-measurable},\ \int_E{|K( e ) |^2\lambda(\de)} <\infty \biggr\}.
\end{align*}}
If \(t\), \(T\) and \(H\) are clear from the context, we denote these spaces by \(L^{p}_{\Fcal_t}\), \(\mathcal{S}^{p}\), \(\mathcal{L}^{p}\) and \(\mathcal{L}^{p}_{\lambda}\), respectively.

Let \(U\subset \R^{d_v}\) be a nonempty compact set. For \(t\in [0,T]\), let \(\Ucal_{t,T}\) be the set of all \(U\)-valued \(\mathbb{F}\)-predictable processes on \([t,T]\).
We refer to elements of \(\Ucal_{t,T}\) as admissible control processes.

For \(t\in [0,T],\ \eta \in L^{p}_{\Fcal_t}(\Omega ; \R^{d'}),\ v\in \Ucal_{t,T}\),
we consider the following forward-backward stochastic differential equation (FBSDE) with jumps defined on \([t,T]\):
\begin{align}\label{eq:controlled-fbsde}
    \begin{dcases}
        X^{t,\eta;v}_s&=\eta +\int_t^s{b\bigl( r,X^{t,\eta;v}_r,v_r \bigr) \dr} +\int_t^s{\sigma \bigl( r,X^{t,\eta;v}_r,v_r \bigr) \dW_r} \\
        &\quad +\int_t^s{\int_E{\gamma\bigl( r,X^{t,\eta;v}_{r-},v_r,e \bigr) \tilmu( \dr,\de )}},\ \ \ s\in [t,T]; \\
        Y^{t,\eta;v}_s&=h\bigl(X^{t,\eta;v}_T\bigr)+\int_s^T{f\biggl( r,X^{t,\eta;v}_r,Y^{t,\eta;v}_r,Z^{t,\eta;v}_r,\int_E{\ell_r(e)K^{t,\eta;v}_r(e) \lambda(\de)},v_r \biggr) \dr} \\
        &\quad -\int_s^T{\ip[\big]{Z^{t,\eta;v}_r}{\dW_r}}-\int_s^T{\int_E{K^{t,\eta;v}_r( e ) \tilmu( \dr,\de )}},\ \ \ s\in [t,T],
    \end{dcases}
\end{align}
where \(b:[0,T]\times \R^{d'}\times U\rightarrow \R^{d'}\), \(\sigma :[0,T]\times \R^{d'}\times U\rightarrow \R^{d'\times d}\), \(\gamma :[0,T]\times \R^{d'}\times U\times E\rightarrow \R^{d'}\),
\(f:[0,T]\times \R^{d'}\times \R\times \R^d\times \R\times U\rightarrow \R\), \(\ell:[0,T]\times E\rightarrow \R\) and \(h:\R^{d'}\rightarrow \R\) are given measurable functions.
\begin{assumption}\label{ass:standing}
\leavevmode\par\smallskip
\begin{enumerate}[label=\textup{(A\arabic*)},ref=\theassumption\textup{(A\arabic*)},wide=0pt]

    \item\label[assumption]{ass:A1} For any $(t,x,y,z,k,v,e)\in[0,T]\times \R^{d'}\times \R\times \R^d\times \R\times U\times E$, the functions \(b,\sigma\) are continuous with respect to
    \(t\), with \((x,v)\) fixed, and \(\gamma\) is continuous with respect to \(t\), with \((x,v,e)\) fixed; the generator \(f\) is jointly
    continuous with respect to \((t,y)\), with \((x,z,k,v)\) fixed;
    \(\ell\) is continuous with respect to \(t\), with \(e\) fixed.
        
    \item\label[assumption]{ass:A2} There exists a constant \(L_1>0\) such that for any \((t,x,x',v,v')\in [0,T]\times \R^{d'}\times \R^{d'}\times U\times U\),
    \begin{align*}
        |b(t,x,v)-b(t,x',v')|+|\sigma (t,x,v)-\sigma (t,x',v')|&\le L_1(|x-x'|+|v-v'|), \\
        |b(t,x,v)|+|\sigma (t,x,v)|&\le L_1(1+|x|).
    \end{align*}

    \item\label[assumption]{ass:A3} There exists a constant \(L_2>0\) such that for any \((t,x,x',v,v',e)\in [0,T]\times \R^{d'}\times \R^{d'}\times U\times U\times E\),
    \begin{align*}
        |\gamma (t,x,v,e)-\gamma (t,x',v',e)|&\le L_2(1 \land |e|)(|x-x'|+|v-v'|), \\
        |\gamma (t,x,v,e)|&\le L_2(1 \land |e|)(1+|x|), \\
        0\le \ell(t,e)&\le L_2(1 \land |e|).
    \end{align*}

    \item\label[assumption]{ass:A4} There exists a constant \(L_3>0\) such that for any \((t,x,x',y,z,k,v,v')\in [0,T]\times \R^{d'}\times \R^{d'}\times \R\times \R^d\times \R\times U\times U\),
    \begin{align*}
        |h(x)-h(x')| &\le L_3|x-x'|.
    \end{align*}

    \item\label[assumption]{ass:A5} There exists a constant \(\kappa \in \R\) such that for any \((t,x,y,y',z,k,v)\in [0,T]\times \R^{d'}\times \R\times \R\times \R^d\times \R\times U\),
    \begin{align*}
        (y-y')(f(t,x,y,z,k,v)-f(t,x,y',z,k,v))\le \kappa |y-y'|^2.
    \end{align*}
    Moreover, there exist constants \(L_4>0\) and \(p_0\ge 1\) such that
    \begin{align*}
        |f(t,0,y,0,0,v)-f(t,0,0,0,0,v)|\le L_4(1+|y|^{p_0}).
    \end{align*}

    \item\label[assumption]{ass:A6} For any \((t,x,y,z,v)\in [0,T]\times \R^{d'}\times \R\times \R^d\times U\), the $\R\to\R$ mapping \(k \rightarrow f(t,x,y,z,k,v)\) is nondecreasing.

    \item\label[assumption]{ass:A7} 
    There exist a continuous nondecreasing function
    \(\varpi (\cdot ):\R_+\to\R_+\) and a fixed constant \(A>0\), with
    \(\varpi(0)=0\) and \(\varpi(r)\le A(1+r)\), such that, for all
    \((t,x,x',y,z,z',k,k',v,v')\in[0,T]\times(\R^{d'})^2\times\R\times(\R^d)^2\times\R^2\times U^2\),
    we have
    \begin{align*}
    |f(t,x,y,z,k,v)-f(t,x',y,z',k',v')|\le
    \varpi(|x-x'|)+\varpi(|z-z'|)+\varpi(|k-k'|)+\varpi(|v-v'|).
    \end{align*}
\end{enumerate}
\end{assumption}

Set \(L_0:=|h(0)|+\sup_{(t,v)\in[0,T]\times U}|f(t,0,0,0,0,v)|\).
For simplicity, we write \(L:=\max\{L_0,L_1,L_2,L_3,L_4\}\). Throughout the paper,
\(C\) denotes a generic positive constant, which may vary from line to line and depends only on \((T,L,A,\kappa,p_0,\lambda)\).

Moreover, in \assitemref{ass:A5}, by an exponential change of variables
\[\bigl(\bar{Y}_s,\bar{Z}_s,\bar{K}_s(e)\bigr)=(\exp(\beta s)Y_s,\exp(\beta s)Z_s,\exp(\beta s)K_s(e)),\]
with \(\beta \ge \kappa\), one may reduce the monotonicity constant \(\kappa\) to a nonpositive one without loss of generality.
Accordingly, throughout \Cref{sec:preliminaries}, \Cref{sec:dpp} and \Cref{sec:hjb}, we assume \(\kappa = 0\), while in \Cref{sec:uniqueness} we impose a stronger condition \(\kappa < 0\).

The following lemma will be repeatedly used to handle the non-Lipschitz dependence of the generator \(f\) on the jump term.
The proof is trivial and is omitted.
\begin{lemma}\label{lem:monotone-jump-inequality}
    Let \(\nu:\R\rightarrow \R\) be a nondecreasing function such that
    \[\abs{\nu(k_1)-\nu(k_2)}\le \varpi(|k_1-k_2|),\ \forall k_1,k_2\in \R,\]
    where \(\varpi\) is the modulus of continuity introduced in \assitemref{ass:A7}.
    Then, for any \(\rho > 2A\), \(x \ge 0\) and \(y_1,y_2,k_1,k_2\in \R\), we have
    \(\varpi(x)\le \rho x + \varpi(2A/\rho)\) and
    \begin{align*}
        (y_1-y_2)(\nu(k_1)-\nu(k_2))\le \1_{\{(y_1-y_2)(k_1-k_2)>0\}}\rho(y_1-y_2)(k_1-k_2)
        +\varpi\Bigl(\frac{2A}{\rho}\Bigr)|y_1-y_2|.
    \end{align*}
\end{lemma}

We next introduce the inf-convolution of the generator \(f\) (see \cite[Lemma~1]{LepeltierSanMartin1997}).
\begin{lemma}\label{lem:lipschitz-approximation}

    For \(m>A\), define
    \[
    f_m( t,x,y,z,k,v ) := \inf_{(u_1,u_2)\in \R^d\times \R}[ f( t,x,y,u_1,u_2,v ) +m( |z-u_1|+|k-u_2| ) ].
    \]
If \(f\) satisfies \Cref{ass:A1,ass:A5,ass:A6,ass:A7}, then \(f_m\) satisfies the following properties:

    \begin{enumerate}[label=\textup{(B\arabic*)},ref=\textup{(B\arabic*)},wide=0pt]
    \item \(f_m\) satisfies \Cref{ass:A1,ass:A5,ass:A6,ass:A7} 
    with the same constants and modulus of continuity as \(f\), 
    except that the constant \(L_4\) in \assitemref{ass:A5} is replaced by \(L_4+2A\).
    Moreover, 
    for any
    \[(t,x,x',y,z,z',k,k',v,v')\in[0,T]\times(\R^{d'})^2\times\R\times(\R^d)^2\times\R^2\times U^2,\]
    \begin{align*}
        |f_m(t,x,y,z,k,v)-f_m(t,x',y,z',k',v')|\le\varpi(|x-x'|)+\varpi(|v-v'|)+m(|z-z'|+|k-k'|).
    \end{align*}

    \item 
    For any \((t,x,y,z,k,v)\in [0,T]\times \R^{d'}\times \R\times \R^d\times \R\times U,\)
    \begin{align*}
        |f_m(t,x,y,z,k,v)-f(t,x,y,z,k,v)|\le 2\varpi\Bigl(\frac{2A}{m-A}\Bigr).
    \end{align*}
    
    \end{enumerate}
\end{lemma}

\begin{lemma}[Composition with a linear-growth modulus]
\label{lem:modulus-uniform-integrability}
Let \(\omega:\R_+\to\R_+\) be continuous and nondecreasing, with
\(\omega(0)=0\) and
\(\omega(r)\le C_\omega(1+r)\). Let \(R\) be a nonnegative
progressively measurable process. Then, for any
\(\delta\in(0,1]\),
\begin{equation}\label{eq:modulus-truncation}
\begin{aligned}
\E^{\Fcal_t}\int_t^T\omega^2(R_s)\,\ds
\le
T\omega^2(\delta)
+
2C_\omega^2(1+\delta^{-2})
\E^{\Fcal_t}\int_t^T R_s^2\,\ds.
\end{aligned}
\end{equation}

Moreover, 
\[
\E^{\Fcal_t}\int_t^T\omega^2(R_s)\,\ds
\le \omega_{\ast}\biggl(\E^{\Fcal_t}\int_t^T R_s^2\,\ds\biggr),
\]
where
\(
\omega_{\ast}(r)
:=
\inf_{0<\delta\le1}
\bigl[
T\omega^2(\delta)
+2C_\omega^2(1+\delta^{-2})r
\bigr].
\)
In particular, \(\omega_{\ast}(r)\downarrow 0\) as \(r\downarrow0\), and
\(\omega_{\ast}\) has at most linear growth and is continuous in \(r\).
\end{lemma}

\begin{proof}
By the linear-growth condition, we have
\[
\omega^2(R_s)
\le
\omega^2(\delta)+\omega^2(R_s)\1_{\{R_s>\delta\}}
\le
\omega^2(\delta)+2C_\omega^2(1+R_s^2)\1_{\{R_s>\delta\}}
\le
\omega^2(\delta)+2C_\omega^2(1+\delta^{-2})R_s^2.
\]
Integrating and taking conditional expectations, we further have
\eqref{eq:modulus-truncation}. The remaining conclusions are immediate.
\end{proof}

Then we introduce a compensator inequality which will be used repeatedly without further comment.
\begin{remark}\label{rem:compensator-inequality}
Let \(H:\Omega\times[0,T]\times E\to[0,\infty)\) be a predictable process, and set
\[
    A_T:=\int_0^T\int_E H_s(e)\lambda(\de)\ds\ \ \ \text{and}\ \ \
    N_T:=\int_0^T\int_E H_s(e)\mu(\ds,\de).
\]
Then, for any \(q\ge 1\),
\begin{align*}
    \E\bigl[A_T^q\bigr]\le q^q\,\E\bigl[N_T^q\bigr].
\end{align*}
\end{remark}

    We first consider a general decoupled BSDE with jumps:
\begin{equation}\label{eq:bsde-general}
\begin{aligned}
        Y_t = \xi + \int_t^T{g\biggl(s,Y_s,Z_s,\int_E{\ell_s(e)K_s(e) \lambda(\de)}\biggr) \ds} - \int_t^T{\ip{Z_s}{\dW_s}} - \int_t^T{\int_E{K_s(e) \tilmu( \ds,\de )}},
\end{aligned}
    \end{equation}
where the generator \(g:\Omega \times [0,T]\times \R\times \R^d\times \R\rightarrow \R\)
is \(\operatorname{Prog}(\mathbb F)\otimes\mathcal B(\R\times\R^d\times\R)\text{-measurable}\), 
\(\xi:\Omega\rightarrow\R\) is \(\mathcal F_T\)-measurable, 
and \(\ell:[0,T]\times E\rightarrow \R\) is measurable.

We next introduce a set of assumptions needed for BSDE \eqref{eq:bsde-general}.
\begin{assumption}\label{ass:bsde-conditions}
\leavevmode\par\smallskip
\begin{enumerate}[label=\textup{(H\arabic*)},ref=\theassumption\textup{(H\arabic*)},wide=0pt]

    \item\label[assumption]{H1} There exists a constant \(L_2>0\) such that for any \((t,e)\in [0,T]\times E\),
    \begin{align*}
        0\le \ell(t,e)&\le L_2(1 \land |e|).
    \end{align*}

    \item\label[assumption]{H2} \(\E \bigl( |\xi |^p + \int_0^T{|g(t,0,0,0)|^p\dt} \bigr)<\infty\).

    \item\label[assumption]{H3} For any \((\omega ,t,z,k) \in \Omega \times [0,T]\times \R^d\times \R\),
    the \(\R\to\R\) mapping \(y\rightarrow g(\omega ,t,y,z,k)\) is continuous.
    Moreover, there exists a constant \(\kappa \in \R\) such that,
    for any \((\omega ,t,y,y',z,k)\in \Omega \times [0,T]\times \R\times \R\times \R^d\times \R\),
    \begin{align*}
        (y-y')(g(\omega ,t,y,z,k)-g(\omega ,t,y',z,k))\le \kappa |y-y'|^2.
    \end{align*}

    \item\label[assumption]{H4} There exist a continuous and nondecreasing function \(\varpi' (\cdot ) : \R_{+}\rightarrow \R_{+}\) and a fixed constant \(A>0\) with \(\varpi' (0)=0\) and \(\varpi' (x)\le A(1+x)\), such that for any \((\omega ,t,y,z,z',k,k')\in \Omega \times [0,T]\times \R\times \R^d\times \R^d\times \R\times \R\),
    \begin{align*}
        |g(\omega ,t,y,z,k)-g(\omega ,t,y,z',k')|\le \varpi' (|z-z'|) + \varpi' (|k-k'|).
    \end{align*}

    \item\label[assumption]{H5} For any \((\omega ,t,y,z)\in \Omega \times [0,T]\times \R\times \R^d\), the \(\R\to\R\) mapping \(k \rightarrow g(\omega ,t,y,z,k)\) is nondecreasing.

     \item\label[assumption]{H6}   For any \(r>0\), \(\psi_r \in L^1([0,T]\times \Omega ,\dt\otimes \Pp)\), where \[\psi_r(\omega,t):=\sup_{|y|\le r}|g(\omega ,t,y,0,0)-g(\omega ,t,0,0,0)|.\]
\end{enumerate}
\end{assumption}

\begin{lemma}\label{lem:bsde-existence-comparison}

    For a given \(p\ge2\), we say that a triple \((Y,Z,K)\in \mathcal{S}^p(0,T;\R)\times \mathcal{L}^p(0,T;\R^d)\times \mathcal{L}^p_{\lambda}(0,T;\R)\) is an \(L^p\)-solution of BSDE~\eqref{eq:bsde-general} with parameters \((\xi,g)\), 
    if it satisfies BSDE~\eqref{eq:bsde-general} \(\Pp\text{-a.s.}\) for any \(t\in [0,T]\). Then we have the following three conclusions:

    \begin{enumerate}[label=\textup{(\roman*)},ref=\textup{(\roman*)},wide=0pt]
    \item\label{lem:bsde-existence-comparison-i} Let \((\xi,g)\) and \((\xi',g')\) be two sets of parameters both satisfying \Cref{ass:bsde-conditions}, and \((Y,Z,K)\) and \((Y',Z',K')\) be the corresponding \(L^p\)-solutions of BSDE~\eqref{eq:bsde-general} with parameters \((\xi,g)\) and \((\xi',g')\), respectively.
    If \(\xi \le \xi',\ \Pp\text{-a.s.}\), and
    \[g\biggl(t,Y'_t,Z'_t,\int_E{\ell_t(e)K'_t(e) \lambda(\de)}\biggr)\le g'\biggl(t,Y'_t,Z'_t,\int_E{\ell_t(e)K'_t(e) \lambda(\de)}\biggr),\]
    \((\Pp\otimes\dt)\text{-a.e.}\), then \(Y_t\le Y'_t,\ \Pp\text{-a.s.}\) for any \(t\in [0,T]\).

    \item\label{lem:bsde-existence-comparison-ii} If the parameters \((\xi,g)\) satisfy \Cref{ass:bsde-conditions}, BSDE~\eqref{eq:bsde-general} has a unique \(L^p\)-solution.

    \item\label{lem:bsde-existence-comparison-iii} If the parameters \((\xi,g)\) satisfy \Cref{ass:bsde-conditions}, the following a priori estimate holds for the \(L^p\)-solution of BSDE~\eqref{eq:bsde-general}: there exists a constant \(C(p)>0\) such that
    \begin{align*}
        &\quad \E\Biggl[ \sup_{t\in [0,T]}|Y_t|^p + \biggl( \int_0^T{|Z_s|^2\ds} \biggr) ^{p/2} + \biggl( \int_0^T{\int_E{|K_s(e)|^2\lambda(\de) \ds} } \biggr) ^{p/2} \Biggr] \\
        &\le C(p)\biggl(1+\E|\xi |^p + \E\int_0^T{|g(s,0,0,0)|^p\ds} \biggr).
    \end{align*}
    \end{enumerate}
\end{lemma}
\begin{proof}
    The proof uses a similar approach to Ma, Fan and Song \cite{MaFanSong2013}, with the necessary modifications.
    By the standard exponential transformation, we may assume without loss of generality that \(\kappa \le 0\) in \assitemref{H3}.

    \begin{enumerate}[label=\textup{(\roman*)},ref=\textup{(\roman*)},wide=0pt]
    \item Set \(\hat{Y}=Y-Y',\ \hat{Z}=Z-Z',\ \hat{K}=K-K',\ \hat{\xi}=\xi -\xi'\) and
    \[\hat{g}_s = g\biggl(s,Y_s,Z_s,\int_E{\ell_s(e)K_s(e) \lambda(\de)}\biggr)
    - g'\biggl(s,Y'_s,Z'_s,\int_E{\ell_s(e)K'_s(e) \lambda(\de)}\biggr).\]
    Applying Tanaka's formula to \(\hat{Y}^+_s\) on \([t,T]\), we obtain
    \begin{equation}\label{eq:tanaka-comparison-estimate}
    \begin{aligned}
        \hat{Y}^+_t \le \int_t^T{\1_{\{\hat{Y}_s>0\}}\hat{g}_s\ds} - \int_t^T{\1_{\{\hat{Y}_{s-}>0\}}\ip[\big]{\hat{Z}_s}{\dW_s}}
        - \int_t^T{\int_E{\1_{\{\hat{Y}_{s-}>0\}}\hat{K}_s(e) \tilmu( \ds,\de )}}.
    \end{aligned}
    \end{equation}
    Using \Cref{H3,H4,H5} and \Cref{lem:monotone-jump-inequality}, we obtain the following estimate:
    \begin{equation}\label{eq:generator-comparison-estimate}
    \begin{aligned}
        &\quad \1_{\{\hat{Y}_s>0\}} \Biggl(g\biggl(s,Y_s,Z_s,\int_E{\ell_s(e)K_s(e) \lambda(\de)}\biggr) - g'\biggl(s,Y'_s,Z'_s,\int_E{\ell_s(e)K'_s(e) \lambda(\de)}\biggr)\Biggr) \\        
        &\le \1_{\{\hat{Y}_s>0\}}\rho|\hat{Z}_s| + \1_{\{\hat{Y}_s>0\}}\1_{\{\hat{Y}_s\int_E{\ell_s(e)\hat{K}_s(e) \lambda(\de)}>0\}}\rho \int_E{\ell_s(e)\hat{K}_s(e) \lambda(\de)} + 2\varpi'\Bigl(\frac{2A}{\rho}\Bigr).
    \end{aligned}
    \end{equation}
    Combining \eqref{eq:tanaka-comparison-estimate} with \eqref{eq:generator-comparison-estimate}, we have
    \begin{equation}\label{eq:comparison-girsanov-inequality}
    \begin{aligned}
        \hat{Y}^+_t 
        &\le 2T\varpi'\Bigl(\frac{2A}{\rho}\Bigr) - \int_t^T{\1_{\{\hat{Y}_{s-}>0\}}\ip[\bigg]{\hat{Z}_s}{\dW_s - \rho \1_{\{\hat{Z}_s \neq 0\}} \frac{\hat{Z}_s}{|\hat{Z}_s|}\ds}} \\
        &\quad - \int_t^T{\int_E{\1_{\{\hat{Y}_{s-}>0\}}\hat{K}_s(e)\bigl(\tilmu( \ds,\de ) - \1_{\{\hat{Y}_{s-}\int_E{\ell_s(e)\hat{K}_s(e) \lambda(\de)}>0\}}\rho \ell_s(e) \lambda(\de) \ds\bigr)}}.
    \end{aligned}
    \end{equation}

    Define \(\theta^{\rho,c}_s:=\rho \1_{\{\hat{Z}_s \neq 0\}} \hat{Z}_s/|\hat{Z}_s|\) and \(\theta^{\rho,\mu}_s(e):=\1_{\{\hat{Y}_{s-}\int_E{\ell_s(e)\hat{K}_s(e) \lambda(\de)}>0\}}\rho \ell_s(e)\).
    Here and below, every normalized-vector quotient is defined to be zero when its denominator is zero.
    It is immediate that \(|\theta^{\rho,c}_s|\le \rho\) and \(0\le \theta^{\rho,\mu}_s(e)\le \rho L_2(1\wedge|e|)\). Then set
    \[
    M^\rho_t := \int_0^t \ip[\big]{\theta^{\rho,c}_s}{\dW_s} +
    \int_0^t\int_E\theta^{\rho,\mu}_s(e)\,\tilmu(\ds,\de),
    \]
    and define \(\Theta^\rho_t:=\mathcal E(M^\rho)_t\),
    namely the Doléans-Dade exponential satisfying
    \(\dd\Theta^\rho_t = \Theta^\rho_{t-} \dd M^\rho_t\) and \(\Theta^\rho_0=1\).
    
    For any \(q>1\), set
    \[a_s^{(q)}:={}\frac{q(q-1)}{2}|\theta_s^{\rho,c}|^2
    +\int_E\bigl[(1+\theta_s^{\rho,\mu}(e))^q-q\theta_s^{\rho,\mu}(e)-1\bigr]\lambda(\de).\]
    It is immediate that
    \(0\le a_s^{(q)}\le C_{q,\rho},\)
    where \(C_{q,\rho}<\infty\) is a deterministic constant.
    On the other hand, by It\^o's formula,
    \[
    \dd \bigl[\exp(-C_{q,\rho}t)(\Theta_t^\rho)^q\bigr]
    =\exp(-C_{q,\rho}t)(\Theta_{t-}^\rho)^q\biggl[q\,\ip{\theta_t^{\rho,c}}{\dW_t}
    +\int_E\bigl[(1+\theta_t^{\rho,\mu}(e))^q-1\bigr]\tilmu(\dt,\de)
    +(a_t^{(q)}-C_{q,\rho})\dt\biggr].
    \]
    Consequently,
    \(\exp(-C_{q,\rho}t)(\Theta_t^\rho)^q\)
    is a nonnegative local supermartingale, and hence a supermartingale. 
    Therefore, for \(0\le t\le s\le T\),
    \[\E\bigl[(\Theta_s^\rho)^q\mid\Fcal_t\bigr]\le \exp(C_{q,\rho}(s-t))(\Theta_t^\rho)^q.\]
    In particular,
    \(\sup_{0\le t\le T}\E(\Theta_t^\rho)^q\le \exp(C_{q,\rho}T),\)
    which implies that \(\{\Theta_t^\rho\}_{0\le t\le T}\) are uniformly integrable.
    On the other hand, since \(\Theta^\rho\) is a local martingale, 
    let \(\{\tau_n\}_{n\ge 1}\) be its localizing sequence. By optional sampling,
    \(\sup_n \E(\Theta_{t\wedge\tau_n}^\rho)^q \le \exp(C_{q,\rho}T),\) for \(0\le t\le T\).
    Thus \(\{(\Theta_{t\wedge\tau_n}^\rho)^q\}_{n\ge 1}\) are uniformly
    integrable for \(0\le t\le T\) and \(\Theta^\rho\) is a uniformly integrable martingale.

    Define a probability measure \(\Pp^\rho\) on \((\Omega,\Fcal)\) by
    \(\dP^\rho=\Theta^\rho_T\dP\).
    By Girsanov's theorem, \(W^\rho_t:=W_t-\int_0^t\theta^{\rho,c}_s\ds\)
    is a \((\mathbb F,\Pp^\rho)\)-Brownian motion, and the \(\Pp^\rho\)-predictable compensator of \(\mu\) is
    \( \nu^\rho(\ds,\de) = \bigl(1+\theta^{\rho,\mu}_s(e)\bigr)\lambda(\de)\ds\).
    Consequently, \(\tilmu^\rho(\ds,\de) := \tilmu(\ds,\de)-\theta^{\rho,\mu}_s(e)\lambda(\de)\ds\) is a \((\mathbb F,\Pp^\rho)\)-random martingale measure.

    We next verify that the stochastic integrals in
    \eqref{eq:comparison-girsanov-inequality} are martingales under \(\Pp^\rho\).
    For the Brownian part, since \(\Theta^\rho_T\in L^2(\Pp)\), by Cauchy's
    inequality,
    \[
    \E^{\Pp^\rho}\Biggl[\biggl(\int_0^T \1_{\{\hat Y_{s-}>0\}}|\hat Z_s|^2\ds\biggr)^{1/2}\Biggr]
    \le \bigl(\E\bigl[\bigl(\Theta^\rho_T\bigr)^2\bigr]\bigr)^{1/2}
    \biggl(\E\int_0^T |\hat Z_s|^2\ds \biggr)^{1/2} <\infty.
    \]
Similarly, for the jump part,
    \begin{align*}
    \E^{\Pp^\rho}\Biggl[ \biggl( \int_0^T\int_E  \1_{\{\hat Y_{s-}>0\}}|\hat K_s(e)|^2\mu(\ds,\de) \biggr)^{1/2} \Biggr]
    \le \bigl(\E\bigl[\bigl(\Theta^\rho_T\bigr)^2\bigr]\bigr)^{1/2} \biggl( \E\int_0^T\int_E |\hat K_s(e)|^2\lambda(\de)\ds \biggr)^{1/2} < \infty.
    \end{align*}
Hence \(\int_0^\cdot \1_{\{\hat Y_{s-}>0\}}\ip[\big]{\hat Z_s}{\dW^\rho_s}\) and
    \(\int_0^\cdot\int_E \1_{\{\hat Y_{s-}>0\}}\hat K_s(e)\,\tilmu^\rho(\ds,\de)\)
    are both true \(\Pp^\rho\)-martingales.
    Taking the conditional expectation with respect to \(\Fcal_t\) under
    \(\Pp^\rho\) in \eqref{eq:comparison-girsanov-inequality}, we obtain
    \begin{align*}
        \hat{Y}^+_t \le C\varpi'\Bigl(\frac{2A}{\rho}\Bigr).
    \end{align*}

    Letting \(\rho \rightarrow \infty\), we obtain \(\hat{Y}^+_t = 0\) for any \(t\in [0,T]\),
    which implies \(Y_t \le Y'_t\), \(\Pp\text{-a.s.}\) for any \(t\in [0,T]\).

    \item It remains to prove existence, since uniqueness follows immediately from \ref{lem:bsde-existence-comparison-i}.

    For a positive integer \(n>A\), define \(g_n\) by the approximation used in \Cref{lem:lipschitz-approximation}. 
    Then \(g_n\) satisfies \Cref{H2,H3,H4,H5,H6} with the same constants as \(g\) 
    except the local growth bound. \(g_n\) is also Lipschitz continuous in \((z,k)\).
    From \cite[Theorems~1 and 2]{KrusePopier2016}, we can deduce that BSDE~\eqref{eq:bsde-general} with parameters \((\xi,g_n)\) has a unique \(L^p\)-solution \((Y^n,Z^n,K^n)\).

    Our goal is to show that the sequence \((Y^n,Z^n,K^n)\) is a Cauchy sequence in \(\mathcal{S}^p\times \mathcal{L}^p\times \mathcal{L}^p_{\lambda}\), whose limit is the desired solution of BSDE~\eqref{eq:bsde-general} with parameters \((\xi,g)\).
    Set \(\hat{Y}^{m,n}=Y^m-Y^n,\ \hat{Z}^{m,n}=Z^m-Z^n,\ \hat{K}^{m,n}=K^m-K^n\), and 
    \[
    \begin{aligned}
    \hat g_s^{m,n}:={}g_m\biggl(s,Y_s^m,Z_s^m,\int_E \ell_s(e)K_s^m(e)\lambda(\de)\biggr)
    -g_n\biggl(s,Y_s^n,Z_s^n,\int_E \ell_s(e)K_s^n(e)\lambda(\de)\biggr).
    \end{aligned}
    \]
    For \(m\ge n>A\), we have

    \begin{equation}\label{eq:cauchy-decomposition}
    \begin{aligned}
        &\quad \sgn\bigl(\hat{Y}^{m,n}_s\bigr)\Biggl( g_m\biggl(s,Y^m_s,Z^m_s,\int_E{\ell_s(e)K^m_s(e) \lambda(\de)}\biggr) - g_n\biggl(s,Y^n_s,Z^n_s,\int_E{\ell_s(e)K^n_s(e) \lambda(\de)}\biggr) \Biggr) \\
        &\le \1_{\{\hat{Y}^{m,n}_s\neq0\}}\rho|\hat{Z}^{m,n}_s| + \sgn\bigl(\hat{Y}^{m,n}_s\bigr)\1_{\{\hat{Y}^{m,n}_s\int_E{\ell_s(e)\hat{K}^{m,n}_s(e) \lambda(\de)}>0\}}\rho \int_E{\ell_s(e)\hat{K}^{m,n}_s(e) \lambda(\de)} \\
        &\quad + 2\varpi'\Bigl(\frac{2A}{\rho}\Bigr) +4\varpi'\Bigl(\frac{2A}{n-A}\Bigr),
    \end{aligned}
    \end{equation}
    Applying Tanaka's formula to \(|\hat{Y}^{m,n}_s|\) on \([t,T]\), we have
    \begin{equation}\label{eq:cauchy-sign-estimate}
    \begin{aligned}
        |\hat{Y}^{m,n}_t| \le \int_t^T{\sgn\bigl(\hat{Y}^{m,n}_s\bigr)\hat{g}^{m,n}_s\ds} - \int_t^T{\sgn\bigl(\hat{Y}^{m,n}_{s-}\bigr)\ip[\big]{\hat{Z}^{m,n}_s}{\dW_s}} 
        -\int_t^T{\int_E{\sgn\bigl(\hat{Y}^{m,n}_{s-}\bigr)\hat{K}^{m,n}_s(e) \tilmu( \ds,\de )}}.
    \end{aligned}
    \end{equation}
 Combining \eqref{eq:cauchy-decomposition} with \eqref{eq:cauchy-sign-estimate}, we have
    \begin{equation}\label{eq:cauchy-girsanov-inequality}
    \begin{aligned}
        |\hat{Y}^{m,n}_t| 
        &\le T\Bigl(2\varpi'\Bigl(\frac{2A}{\rho}\Bigr)+4\varpi'\Bigl(\frac{2A}{n-A}\Bigr)\Bigr)\\
        &\quad - \int_t^T{\sgn\bigl(\hat{Y}^{m,n}_{s-}\bigr)\ip[\bigg]{\hat{Z}^{m,n}_s}{\dW_s - \rho\sgn\bigl(\hat{Y}^{m,n}_{s-}\bigr)\1_{\{\hat{Z}^{m,n}_s\neq0\}}\frac{\hat{Z}^{m,n}_s}{|\hat{Z}^{m,n}_s|}\ds}} \\
        &\quad - \int_t^T{\int_E{\sgn\bigl(\hat{Y}^{m,n}_{s-}\bigr)\hat{K}^{m,n}_s(e)\bigl(\tilmu( \ds,\de ) - \1_{\{\hat{Y}^{m,n}_{s-}\int_E{\ell_s(e)\hat{K}^{m,n}_s(e) \lambda(\de)}>0\}}\rho \ell_s(e) \lambda(\de) \ds\bigr)}}.
    \end{aligned}
    \end{equation}
    Using the same Girsanov transformation argument as in \ref{lem:bsde-existence-comparison-i} and taking conditional expectations on both sides of \eqref{eq:cauchy-girsanov-inequality} under the new probability measure, we have
    \begin{align}\label{eq:cauchy-y-estimate}
        |\hat{Y}^{m,n}_t| \le C\Bigl(\varpi'\Bigl(\frac{2A}{\rho}\Bigr)+\varpi'\Bigl(\frac{2A}{n-A}\Bigr)\Bigr).
    \end{align}
    Since \eqref{eq:cauchy-y-estimate} holds uniformly with respect to \(t\), by taking the supremum over \(t\in [0,T]\),
    we can deduce that \(\{Y^n\}_{n>A}\) is a Cauchy sequence in \(\mathcal{S}^p\).

    We then turn to proving that \(\{Z^n\}_{n>A}\) and \(\{K^n\}_{n>A}\) are Cauchy sequences in \(\mathcal{L}^p\) and \(\mathcal{L}^p_{\lambda}\), respectively.
    First, we establish the following estimate:
    \begin{equation}\label{eq:cauchy-zk-ito-estimate}
    \begin{aligned}
        &\quad \hat{Y}^{m,n}_s \Biggl(g_m\biggl(s,Y^m_s,Z^m_s,\int_E{\ell_s(e)K^m_s(e) \lambda(\de)}\biggr) - g_n\biggl(s,Y^n_s,Z^n_s,\int_E{\ell_s(e)K^n_s(e) \lambda(\de)}\biggr)\Biggr) \\
        &\le |\hat{Y}^{m,n}_s|\biggl(\rho|\hat{Z}^{m,n}_s|+\rho\int_E{\ell_s(e)|\hat{K}^{m,n}_s(e)| \lambda(\de)} + 2\varpi'\Bigl(\frac{2A}{\rho}\Bigr) +4\varpi'\Bigl(\frac{2A}{n-A}\Bigr)\biggr) \\
        &\le C(\rho,p)|\hat{Y}^{m,n}_s|^2 + \frac{1}{6}|\hat{Z}^{m,n}_s|^2 + \frac{1}{C(p)}\int_E{|\hat{K}^{m,n}_s(e)|^2 \lambda(\de)} + \Bigl(\varpi' \Bigl(\frac{2A}{\rho}\Bigr)\Bigr)^2+\Bigl(\varpi'\Bigl(\frac{2A}{n-A}\Bigr)\Bigr)^2.
    \end{aligned}
    \end{equation}
Applying It\^o's formula to \(|\hat{Y}^{m,n}_t|^2\), together with the above estimate \eqref{eq:cauchy-zk-ito-estimate}, we have
    \begin{equation*}
    \begin{aligned}
        &\quad \int_0^T{|\hat{Z}^{m,n}_s|^2\ds} + \int_0^T{\int_E{|\hat{K}^{m,n}_s(e)|^2 \mu( \ds,\de )}} \\
        &\le \int_0^T{2\hat{Y}^{m,n}_s\hat{g}^{m,n}_s\ds} - \int_0^T{2\hat{Y}^{m,n}_{s-}\ip[\big]{\hat{Z}^{m,n}_s}{\dW_s}} - \int_0^T{\int_E{2\hat{Y}^{m,n}_{s-}\hat{K}^{m,n}_s(e) \tilmu( \ds,\de )}} \\
        &\le \int_0^T\left[
        C(\rho,p)|\hat{Y}^{m,n}_s|^2 + \frac{1}{6}|\hat{Z}^{m,n}_s|^2 + \frac{1}{C(p)}\int_E{|\hat{K}^{m,n}_s(e)|^2 \lambda(\de)}
        + \Bigl(\varpi' \Bigl(\frac{2A}{\rho}\Bigr)\Bigr)^2+\Bigl(\varpi'\Bigl(\frac{2A}{n-A}\Bigr)\Bigr)^2
        \right]\ds \\
        &\quad + \biggl|\int_0^T{2\hat{Y}^{m,n}_{s-}\ip[\big]{\hat{Z}^{m,n}_s}{\dW_s}}\biggr| + \biggl| \int_0^T{\int_E{2\hat{Y}^{m,n}_{s-}\hat{K}^{m,n}_s(e) \tilmu( \ds,\de )}} \biggr|.
    \end{aligned}
    \end{equation*}
 Thus
    \begin{equation}\label{eq:apriori-ito-formula}
    \begin{aligned}
        &\quad \E\Biggl[ \biggl(\int_0^T{|\hat{Z}^{m,n}_s|^2\ds}\biggr)^{p/2}\Biggr] + \E\Biggl[ \biggl(\int_0^T{\int_E{|\hat{K}^{m,n}_s(e)|^2 \mu( \ds,\de )}}\biggr)^{p/2} \Biggr] \\
        &\le C(\rho,p)\E\sup_{0\le t\le T}{|\hat{Y}^{m,n}_t|^p} + C(p)\Bigl( \varpi'\Bigl(\frac{2A}{\rho}\Bigr)+\varpi'\Bigl(\frac{2A}{n-A}\Bigr) \Bigr)^p \\
        &\quad + C(p)\E\biggl|\int_0^T {\hat{Y}^{m,n}_{s-}\ip[\big]{\hat{Z}^{m,n}_s}{\dW_s}}\biggr|^{p/2} + C(p)\E\biggl| \int_0^T{\int_E{\hat{Y}^{m,n}_{s-}\hat{K}^{m,n}_s(e) \tilmu( \ds,\de )}} \biggr|^{p/2}.
    \end{aligned}
    \end{equation}
By the Burkholder--Davis--Gundy (BDG) inequality and Young's inequality, we have
    \begin{equation}\label{eq:cauchy-bdg-estimate-brownian}
    \begin{aligned}
        &\quad C(p)\E\biggl|\int_0^T {\hat{Y}^{m,n}_{s-}\ip[\big]{\hat{Z}^{m,n}_s}{\dW_s}}\biggr|^{p/2}+C(p)\E\biggl| \int_0^T{\int_E{\hat{Y}^{m,n}_{s-}\hat{K}^{m,n}_s(e) \tilmu( \ds,\de )}} \biggr|^{p/2} \\
        &\le C(p)\E\sup_{0\le t\le T}{|\hat{Y}^{m,n}_t|^p}+\frac{1}{2}\E\Biggl[ \biggl(\int_0^T{|\hat{Z}^{m,n}_s|^2\ds}\biggr)^{p/2} \Biggr]+ \frac{1}{2}\E\Biggl[ \biggl(\int_0^T{\int_E{|\hat{K}^{m,n}_s(e)|^2 \mu( \ds,\de )}}\biggr)^{p/2} \Biggr].
    \end{aligned}
    \end{equation}
Combining \eqref{eq:apriori-ito-formula} with \eqref{eq:cauchy-bdg-estimate-brownian}, we have
    \begin{equation}\label{eq:cauchy-zk-convergence}
    \begin{aligned}
        &\quad \E\Biggl[ \biggl(\int_0^T{|\hat{Z}^{m,n}_s|^2\ds}\biggr)^{p/2}\Biggr] + \E\Biggl[ \biggl(\int_0^T{\int_E{|\hat{K}^{m,n}_s(e)|^2 \mu( \ds,\de )}}\biggr)^{p/2} \Biggr] \\
        &\le C(\rho,p)\E\sup_{0\le t\le T}{|\hat{Y}^{m,n}_t|^p} + C(p)\Bigl( \varpi'\Bigl(\frac{2A}{\rho}\Bigr)+\varpi'\Bigl(\frac{2A}{n-A}\Bigr) \Bigr)^p.
    \end{aligned}
    \end{equation}
Take \(m,n\rightarrow \infty\), and then \(\rho \rightarrow \infty\), which together with \Cref{rem:compensator-inequality} implies that \(\{Z^n\}_{n>A}\) and \(\{K^n\}_{n>A}\) are Cauchy sequences in \(\mathcal{L}^p\) and \(\mathcal{L}^p_{\lambda}\), respectively.

    Now, since \(\{Z^n\}_{n>A}\) and \(\{K^n\}_{n>A}\) are Cauchy sequences in \(\mathcal{L}^p\) and \(\mathcal{L}^p_{\lambda}\), respectively, passing to a subsequence if necessary, we can assume that the convergence also holds \((\dt\otimes\Pp)\text{-a.e.}\)
    and \(\sup_{n>A}{|Z^n_t|} + \sup_{n>A}{\int_E{\ell_t(e)|K^n_t(e)| \lambda(\de)}} \in L^1([0,T]\times\Omega,\dt\otimes\Pp)\). 
Therefore, by the BDG inequality and the dominated convergence theorem, we obtain
    \begin{equation}\label{eq:bsde-uniqueness-limit-convergence}
    \begin{aligned}
        &\E\sup_{0\le t\le T}\biggl|\int_t^T{\ip[\big]{\bigl(Z^n_s-Z_s\bigr)}{\dW_s}}\biggr|+\E\sup_{0\le t\le T}\biggl|\int_t^T{\int_E{\bigl(K^n_s(e)-K_s(e)\bigr) \tilmu( \ds,\de )}}\biggr| \rightarrow 0 \\
        &\E\int_0^T{\biggl|g_n\biggl(t,Y^n_t,Z^n_t,\int_E{\ell_t(e)K^n_t(e) \lambda(\de)}\biggr) - g\biggl(t,Y^n_t,Z_t,\int_E{\ell_t(e)K_t(e) \lambda(\de)}\biggr)\biggr| \dt} \rightarrow 0 
    \end{aligned}
    \end{equation}

Since \(g_n\) is increasing as \(n\to\infty\),
\ref{lem:bsde-existence-comparison-i} implies that
\(Y_t^n\) is nondecreasing in \(n\). Together with the convergence of
\(Y^n\) in \(\mathcal{S}^p\), we have \(Y^n_t \uparrow Y_t\), for any \(t\in [0,T]\), \(\Pp\text{-a.s.}\)

Set \(\bar K_t:=\int_E \ell_t(e)K_t(e)\lambda(\de)\).
We first note that
\[
    \int_0^T
    |g(t,Y_t,Z_t,\bar K_t)|\dt
    <\infty,
    \qquad \Pp\text{-a.s.}
\]
Indeed, since \(Y\in\mathcal{S}^p\), its path is c\`adl\`ag and bounded
on \([0,T]\), \(\Pp\)-a.s. For \(R>0\), on the event
\(B_R:=\{\sup_{0\le t\le T}|Y_t|\le R\}\)
we have, by \assitemref{H4} and \assitemref{H6},
\[
\begin{aligned}
    |g(t,Y_t,Z_t,\bar K_t)|
    &\le
    |g(t,0,0,0)|
    +\psi_R(t)
    +\varpi'(|Z_t|)
    +\varpi'(|\bar K_t|).
\end{aligned}
\]

Hence the right-hand side is integrable on
\([0,T]\times\Omega\), locally on \(B_R\). Since
\(B_R\uparrow\Omega\), \(\Pp\)-a.s., the claimed pointwise integrability
follows.

We now prove the convergence in the \(y\)-variable. 
Set \(n_0:=\lfloor 2A+1 \rfloor\). Define
\[
    D_n:=\int_0^T
    |g(t,Y_t^n,Z_t,\bar K_t)-g(t,Y_t,Z_t,\bar K_t)|\dt .
\]
For \(R>0\), let
\(A_R:=\bigl\{\sup_{0\le t\le T}|Y_t^{n_0}|\vee\sup_{0\le t\le T}|Y_t|\le R\bigr\}.\)
Since \(Y^{n_0}\le Y^n\le Y\), on \(A_R\) we have
\[
    |Y_t^n|\le R,
    \qquad
    |Y_t|\le R,
    \qquad
    0\le t\le T,
\]
for every \(n\ge n_0\). Therefore, by \assitemref{H4} and
\assitemref{H6},
\[
\begin{aligned}
    \1_{A_R}|g(t,Y_t^n,Z_t,\bar K_t)-g(t,Y_t,Z_t,\bar K_t)|
    \le 2\1_{A_R}\bigl[\psi_R(t)+\varpi'(|Z_t|)+\varpi'(|\bar K_t|)\bigr].
\end{aligned}
\]
The right-hand side belongs to
\(L^1([0,T]\times\Omega,\dt\otimes\Pp)\). 
Hence, by the dominated convergence theorem, for every fixed \(R>0\),
we have that \(\E[\1_{A_R}D_n]\longrightarrow 0.\)
Consequently, for every \(\varepsilon>0\),
\[
\begin{aligned}
    \Pp(D_n>\varepsilon)
    \le\Pp\bigl(A_R^c\bigr)+\Pp(\1_{A_R}D_n>\varepsilon)
    \le\Pp\bigl(A_R^c\bigr)+\frac{1}{\varepsilon}\E[\1_{A_R}D_n].
\end{aligned}
\]
Letting first \(n\to\infty\) and then \(R\to\infty\), and using
\(Y^{n_0},Y\in\mathcal{S}^p\), we obtain
\begin{equation}\label{eq:cauchy-limit-convergence}
    \int_0^T|g(t,Y_t^n,Z_t,\bar K_t)-g(t,Y_t,Z_t,\bar K_t)|\dt
    \longrightarrow 0 \quad\text{in probability}.
\end{equation}

Combining \eqref{eq:bsde-uniqueness-limit-convergence} with
\eqref{eq:cauchy-limit-convergence}, we have
\[
\begin{aligned}
    \int_0^T
    \biggl|g_n\biggl(t,Y_t^n,Z_t^n,\int_E \ell_t(e)K_t^n(e)\lambda(\de)\biggr) 
    -g\biggl(t,Y_t,Z_t,\int_E \ell_t(e)K_t(e)\lambda(\de)\biggr)\biggr|\dt
    \longrightarrow 0 \quad\text{in probability}.
\end{aligned}
\]
Together with the convergence of \(Y^n\), \(Z^n\), \(K^n\) and
\eqref{eq:bsde-uniqueness-limit-convergence}, we may extract a
subsequence along which all the above convergences hold
\(\Pp\)-a.s. Passing to the limit in the BSDE for this subsequence,
uniformly in \(t\in[0,T]\), yields
that \((Y,Z,K)\) is indeed an \(L^p\)-solution of \eqref{eq:bsde-general} with parameters \((\xi,g)\).

    \item Since \(g_n\) is Lipschitz continuous in \((z,k)\) with Lipschitz constant \(n\),
    by \cite[Theorems~1 and 2]{KrusePopier2016}, we have the following estimate for \(n\ge 2A\):
    \begin{equation}\label{eq:apriori-lipschitz-estimate}
    \begin{aligned}
        &\quad \E \sup_{0\le t\le T}{|Y^n_t|^p} + \E\biggl(\int_0^T{|Z^n_s|^2\ds}\biggr)^{p/2} + \E\biggl(\int_0^T{\int_E{|K^n_s(e)|^2 \lambda(\de)\ds}}\biggr)^{p/2} \\
        &\le C(p,n)\biggl(1 + \E |\xi|^p + \E\int_0^T{|g(s,0,0,0)|^p\ds} \biggr).
    \end{aligned}
    \end{equation}
    On the other hand, adding \eqref{eq:cauchy-y-estimate} and \eqref{eq:cauchy-zk-convergence}, we obtain
    \begin{equation}\label{eq:apriori-cauchy-estimate}
    \begin{aligned}
        &\quad \E \sup_{0\le t\le T}{|\hat{Y}^{m,n}_t|^p} + \E\biggl(\int_0^T{|\hat{Z}^{m,n}_s|^2\ds}\biggr)^{p/2} + \E\biggl(\int_0^T{\int_E{|\hat{K}^{m,n}_s(e)|^2 \lambda(\de)\ds}}\biggr)^{p/2} \\
        &\le C(\rho,p)\Bigl( \varpi'\Bigl(\frac{2A}{\rho}\Bigr)+\varpi'\Bigl(\frac{2A}{n-A}\Bigr) \Bigr)^p.
    \end{aligned}
    \end{equation}
Fix \(\rho=2A+1\) and \(n=\lfloor 2A+1 \rfloor\). Taking \(m\rightarrow \infty\) in \eqref{eq:apriori-cauchy-estimate},
    we obtain an estimate for the difference between \((Y,Z,K)\) and \((Y^{n},Z^{n},K^{n})\).
   A combination of this estimate with \eqref{eq:apriori-lipschitz-estimate} yields
    \begin{equation*}
    \begin{aligned}
        &\quad \E \sup_{0\le t\le T}{|Y_t|^p} + \E\biggl(\int_0^T{|Z_s|^2\ds}\biggr)^{p/2} + \E\biggl(\int_0^T{\int_E{|K_s(e)|^2 \lambda(\de)\ds}}\biggr)^{p/2} \\
        &\le C(p)\biggl(1 + \E |\xi|^p + \E\int_0^T{|g(s,0,0,0)|^p\ds} \biggr).
    \end{aligned}
    \end{equation*}
    \end{enumerate}
Thus the desired a priori estimate for the solution of BSDE~\eqref{eq:bsde-general} is obtained.

\end{proof}

Next we consider FBSDE \eqref{eq:controlled-fbsde}. Suppose that \Cref{ass:standing} holds. Then, for any \(p \geq 2\), \(\eta \in L_{\Fcal_t}^{p}\) and \(v\in \Ucal_{t,T}\), the mappings \(\xi:=h\bigl(X^{t,\eta,v}_T\bigr)\) and
\(g(\omega,s,y,z,k):=f\bigl(s,X^{t,\eta,v}_s,y,z,k,v_s\bigr)\) satisfy \Cref{ass:bsde-conditions}.
It is immediate from \Cref{lem:bsde-existence-comparison} that 
FBSDE \eqref{eq:controlled-fbsde} has a unique \(L^p\)-solution.

We next deduce some useful estimates.
\begin{theorem}\label{thm:bsde-continuous-dependence}
    Suppose that \Cref{ass:standing} holds.
    Let \({\{(\eta^n,v^n)\}_{n \in \mathbb{N}} \subset L_{\Fcal_t}^{p} \times \Ucal_{t,T}}\) with \(p > 2\) be such that
    \[\eta ^n\rightarrow \eta \ \text{in\ }L_{\Fcal_t}^{2}(\Omega;\R^{d'})\ \ \ \text{and}\ \ \ v^n\rightarrow v\ \text{in\ }\mathcal{L}^2,\]
and \(( X^{t,\eta ^n,v^n},Y^{t,\eta ^n,v^n},Z^{t,\eta ^n,v^n},K^{t,\eta ^n,v^n}( e ) ) \) and
    \(( X^{t,\eta ,v},Y^{t,\eta ,v},Z^{t,\eta ,v},K^{t,\eta ,v}( e ) ) \)
   be the corresponding solutions of \eqref{eq:controlled-fbsde} with parameters \(( \eta ^n,v^n ) \text{ and }( \eta ,v )\), respectively.
Then we have the following continuous dependence result:
    \begin{align*}
        \E\Bigl[ \sup_{t\le s\le T}| Y_{s}^{t,\eta ^n,v^n}-Y_{s}^{t,\eta ,v} |^2 \Bigr] \rightarrow 0,\quad\text{as }n\rightarrow \infty.
    \end{align*}

\end{theorem}
\begin{proof}
    Let \(f_m\) be the approximation sequence introduced in \Cref{lem:lipschitz-approximation}.
    Denote by \(\bigl(Y^{t,\eta^n,v^n,m},Z^{t,\eta^n,v^n,m},K^{t,\eta^n,v^n,m}(e)\bigr)\) and \((Y^{t,\eta,v,m},Z^{t,\eta,v,m},K^{t,\eta,v,m}(e))\)
    the solutions of \eqref{eq:controlled-fbsde} with parameters \((\eta^n,v^n,f_m)\) and \((\eta,v,f_m)\), respectively.

    In order to derive the estimate of \(| Y_{s}^{t,\eta ^n,v^n}-Y_{s}^{t,\eta ,v} |^2\), we can decompose it as follows:
    \begin{equation*}
    \begin{aligned}
	    | Y_{s}^{t,\eta ^n,v^n}-Y_{s}^{t,\eta ,v} |^2
	    \le C\bigl( | Y_{s}^{t,\eta ^n,v^n}-Y_{s}^{t,\eta ^n,v^n,m} |^2+| Y_{s}^{t,\eta ^n,v^n,m}-Y_{s}^{t,\eta ,v,m} |^2+| Y_{s}^{t,\eta ,v,m}-Y_{s}^{t,\eta ,v} |^2 \bigr).
    \end{aligned}
    \end{equation*}
 Set \(\bar{Y}^m_s=Y^{t,\eta,v}_s-Y^{t,\eta,v,m}_s,\ \bar{Z}^m_s=Z^{t,\eta,v}_s-Z^{t,\eta,v,m}_s,\ \bar{K}^m_s(e)=K^{t,\eta,v}_s(e)-K^{t,\eta,v,m}_s(e)\).    For the first and the third terms, we argue as in the proof of \Cref{lem:bsde-existence-comparison}\,\ref{lem:bsde-existence-comparison-ii} to show that they converge to 0 as \(m\rightarrow \infty\) uniformly with respect to \((\eta,v)\).
    That is to say for any \(\eta \in L_{\Fcal_t}^{2},\ v\in \Ucal_{t,T},\)
    \[\sup_{t\le s\le T}|\bar{Y}^m_s |^2\le C\varpi^2 \Bigl( \frac{2A}{m-A} \Bigr).\]
Applying It\^o's formula to \(\exp(\alpha s){|\bar{Y}_{s}^{m}|}^{2}\) on \([s,T]\), we obtain
    \begin{align*}
        &\quad \exp(\alpha s)|\bar{Y}_{s}^{m}|^2+\int_s^T{\alpha \exp(\alpha r)|\bar{Y}_r^{m}|^2\dr}
        +\int_s^T{2\exp(\alpha r)\bar{Y}_{r-}^{m}\ip[\big]{\bar{Z}_r^{m}}{\dW_r}}+\int_s^T{\int_E{2\exp(\alpha r)\bar{Y}_{r-}^{m} \bar{K}_r^{m}( e ) \tilmu( \dr,\de )}}\\*
        &\le \int_s^T{2\exp(\alpha r)\bar{Y}_r^{m}\bigl( f\bigl( r,X_r^{t,\eta ,v},Y_r^{t,\eta ,v},Z_r^{t,\eta ,v},\widehat K_{r}^{t,\eta ,v},v_r \bigr) 
        - f_m\bigl( r,X_r^{t,\eta ,v},Y_r^{t,\eta ,v ,m},Z_r^{t,\eta ,v, m},\widehat K_{r}^{t,\eta ,v,m},v_r \bigr) \bigr) \dr},
    \end{align*}
    where \(\widehat K_{r}^{t,\eta ,v}:=\int_E \ell_r(e)K_{r}^{t,\eta ,v}(e)\lambda(\de)\) and
    \(\widehat K_{r}^{t,\eta ,v,m}:=\int_E \ell_r(e)K_{r}^{t,\eta ,v,m}(e)\lambda(\de)\).
    Moreover, with \Cref{lem:monotone-jump-inequality} we have
    \begin{align*}
        &\quad \bar{Y}_{r}^{m}\bigl(
        f\bigl( r,X_{r}^{t,\eta ,v},Y_{r}^{t,\eta ,v},Z_{r}^{t,\eta ,v},
        \widehat K_{r}^{t,\eta ,v},v_r \bigr)
        -f_m\bigl( r,X_{r}^{t,\eta ,v},Y_{r}^{t,\eta ,v,m},Z_{r}^{t,\eta ,v,m},
        \widehat K_{r}^{t,\eta ,v,m},v_r \bigr)
        \bigr)\\
        &\le\rho |\bar{Y}_{r}^{m}||\bar{Z}_{r}^{m}|+
        \1_{\{ \bar{Y}_{r}^{m}\int_E{\ell_r( e ) \bar{K}_{r}^{m}( e ) \lambda(\de)}>0 \}}\rho \bar{Y}_{r}^{m}\int_E{\ell_r( e ) \bar{K}_{r}^{m}( e ) \lambda(\de)}
        +\Bigl( 2\varpi \Bigl( \frac{2A}{\rho} \Bigr) +2\varpi \Bigl( \frac{2A}{m-A} \Bigr) \Bigr) |\bar{Y}_{r}^{m}|.
    \end{align*}
    Thus, we can deduce that
    \begin{equation}\label{eq:continuous-dependence-approx-girsanov}
    \begin{aligned}
	    &\quad \exp(\alpha s)|\bar{Y}_{s}^{m}|^2+\int_s^T{\alpha \exp(\alpha r)|\bar{Y}_{r}^{m}|^2\dr}\\
	    &\le C( \alpha ) \Bigl( \varpi ^2\Bigl( \frac{2A}{\rho} \Bigr) +\varpi ^2\Bigl( \frac{2A}{m-A} \Bigr) \Bigr) \\
        &\quad +\int_s^T{2\exp(\alpha r)|\bar{Y}_{r}^{m}|^2\dr}-\int_s^T{2\exp(\alpha r)\bar{Y}_{r-}^{m}\ip[\bigg]{\bar{Z}_{r}^{m}}{\dW_r-\rho\frac{\bar{Y}_{r-}^{m}\bar{Z}_{r}^{m}}{| \bar{Y}_{r-}^{m} || \bar{Z}_{r}^{m} |}\dr}}\\
	    &\quad -\int_s^T{\int_E{2\exp(\alpha r)\bar{Y}_{r-}^{m}\bar{K}_{r}^{m}( e ) \bigl( \tilmu( \dr,\de ) -\1_{\{ \bar{Y}_{r-}^{m}\int_E{\ell_r( e ) \bar{K}_{r}^{m}( e ) \lambda(\de)}>0 \}}\rho \ell_r( e ) \lambda(\de) \dr \bigr)}}.
    \end{aligned}
    \end{equation}
Choose \(\alpha=2\) in \eqref{eq:continuous-dependence-approx-girsanov}. By a Girsanov transformation argument similar to that used in \Cref{lem:bsde-existence-comparison}\,\ref{lem:bsde-existence-comparison-i}, we have
    \begin{equation*}
    |\bar{Y}_{s}^{m}|^2\le C\Bigl( \varpi ^2\Bigl( \frac{2A}{\rho} \Bigr) +\varpi ^2\Bigl( \frac{2A}{m-A} \Bigr) \Bigr) \xrightarrow{\rho\rightarrow \infty}C \varpi ^2\Bigl( \frac{2A}{m-A} \Bigr).
    \end{equation*}
 That is to say for any \(\eta \in L_{\Fcal_t}^{2},\ v\in \Ucal_{t,T}\), we have
    \begin{equation}\label{eq:continuous-dependence-uniform-approx}
    \sup_{t\le s\le T}| Y_{s}^{t,\eta ,v}-Y_{s}^{t,\eta ,v,m} |^2\le C\varpi^2 \Bigl( \frac{2A}{m-A} \Bigr).
    \end{equation}

    For the second term, our goal is to prove that for any fixed \(m\), as \(n\rightarrow \infty\),
    \[\E\Bigl[ \sup_{t\le s\le T}| Y_{s}^{t,\eta ^n,v^n,m}-Y_{s}^{t,\eta ,v,m} |^2 \Bigr] \rightarrow 0.\]
Set \[\bar{X}_{s}^{n}=X_{s}^{t,\eta ^n,v^n}-X_{s}^{t,\eta ,v},\ \bar{Y}_{s}^{n,m}=Y_{s}^{t,\eta ^n,v^n,m}-Y_{s}^{t,\eta ,v,m},\ \bar{Z}_{s}^{n,m}=Z_{s}^{t,\eta ^n,v^n,m}-Z_{s}^{t,\eta ,v,m},\]
    \[\bar{K}_{s}^{n,m}( e ) =K_{s}^{t,\eta ^n,v^n,m}( e ) -K_{s}^{t,\eta ,v,m}( e ),\ \bar{v}_{r}^{n} = v_{r}^{n} - v_{r},\ \bar{\eta}^n = \eta ^n - \eta.\]
Applying It\^o's formula to \(\exp(\alpha s){|\bar{Y}_{s}^{n,m}|}^{2}\) on \([s,T]\), we obtain
    \begin{equation}\label{eq:continuous-dependence-generator-ito}
    \begin{aligned}
	    &\quad \exp(\alpha s)|\bar{Y}_{s}^{n,m}|^2+\int_s^T{\alpha \exp(\alpha r)|\bar{Y}_{r}^{n,m}|^2\dr}+\int_s^T{\exp(\alpha r)|\bar{Z}_{r}^{n,m}|^2\dr}+\int_s^T{\int_E{\exp(\alpha r)|\bar{K}_{r}^{n,m}( e ) |^2\lambda(\de) \dr}}\\
	    &\quad +\int_s^T{\int_E{\exp(\alpha r)\bigl( |\bar{K}_{r}^{n,m}( e ) |^2+2\bar{Y}_{r-}^{n,m}\bar{K}_{r}^{n,m}( e ) \bigr) \tilmu( \dr,\de )}}+\int_s^T{2\exp(\alpha r)\bar{Y}_{r-}^{n,m}\ip[\big]{\bar{Z}_{r}^{n,m}}{\dW_r}}\quad\\
	    &=\exp(\alpha T)|h( X_{T}^{t,\eta ^n,v^n} ) -h( X_{T}^{t,\eta ,v} ) |^2\\
        &\quad +\int_s^T{2\exp(\alpha r)\bar{Y}_{r}^{n,m}\Biggl( f_m\biggl( r,X_{r}^{t,\eta ^n,v^n},Y_{r}^{t,\eta ^n,v^n,m},Z_{r}^{t,\eta ^n,v^n,m},\int_E{\ell_r( e ) K_{r}^{t,\eta ^n,v^n,m}( e ) \lambda(\de)},v_{r}^{n} \biggr)}\\
        &\quad - f_m\biggl( r,X_{r}^{t,\eta ,v},Y_{r}^{t,\eta ,v,m},Z_{r}^{t,\eta ,v,m},\int_E{\ell_r( e ) K_{r}^{t,\eta ,v,m}( e ) \lambda(\de)},v_r \biggr) \Biggr) \dr.
    \end{aligned}
    \end{equation}
Due to the definition of \(f_m(\cdot)\), we have
    \begin{equation}\label{eq:continuous-dependence-generator-lipschitz}
    \begin{aligned}
	    &\quad \bar{Y}_{r}^{n,m}\Biggl( f_m\biggl( r,X_{r}^{t,\eta ^n,v^n},Y_{r}^{t,\eta ^n,v^n,m},Z_{r}^{t,\eta ^n,v^n,m},\int_E{\ell_r( e ) K_{r}^{t,\eta ^n,v^n,m}( e ) \lambda(\de)},v_{r}^{n} \biggr)  \\
        &\quad -f_m\biggl( r,X_{r}^{t,\eta ,v},Y_{r}^{t,\eta ,v,m},Z_{r}^{t,\eta ,v,m},\int_E{\ell_r( e ) K_{r}^{t,\eta ,v,m}( e ) \lambda(\de)},v_{r} \biggr) \Biggr)\\
	    &\le |\bar{Y}_{r}^{n,m}|\Biggl( \varpi(|\bar{X}_{r}^{n}|)+m\biggl( |\bar{Z}_{r}^{n,m}|+\int_E{\ell_r( e ) |\bar{K}_{r}^{n,m}( e ) |\lambda(\de)} \biggr) +\varpi(|\bar{v}^n_r|) \Biggr)\\
	    &\le \varpi^2(|\bar{X}_{r}^{n}|)+C(m) |\bar{Y}_{r}^{n,m}|^2+\frac{1}{4}\biggl( |\bar{Z}_{r}^{n,m}|^2+\int_E{|\bar{K}_{r}^{n,m}( e ) |^2\lambda(\de)} \biggr)+\varpi^2(|\bar{v}^n_r|).
    \end{aligned}
    \end{equation}
Combining \eqref{eq:continuous-dependence-generator-ito} with \eqref{eq:continuous-dependence-generator-lipschitz}, and then choosing a sufficiently large \(\alpha(m,L)\), we have
    \begin{equation}\label{eq:continuous-dependence-combined}
    \begin{aligned}
	    &\quad \exp(\alpha s)|\bar{Y}_{s}^{n,m}|^2+\int_s^T{\exp(\alpha r)|\bar{Z}_{r}^{n,m}|^2\dr}+\int_s^T{\int_E{\exp(\alpha r)|\bar{K}_{r}^{n,m}( e ) |^2\lambda(\de)\dr}}\\
	    &\quad +\int_s^T{\int_E{\exp(\alpha r)\bigl(|\bar{K}_{r}^{n,m}( e ) |^2+ 2\bar{Y}_{r-}^{n,m}\bar{K}_{r}^{n,m}( e ) \bigr)  \tilmu( \dr,\de )}}+2\int_s^T{\exp(\alpha r)\bar{Y}_{r-}^{n,m}\ip[\big]{\bar{Z}_{r}^{n,m}}{\dW_r}}\quad\\
        &\le C \exp(\alpha T)\sup_{t\le s\le T}|\bar{X}_{s}^{n}|^2+\frac{1}{2}\int_s^T{\exp(\alpha r)\biggl( |\bar{Z}_{r}^{n,m}|^2+\int_E{|\bar{K}_{r}^{n,m}( e ) |^2\lambda(\de)} \biggr) \dr}\\
        &\quad +C\int_s^T{\exp(\alpha r)\bigl[\varpi^2(|\bar{X}_{r}^{n}|)+\varpi^2(|\bar{v}^n_{r}|)\bigr]\dr}.
    \end{aligned}
    \end{equation}
    Taking expectations on both sides of \eqref{eq:continuous-dependence-combined}, we have
    \begin{equation}\label{eq:continuous-dependence-zk-estimate}
    \begin{aligned}
        &\quad \E\biggl[ \int_s^T{\exp(\alpha r)|\bar{Z}_{r}^{n,m}|^2\dr} \biggr] +\E\biggl[ \int_s^T{\int_E{\exp(\alpha r)|\bar{K}_{r}^{n,m}( e ) |^2\lambda(\de) \dr}} \biggr] \\
        &\le C \exp(\alpha T)\E\Bigl[ \sup_{t\le s\le T}|\bar{X}_{s}^{n}|^2 \Bigr] +C\E\int_s^T{\exp(\alpha r)\bigl[\varpi^2(|\bar{X}_{r}^{n}|)+\varpi^2(|\bar{v}^n_{r}|)\bigr]\dr}.
    \end{aligned}
    \end{equation}
Applying the BDG inequality to \eqref{eq:continuous-dependence-combined} and then combining the result with \eqref{eq:continuous-dependence-zk-estimate}, we further have
    \begin{equation}\label{eq:continuous-dependence-bdg-estimate}
    \begin{aligned}
	    &\quad \E\Bigl[ \sup_{t\le s\le T}\exp(\alpha s)|\bar{Y}_{s}^{n,m}|^2 \Bigr]\\
	    &\le C \exp(\alpha T)\E\Bigl[ \sup_{t\le s\le T}|\bar{X}_{s}^{n}|^2 \Bigr] +C\E\int_t^T{\exp(\alpha r)\bigl[\varpi^2(|\bar{X}_{r}^{n}|)+\varpi^2(|\bar{v}^n_{r}|)\bigr]\dr}\\
	    &\quad +\E\sup_{t\le s\le T}\biggl| \int_s^T{\exp(\alpha r)\bar{Y}_{r-}^{n,m}\ip[\big]{\bar{Z}_{r}^{n,m}}{\dW_r}} \biggr| +\E\sup_{t\le s\le T}\biggl| \int_s^T{\int_E{\exp(\alpha r)\bar{Y}_{r-}^{n,m}\bar{K}_{r}^{n,m}( e ) \tilmu( \dr,\de )}} \biggr|\\
	    &\le C \exp(\alpha T)\E\Bigl[ \sup_{t\le s\le T}|\bar{X}_{s}^{n}|^2 \Bigr] +C\E\int_t^T{\exp(\alpha r)\bigl[\varpi^2(|\bar{X}_{r}^{n}|)+\varpi^2(|\bar{v}^n_{r}|)\bigr]\dr} +\frac{1}{2}\E\Bigl[ \sup_{t\le s\le T}\exp(\alpha s)|\bar{Y}_{s}^{n,m}|^2 \Bigr]\\
	    &\quad +C\E\biggl[ \int_t^T{\exp(\alpha r)|\bar{Z}_{r}^{n,m}|^2\dr} \biggr] +C\E\biggl[ \int_t^T{\int_E{\exp(\alpha r)|\bar{K}_{r}^{n,m}( e ) |^2\lambda(\de) \dr}} \biggr].
    \end{aligned}
    \end{equation}
By \eqref{eq:continuous-dependence-zk-estimate} and \eqref{eq:continuous-dependence-bdg-estimate}, we have
    \begin{equation}\label{eq:continuous-dependence-sup-estimate}
    \quad \E\Bigl[ \sup_{t\le s\le T}\exp(\alpha s)|\bar{Y}_{s}^{n,m}|^2 \Bigr] \le C( m) \biggl( \E\Bigl[ \sup_{t\le s\le T}|\bar{X}_{s}^{n}|^2 \Bigr] +\E\int_t^T{\bigl(\varpi^2(|\bar{X}_{r}^{n}|)+\varpi^2(|\bar{v}^n_{r}|)\bigr)\dr} \biggr).
    \end{equation}
    On the other hand, we have the following continuous dependence property of SDEs:
    \begin{equation}\label{eq:sde-continuous-dependence}
    \E\Bigl[ \sup_{t\le s\le T}|\bar{X}_{s}^{n}|^2 \Bigr] \le C\biggl[ \E| \eta ^n-\eta |^2+\E\int_t^T{| v_{r}^{n}-v_r |^2\dr} \biggr].
    \end{equation}

    Therefore, combining \eqref{eq:continuous-dependence-uniform-approx}, \eqref{eq:continuous-dependence-sup-estimate}, \eqref{eq:sde-continuous-dependence} and \Cref{lem:modulus-uniform-integrability}, we arrive at
    \begin{equation}\label{eq:continuous-dependence-final}
    \begin{aligned}
        &\quad \E\Bigl[ \sup_{t\le s\le T}|Y_{s}^{t,\eta ^n,v^n}-Y_{s}^{t,\eta ,v}|^2 \Bigr] \\
        &\le 3\E\Bigl[ \sup_{t\le s\le T}|Y_{s}^{t,\eta ^n,v^n}-Y_{s}^{t,\eta ^n,v^n,m}|^2 \Bigr] +3\E\Bigl[ \sup_{t\le s\le T}|Y_{s}^{t,\eta ^n,v^n,m}-Y_{s}^{t,\eta ,v,m}|^2 \Bigr]
        +3\E\Bigl[ \sup_{t\le s\le T}|Y_{s}^{t,\eta ,v,m}-Y_{s}^{t,\eta ,v}|^2 \Bigr] \\
        &\le C\varpi^2 \Bigl( \frac{2A}{m-A} \Bigr) +C(m)\Biggl( \E| \bar{\eta} ^n |^2+\E\int_t^T{| \bar{v}_{r}^{n} |^2\dr} +\varpi_\ast\biggl( C\E| \bar{\eta} ^n |^2+C\E\int_t^T{| \bar{v}_{r}^{n} |^2\dr} \biggr) \Biggr).
    \end{aligned}
    \end{equation}
    Here \(\varpi_\ast\) denotes the function \(\omega_\ast\) defined in \Cref{lem:modulus-uniform-integrability} with \(\omega=\varpi\).
    Taking \(n\rightarrow \infty\) for fixed \(m\) first, and then \(m\rightarrow \infty\), we complete the proof of \Cref{thm:bsde-continuous-dependence}.

\end{proof}

\begin{theorem}[Uniform initial stability]\label{thm:bsde-uniform-initial-stability}
    Under the same assumptions as \Cref{thm:bsde-continuous-dependence}, we have the following two conclusions.

    \begin{enumerate}[label=\textup{(\roman*)},ref=\textup{(\roman*)},wide=0pt]
    \item\label{thm:bsde-uniform-initial-stability-i} For any \(m>A\), there exists a deterministic nondecreasing
    function \(\Omega_m:\R_+\to\R_+\), independent of \(t\) and \(v\), such
    that \( \Omega_m(r)\to0\) as \(r\downarrow0 \), and the following estimate holds uniformly in \(t\in[0,T]\):
    \begin{equation*}
    \begin{aligned}
        &\quad \esssup_{v\in \Ucal_{t,T}}| Y_{t}^{t,\eta',v}-Y_{t}^{t,\eta ,v} |\\
        &\le \esssup_{v\in \Ucal_{t,T}}| Y_{t}^{t,\eta',v}-Y_{t}^{t,\eta',v,m} |+\esssup_{v\in \Ucal_{t,T}}| Y_{t}^{t,\eta',v,m}-Y_{t}^{t,\eta ,v,m} |+\esssup_{v\in \Ucal_{t,T}}| Y_{t}^{t,\eta ,v,m}-Y_{t}^{t,\eta ,v} |\\
        &\le C\varpi\Bigl(\frac{2A}{m-A}\Bigr)+\Omega_m(|\eta'-\eta|),\ \Pp\text{-a.s.}
    \end{aligned}
    \end{equation*}

    \item\label{thm:bsde-uniform-initial-stability-ii} Let \(\xi \in L_{\Fcal_T}^{2}(\Omega;\R)\)
    and \(\{ \xi _n \} _{n\ge 1} \subset L_{\Fcal_T}^{2}(\Omega;\R)\)
    satisfy \(\xi_n \rightarrow \xi\) in \(L^2_{\Fcal_T}(\Omega;\R)\).
    Consider the following BSDE with jumps:
    \begin{equation}\label{eq:bsde-terminal-stability}
    \begin{aligned}
        Y_s&=\zeta +\int_s^T{f\biggl( r,X_{r}^{t,x,v},Y_r,Z_r,\int_E{\ell_r( e ) K_r( e ) \lambda(\de)},v_r \biggr)\dr}\\
        &\quad -\int_s^T{\ip{Z_r}{\dW_r}}-\int_s^T{\int_E{K_r( e ) \tilmu( \dr,\de )}},\ s\in [ t,T ].
    \end{aligned}
    \end{equation}
If \((Y^n,Z^n,K^n)\) and \((Y,Z,K)\)
    are corresponding solutions of BSDEs \eqref{eq:bsde-terminal-stability} with parameters \((\xi_n,f)\) and \((\xi,f)\), respectively,
then
    \begin{align*}
        \underset{n\rightarrow \infty}{\lim}\E\sup_{t\le s\le T}| Y_{s}^{n}-Y_s |^2\le \underset{m\rightarrow \infty}{\lim}\underset{n\rightarrow \infty}{\lim}\Bigl( C\varpi ^2\Bigl( \frac{2A}{m-A} \Bigr) +C( m ) \E| \xi _n-\xi |^2 \Bigr) =0.
    \end{align*}
    \end{enumerate}
\end{theorem}

\begin{proof}
    Since the proof is similar to that of \Cref{thm:bsde-continuous-dependence}, we only give a sketch of the proof.

    \begin{enumerate}[label=\textup{(\roman*)},ref=\textup{(\roman*)},wide=0pt]
    \item For the first and third terms, we can use the same method as in \Cref{thm:bsde-continuous-dependence} to obtain the following estimates:
    \begin{equation*}
    | Y_{t}^{t,\eta',v}-Y_{t}^{t,\eta',v,m} | + | Y_{t}^{t,\eta ,v,m}-Y_{t}^{t,\eta ,v} | \le C\varpi \Bigl( \frac{2A}{m-A} \Bigr),\quad\text{for any }v\in \Ucal_{t,T}.
    \end{equation*}
As for the second term, an estimate similar to \eqref{eq:continuous-dependence-combined} gives
    \begin{equation}\label{eq:uniform-stability-ito}
    \begin{aligned}
        &\quad |Y_{t}^{t,\eta',v,m}-Y_{t}^{t,\eta ,v,m}|^2 \\
        &\le C \exp(\alpha T)\sup_{t\le s\le T}|X_{s}^{t,\eta',v}-X_{s}^{t,\eta ,v}|^2+C\int_t^T{\exp(\alpha r)\varpi^2(|X_{r}^{t,\eta',v}-X_{r}^{t,\eta ,v}|)\dr}\\
        &\quad -\int_t^T{2\exp(\alpha r)\bigl( Y_{r-}^{t,\eta',v,m}-Y_{r-}^{t,\eta ,v,m} \bigr) \ip[\big]{\bigl( Z_{r}^{t,\eta',v,m}-Z_{r}^{t,\eta ,v,m} \bigr)}{\dW_r}}\\
        &\quad-\int_t^T\int_E \exp(\alpha r)\left[
        |K_{r}^{t,\eta',v,m}(e)-K_{r}^{t,\eta ,v,m}(e)|^2
        +2\bigl(Y_{r-}^{t,\eta',v,m}-Y_{r-}^{t,\eta ,v,m}\bigr)
        \bigl(K_{r}^{t,\eta',v,m}(e)-K_{r}^{t,\eta ,v,m}(e)\bigr)
        \right]\tilmu(\dr,\de).
    \end{aligned}
    \end{equation}
    Taking conditional expectations on both sides of \eqref{eq:uniform-stability-ito}, we obtain
    \begin{equation*}
    |Y_t^{t,\eta';v,m}-Y_t^{t,\eta;v,m}|^2 \le C(m)\E^{\Fcal_t}\biggl[\sup_{t\le s\le T}|X_s^{t,\eta';v}-X_s^{t,\eta;v}|^2+\int_t^T\varpi^2(|X_s^{t,\eta';v}-X_s^{t,\eta;v}|)\,\ds\biggr].
    \end{equation*}

By the conditional SDE estimate and
\Cref{lem:modulus-uniform-integrability}, there exists a deterministic
modulus \(\omega_\ast\), independent of \(t\) and \(v\), such that
\[|Y_t^{t,\eta';v,m}-Y_t^{t,\eta;v,m}|^2\le C(m)(|\eta'-\eta|^2+\omega_\ast(C|\eta'-\eta|^2)).\]
Thus one may take \( \Omega_m(r)=\sqrt{C(m)(r^2+\omega_\ast(Cr^2))}\).

    \item Here we use a similar decomposition
    \begin{align*}
        &\quad\E\Bigl[ \sup_{t\le s\le T}|Y_{s}^{n}-Y_s|^2 \Bigr]\\
         &\le C\E\Bigl[ \sup_{t\le s\le T}|Y_{s}^{n}-Y_{s}^{n,m}|^2 \Bigr] +C\E\Bigl[ \sup_{t\le s\le T}|Y_{s}^{n,m}-Y_{s}^{m}|^2 \Bigr] +C\E\Bigl[ \sup_{t\le s\le T}|Y_{s}^{m}-Y_s|^2 \Bigr],
    \end{align*}
    where \((Y^m,Z^m,K^m)\) and \((Y^{n,m},Z^{n,m},K^{n,m})\) are corresponding solutions of BSDEs \eqref{eq:bsde-terminal-stability} with parameters  \((\xi,f_m)\) and \((\xi_n,f_m)\), respectively.
    For the first and third terms, a similar argument as in the proof of \Cref{thm:bsde-continuous-dependence} leads to
    \[
    \sup_{t\le s\le T}| Y_{s}^{n}-Y_{s}^{n,m} |^2+\sup_{t\le s\le T}| Y_{s}^{m}-Y_s |^2\le C\varpi ^2\Bigl( \frac{2A}{m-A} \Bigr).
    \]
    On the other hand, we can similarly apply It\^o's formula to \(\exp(\alpha s){|Y_{s}^{n,m}-Y_{s}^{m}|}^{2}\) with sufficiently large \(\alpha(m,L)\) to get
    \begin{equation*}
    \E\Bigl[ \sup_{t\le s\le T}|Y_{s}^{n,m}-Y_{s}^{m}|^2 \Bigr] \le C(m) \E| \xi _n-\xi |^2.
    \end{equation*}
Then \ref{thm:bsde-uniform-initial-stability-ii} follows.
    \end{enumerate}
\end{proof}

\section{Dynamic programming principle}\label{sec:dpp}

We return to the controlled FBSDE \eqref{eq:controlled-fbsde}. In the stochastic system \eqref{eq:controlled-fbsde}, the cost functional is defined as \(J(t,x;v):=Y^{t,x;v}_t\). 
The recursive control problem is to maximize the cost functional
\(J(t,x;v)\) over \(v\in\Ucal_{t,T}\), and its value function is
defined by
\begin{equation}\label{eq:value-function}
    u(t,x):=
    \esssup_{v\in\Ucal_{t,T}}J(t,x;v).
\end{equation}

To prove the dynamic programming principle for the control problem \eqref{eq:value-function}, we need to make some preparations. To be specific, we first show that the value function is deterministic and establish its basic growth and continuity properties, and then extend the value function from deterministic initial states to \(\Fcal_t\)-measurable initial data in the following three theorems.

\begin{theorem}[Deterministic value function]\label{thm:value-function-deterministic}
    Suppose that \Cref{ass:standing} holds. Then for any \((t,x)\in [0,T]\times \R^{d'}\),
    the value function \(u(t,x)\) is deterministic.
\end{theorem}
\begin{proof}
    We divide the proof into two steps. To begin with, we introduce the following two classes of admissible controls:
    \begin{align*}
	    \Ucal^t:&=\bigl\{ v\in \Ucal_{t,T}\bigm| v_s\text{ is }\Fcal_{s}^{t}\text{-measurable}  \bigr\},\\
	    \bar{\Ucal}^t:&=\Bigl\{ v\in \Ucal_{t,T}\Bigm| v_s=\sum_{j=1}^N{v_{s}^{j}\1_{A_j}},\,\,v^j\in \Ucal^t  \Bigr\},
    \end{align*}
  where \(\{A_j\}_{j=1}^{N}\subset\Fcal_t\) is a partition of \(\Omega\), and
    \[    
    \Fcal_{s}^{t}=\sigma ( W_r-W_t:r\in [ t,s ] )
    \lor \sigma ( \mu ( ( t,r ] \times A ) :r\in [ t,s ] ,A\in \Bcal( E ) ) \lor \mathcal{N}.
    \]
Denote by \(\mathbb{F}^t=\bigl\{ \Fcal_{s}^{t} \bigr\}_{s\in [t,T]}\) the shifted filtration.

    The first step is to prove
    \begin{equation}\label{eq:esssup-equal-1}
    \esssup_{v\in \Ucal_{t,T}}J(t,x;v)=\esssup_{v\in \bar{\Ucal}^t}J(t,x;v).
    \end{equation}
Obviously, \(\esssup_{v\in \Ucal_{t,T}}J(t,x;v)\ge \esssup_{v\in \bar{\Ucal}^t}J(t,x;v)\) since \(\bar{\Ucal}^t\subset \Ucal_{t,T}\).
    For the other direction, we notice that \(\bar{\Ucal}^t\) is dense in \(\Ucal_{t,T}\) (see \cite[Proposition~3.1]{LiPeng2009}). So for any \(v\in \Ucal_{t,T}\), there exists a sequence \(\{ v^n \} _{n=1}^{\infty}\subset \bar{\Ucal}^t\), such that \(v^n\) converges to \(v\) in \(\mathcal{L}^2\).
    By \Cref{thm:bsde-continuous-dependence}, we have
    \begin{align*}
    \underset{n\rightarrow \infty}{\lim}\E| Y_{t}^{t,x;v^n}-Y_{t}^{t,x;v} |^2\le \underset{n\rightarrow \infty}{\lim}\E\Bigl[ \sup_{t\le s\le T}|Y_{s}^{t,x,v^n}-Y_{s}^{t,x,v}|^2 \Bigr] =0.
    \end{align*}
    Without loss of generality, assume that there exists a subsequence still denoted by \(\{ v^n \} _{n=1}^{\infty}\), such that \(\underset{n\rightarrow \infty}{\lim}Y_{t}^{t,x;v^n} = Y_{t}^{t,x;v}\), \(\Pp\text{-a.s.}\)
    Due to the arbitrariness of \(v \in \Ucal_{t,T}\), we obtain
    \[\esssup_{v\in \Ucal_{t,T}}J(t,x;v)\le \esssup_{v\in \bar{\Ucal}^t}J(t,x;v).\]

    The second step is to prove
    \begin{equation}\label{eq:esssup-equal-2}
    \esssup_{v\in \bar{\Ucal}^t}J(t,x;v)=\esssup_{v\in \Ucal^t}J(t,x;v).
    \end{equation}
    The proof is standard. Here we only show the key computation.
    For any \(v\in \bar{\Ucal}^t\), we can write it as \(v_s=\sum_{j=1}^N{v_{s}^{j}\1_{A_j}}\), where \(v^j\in \Ucal^t\) and \(\{ A_j \} _{j=1}^{N}\subset \Fcal_t\) is a partition of \(( \Omega ,\Fcal_t )\).
    By the uniqueness of the solution to FBSDE \eqref{eq:controlled-fbsde}, we have
    \begin{align*}
        Y_{s}^{t,x;\sum_{j=1}^N{v^j\1_{A_j}}}&=\sum_{j=1}^N{\1_{A_j}Y_{s}^{t,x;v^j}},\ s\in [ t,T ].
    \end{align*}
    Consequently,
    \(
        J(t,x;v)=J\bigl(t,x;\sum_{j=1}^N{v^j\1_{A_j}}\bigr)=\sum_{j=1}^N{\1_{A_j}J(t,x;v^j)}.
    \)
    Since \(v^j\in \Ucal^t\), \(Y_{s}^{t,x;v^j}\) is \(\Fcal_{s}^{t}\)-measurable and \(Y_{t}^{t,x;v^j}\) is a constant. Thus, we have
    \begin{align*}
        J(t,x;v) \le \max_{1 \le j \le N}J(t,x;v^j) \le \esssup_{v\in \Ucal^t}J(t,x;v).
    \end{align*}
Due to the arbitrariness of \(v\in \bar{\Ucal}^t\), we have \(\esssup_{v\in \bar{\Ucal}^t}J(t,x;v)\le \esssup_{v\in \Ucal^t}J(t,x;v)\). The other direction is obvious since \(\Ucal^t\subset \bar{\Ucal}^t\). Thus, we complete the second step.

    Combining \eqref{eq:esssup-equal-1} with \eqref{eq:esssup-equal-2}, we obtain
    \begin{align*}
        \esssup_{v\in \Ucal_{t,T}}J(t,x;v)&=\esssup_{v\in \bar{\Ucal}^t}J(t,x;v)=\esssup_{v\in \Ucal^t}J(t,x;v).
    \end{align*}
    Since \(J(t,x;v)=Y_{t}^{t,x;v}\) is a constant for any \(v\in \Ucal^t\),
    we conclude that \(u(t,x)\) is deterministic.
\end{proof}

\begin{theorem}[Linear growth and spatial continuity of the value function]\label{thm:value-growth-spatial-continuity}
Suppose that \Cref{ass:standing} holds. Then there exists a constant \(C>0\) such that, for any \(t\in [ 0,T ]\) and \(x,x'\in \R^{d'}\),

    \begin{enumerate}[label=\textup{(\roman*)},ref=\textup{(\roman*)},wide=0pt]
    \item\label{thm:value-growth-spatial-continuity-i} \(|u(t,x)|\le C(1+|x|)\);

    \item\label{thm:value-growth-spatial-continuity-ii} \(|u(t,x)-u(t,x')| \rightarrow 0\), as \(|x-x'|\rightarrow 0\).
    \end{enumerate}
\end{theorem}

\begin{proof}
    For any \(\delta >0\), there exist controls \(v_1^\delta ,v_2^\delta \in \Ucal^t\) such that
    \begin{equation*}
    \begin{aligned}
        Y_t^{t,x;v_2^\delta} &\le u(t,x) \le Y_t^{t,x;v_1^\delta}+\delta\ \ \ \text{and}\ \ \ Y_t^{t,x';v_1^\delta} &\le u(t,x') \le Y_t^{t,x';v_2^\delta}+\delta.
    \end{aligned}
    \end{equation*}
    By the a priori estimate in \Cref{lem:bsde-existence-comparison}\,\ref{lem:bsde-existence-comparison-iii}, we obtain
    \[|Y_t^{t,x;v_i^\delta}|\le C(1+|x|),\ i=1,2,\ x \in \R^{d'},\]
    where \(C\) is independent of \(\delta,t,x\) and the choice of \(v_i^\delta\).
    Thus, we have
    \begin{align*}
        |u(t,x)| \le \max \bigl\{ |Y_t^{t,x;v_1^\delta}|+\delta ,|Y_t^{t,x;v_2^\delta}| \bigr\} \le C(1+|x|)+\delta.
    \end{align*}
    Letting \(\delta \downarrow 0\), we obtain \(|u(t,x)| \le C(1+|x|)\). This proves \ref{thm:value-growth-spatial-continuity-i}.

To see \ref{thm:value-growth-spatial-continuity-ii}, by \Cref{thm:bsde-uniform-initial-stability}\,\ref{thm:bsde-uniform-initial-stability-i}, we have
    \begin{equation}\label{eq:value-growth-spatial-continuity}
    \begin{aligned}
	    |u( t,x ) -u( t,x' ) |
        \le \esssup_{v\in \Ucal^t}|Y_{t}^{t,x;v}-Y_{t}^{t,x';v}|\le C\varpi \Bigl( \frac{2A}{m-A} \Bigr) +\Omega_m(|x-x'|).
    \end{aligned}
    \end{equation}
Taking \(|x-x'|\rightarrow 0\) and then \(m\rightarrow \infty\), we obtain the desired result.
\end{proof}

\begin{theorem}[Value function for random initial data]\label{thm:value-random-initial-data}
Suppose that \Cref{ass:standing} holds.  Then for any \(t\in[0,T]\), \(\eta\in L^2_{\Fcal_t}(\Omega;\R^{d'})\)
 and \(v(\cdot) \in \Ucal_{t,T}\), we have
    \begin{align*}
        u(t,\eta) \ge Y^{t,\eta;v}_t,\ \Pp\text{-a.s.}
    \end{align*}
Conversely, for any \(\delta >0\), there exists \(v^\delta (\cdot) \in \Ucal_{t,T}\) such that
    \begin{align*}
        u(t,\eta) \le Y^{t,\eta;v^\delta}_t +\delta ,\ \Pp\text{-a.s.}
    \end{align*}
\end{theorem}
\begin{proof}
    We divide the proof into three steps.

    First, we prove that the conclusion holds for simple random variables \(\zeta = \sum_{i=1}^N{x_i\1_{A_i}}\), where \(x_i \in \R^{d'}\) and \(\{ A_i \} _{i=1}^{N}\subset \Fcal_t\) is a partition of \(\Omega\).
    For any \(v(\cdot) \in \Ucal_{t,T}\), by the uniqueness of the solution to FBSDE \eqref{eq:controlled-fbsde} and the definition of the value function, we have
    \begin{align*}
        Y^{t,\zeta ;v}_t &=\sum_{i=1}^N{\1_{A_i}Y_{t}^{t,x_i;v}} \le \sum_{i=1}^N{\1_{A_i}u(t,x_i)} = u(t,\zeta).
    \end{align*}
    Conversely, for any \(\delta >0\) and each \(x_i \in \R^{d'}\),
    there exists a \(v_i^\delta (\cdot) \in \Ucal^t\) such that
    \(u(t,x_i) \le Y_{t}^{t,x_i;v_i^\delta} +\delta\).
Set \(v^\delta (\cdot)=\sum_{i=1}^N{v_i^\delta (\cdot)\1_{A_i}}\).
    Then \(v^\delta (\cdot) \in \Ucal_{t,T}\) and
    \begin{align*}
        Y^{t,\zeta ;v^\delta}_t + \delta &= \sum_{i=1}^N{\1_{A_i} \bigl( Y_{t}^{t,x_i;v_i^\delta} +\delta \bigr)} \ge \sum_{i=1}^N{\1_{A_i}u(t,x_i)} = u(t,\zeta).
    \end{align*}
Hence we complete the proof for simple random variables.

    We next prove the first inequality for a general random variable \(\eta \in L^2_{\Fcal_t}(\Omega;\R^{d'})\).
    Let \(\{ \zeta _n \} _{n\ge 1} \subset L^2_{\Fcal_t}(\Omega;\R^{d'})\) be a sequence of simple random variables converging to \(\eta\) in \(L^2_{\Fcal_t}(\Omega;\R^{d'})\).
    Then, for any \(v \in \Ucal_{t,T}\), it follows from \Cref{thm:bsde-continuous-dependence} that
    \begin{align*}
        \E|Y_{t}^{t,\zeta _n;v}-Y_{t}^{t,\eta ;v}|^2 &\le \E\Bigl[ \sup_{t\le s\le T}|Y_{s}^{t,\zeta _n;v}-Y_{s}^{t,\eta ;v}|^2 \Bigr] \rightarrow 0,\quad\text{as }n\rightarrow \infty.
    \end{align*}
    Hence, by passing to a subsequence if necessary, we may assume that 
    \(\zeta_n\longrightarrow\eta\) and \(Y_{t}^{t,\zeta _n;v}\rightarrow Y_{t}^{t,\eta ;v}\), \(\Pp\text{-a.s.}\)
    Since \(x\mapsto u(t,x)\) is a deterministic continuous function by
    \Cref{thm:value-growth-spatial-continuity}, it follows that
    \(u(t,\zeta_n)\longrightarrow u(t,\eta),\ \Pp\text{-a.s.}\)
    For every \(n\), the result proved for simple random variables leads to
    \(Y_t^{t,\zeta_n;v}\le u(t,\zeta_n),\ \Pp\text{-a.s.}\)
    Letting \(n\to\infty\), we obtain
    \(Y_t^{t,\eta;v}\le u(t,\eta),\ \Pp\text{-a.s.}\)

    Finally, we prove the second inequality. We first consider the situation of bounded random variables \(\eta \in L^\infty _{\Fcal_t}(\Omega;\R^{d'})\).
    Fix \(\delta >0\). Choose \(m>A\) sufficiently large so that \(C\varpi\bigl(2A/(m-A)\bigr)\le \delta/6\), and then choose \(r_m>0\) so that \(\Omega_m(r_m)\le \delta/6\).
    Since \(\eta\) is bounded, there exists a simple random variable \(\zeta\), such that \(|\eta-\zeta|\le r_m\), \(\Pp\text{-a.s.}\),
    where \(\Omega_m\) is the function in \Cref{thm:bsde-uniform-initial-stability}\,\ref{thm:bsde-uniform-initial-stability-i}.
    By \Cref{thm:bsde-uniform-initial-stability}\,\ref{thm:bsde-uniform-initial-stability-i}, it follows that
    \begin{align*}
        \esssup_{v\in \Ucal_{t,T}}|Y_{t}^{t,\eta ;v}-Y_{t}^{t,\zeta ;v}| &\le C\varpi \Bigl( \frac{2A}{m-A} \Bigr) +\Omega_m(|\eta -\zeta|) \le \frac{\delta}{3},\ \Pp\text{-a.s.}
    \end{align*}
    Applying \eqref{eq:value-growth-spatial-continuity} pointwise with
    \(x=\eta(\omega)\) and \(x'=\zeta(\omega)\), we obtain
    \[
        |u(t,\eta)-u(t,\zeta)|
        \le C\varpi\Bigl(\frac{2A}{m-A}\Bigr) +\Omega_m(|\eta-\zeta|)
        \le \frac{\delta}{3},
        \qquad \Pp\text{-a.s.}
    \]

    Since \(\zeta\) is a simple random variable, there exists a \(v^\delta (\cdot) \in \Ucal_{t,T}\) such that
    \begin{align*}
        u(t,\zeta) \le Y_{t}^{t,\zeta ;v^\delta} +\frac{\delta}{3},\ \Pp\text{-a.s.}
    \end{align*}
    Thus we obtain
    \begin{align*}
        u(t,\eta) \le u(t,\zeta) + \frac{\delta}{3} &\le Y_{t}^{t,\zeta ;v^\delta} + \frac{2\delta}{3} \le Y_t^{t,\eta ;v^\delta} + \delta ,\ \Pp\text{-a.s.}
    \end{align*}
As for general random variables \(\eta \in L^2_{\Fcal_t}(\Omega;\R^{d'})\), we define
    \begin{align*}
        A_i=\{ \omega :i-1\le | \eta ( \omega ) |<i \} ,\ \eta _i=\eta \1_{A_i}( \omega ) ,\ i\in \mathbb{N}.
    \end{align*}
    Then \(\{A_i\}_{i\ge1}\subset\Fcal_t\) is a countable partition of \(\Omega\), and each \(\eta_i\) is bounded.
    Hence, by the bounded random variables case already proved, for any \(\delta >0\), there exists \(v_i^\delta (\cdot) \in \Ucal_{t,T}\) such that
    \begin{align*}
        u(t,\eta _i) \le Y_{t}^{t,\eta _i;v_i^\delta} +\delta ,\ \Pp\text{-a.s.}
    \end{align*}
    Set \(v^\delta (\cdot)=\sum\limits_{i=1}^{\infty}{v_i^\delta (\cdot)\1_{A_i}}\). Then we have \(v^\delta (\cdot) \in \Ucal_{t,T}\) and
    \begin{align*}
        u(t,\eta) &= \sum_{i=1}^{\infty}{\1_{A_i}u(t,\eta _i)} \le \sum_{i=1}^{\infty}{\1_{A_i}\bigl( Y_{t}^{t,\eta _i;v_i^\delta} +\delta \bigr)} = Y^{t,\eta ;v^\delta}_t + \delta ,\ \Pp\text{-a.s.}
    \end{align*}
This completes the proof.
\end{proof}

Our next goal is to establish the dynamic programming principle for the stochastic recursive control problem \eqref{eq:value-function}. To this end, we first introduce the following family of backward semigroups associated with BSDEs with jumps.
Given the initial data \((t,x)\in [0,T]\times \R^{d'}\), a real number \(\delta \in (0,T-t]\), an admissible control \(v(\cdot) \in \Ucal_{t,t+\delta}\) and a random variable \(\eta \in L^2_{\Fcal_{t+\delta}}\), we define the backward semigroup \(G_{t,t+\delta}^{t,x;v}[\eta]\) by
\begin{align*}
    G_{s,t+\delta}^{t,x;v}[\eta] := \tilde{Y}_s^{t,x;v},\ s\in [t,t+\delta],
\end{align*}
where \(\bigl(\tilde{Y}^{t,x;v},\tilde{Z}^{t,x;v},\tilde{K}^{t,x;v}\bigr)\) is the solution of the following BSDE with jumps:
\begin{equation*}
\begin{aligned}
    \tilde{Y}_s^{t,x;v} &= \eta + \int_s^{t+\delta}{f\biggl( r,X_{r}^{t,x;v},\tilde{Y}_r^{t,x;v},\tilde{Z}_r^{t,x;v},\int_E{\ell_r( e ) \tilde{K}_r^{t,x;v}( e )}\lambda(\de),v_r \biggr) \dr} \\
    &\quad - \int_s^{t+\delta}{\ip[\big]{\tilde{Z}_r^{t,x;v}}{\dW_r}} - \int_s^{t+\delta}{\int_E{\tilde{K}_r^{t,x;v}( e ) \tilmu( \dr,\de )}},\ s\in [t,t+\delta].
\end{aligned}
\end{equation*}
By the uniqueness of the solution to the BSDE with jumps, we have
\begin{align*}
    G_{t,T}^{t,x;v}\bigl[h\bigl(X_{T}^{t,x;v}\bigr)\bigr] = G_{t,t+\delta}^{t,x;v}\bigl[Y_{t+\delta}^{t,x;v}\bigr],
\end{align*}
which, together with the definition of the value function, leads to
\begin{align*}
    J(t,x;v)=Y_{t}^{t,x;v}=G_{t,t+\delta}^{t,x;v}\bigl[Y_{t+\delta}^{t,x;v}\bigr]=G_{t,t+\delta}^{t,x;v}\bigl[J\bigl(t+\delta,X_{t+\delta}^{t,x;v};v\bigr)\bigr].
\end{align*}

Now we are in a position to state and prove the dynamic programming principle for our stochastic recursive optimal control problem.
\begin{theorem}[Dynamic programming principle]\label{thm:dpp}
    Suppose that \Cref{ass:standing} holds. Then for any \((t,x)\in [0,T]\times \R^{d'}\) and \(0<\delta \le T-t\), we have
    \begin{align*}
        u(t,x) = \esssup_{v(\cdot) \in \Ucal_{t,t+\delta}}G_{t,t+\delta}^{t,x;v}\bigl[u\bigl(t+\delta,X_{t+\delta}^{t,x;v}\bigr)\bigr].
    \end{align*}
\end{theorem}
\begin{proof}
    We first prove
    \begin{align*}
        u(t,x) \le \esssup_{v(\cdot) \in \Ucal_{t,t+\delta}}G_{t,t+\delta}^{t,x;v}\bigl[u\bigl(t+\delta,X_{t+\delta}^{t,x;v}\bigr)\bigr].
    \end{align*}
    In view of the definition of the value function, we have
    \begin{align*}
        u(t,x)=\esssup_{v(\cdot) \in \Ucal_{t,T}}J(t,x;v)=\esssup_{v(\cdot) \in \Ucal_{t,T}}G_{t,t+\delta}^{t,x;v}\bigl[J\bigl(t+\delta,X_{t+\delta}^{t,x;v};v\bigr)\bigr].
    \end{align*}
In particular, for any \(v(\cdot) \in \Ucal_{t,T}\), \(J\bigl(t+\delta,X_{t+\delta}^{t,x;v};v\bigr) = Y^{t+\delta,X_{t+\delta}^{t,x;v};v}_{t+\delta} \le u\bigl(t+\delta,X_{t+\delta}^{t,x;v}\bigr)\), \(\Pp\)-a.s. Then we verify
the integrability of \(u\bigl(t+\delta,X_{t+\delta}^{t,x;v}\bigr)\). For this, first note that for any \(\varepsilon >0\), there exists \(v^\varepsilon (\cdot) \in \Ucal_{t+\delta,T}\), such that
    \begin{align*}
        Y^{t+\delta,X_{t+\delta}^{t,x;v};v^\varepsilon}_{t+\delta} \le u\bigl(t+\delta,X_{t+\delta}^{t,x;v}\bigr) \le Y^{t+\delta,X_{t+\delta}^{t,x;v};v^\varepsilon}_{t+\delta} +\varepsilon ,\ \Pp\text{-a.s.}
    \end{align*}
    By the a priori estimate in \Cref{lem:bsde-existence-comparison}\,\ref{lem:bsde-existence-comparison-iii}, we have
    \begin{equation*}
    \begin{aligned}
        \E|Y_{t+\delta}^{t+\delta ,X_{t+\delta}^{t,x;v};v^{\varepsilon}}|^2 &\le C\biggl( 1+\E| X_{t+\delta}^{t,x;v} |^2+\E\int_{t+\delta}^T{|f( r,0,0,0,0,v_{r}^{\varepsilon} ) |^2\dr} \biggr)\\
	    &\le C\biggl( 1+| x |^2+\E\int_{t+\delta}^T{|f( r,0,0,0,0,v_{r}^{\varepsilon} ) |^2\dr} \biggr)< \infty.	
    \end{aligned}
    \end{equation*}
    With the integrability of \(u\bigl(t+\delta,X_{t+\delta}^{t,x;v}\bigr)\) checked, by \Cref{lem:bsde-existence-comparison}\,\ref{lem:bsde-existence-comparison-i}, we have
    \begin{align*}
        G_{t,t+\delta}^{t,x;v}\bigl[J\bigl(t+\delta,X_{t+\delta}^{t,x;v};v\bigr)\bigr] \le G_{t,t+\delta}^{t,x;v}\bigl[u\bigl(t+\delta,X_{t+\delta}^{t,x;v}\bigr)\bigr],\ \Pp\text{-a.s.},
    \end{align*}
    for any \(v(\cdot) \in \Ucal_{t,T}\). Taking the essential supremum on both sides, we complete the first part of the proof.

    Next, we prove the other direction:
    \begin{align*}
        u(t,x) \ge \esssup_{v(\cdot) \in \Ucal_{t,t+\delta}}G_{t,t+\delta}^{t,x;v}\bigl[u\bigl(t+\delta,X_{t+\delta}^{t,x;v}\bigr)\bigr].
    \end{align*}

We first prove a simple estimate for the backward semigroup. 
Let \(v\in\Ucal_{t,t+\delta}\) and \(\xi\in L^2_{\Fcal_{t+\delta}}(\Omega;\R)\). 
Therefore, for any \(\varepsilon>0\),
\begin{equation}\label{eq:dpp-terminal-shift}
    0
    \le
    G_{t,t+\delta}^{t,x;v}[\xi]
    -
    G_{t,t+\delta}^{t,x;v}[\xi-\varepsilon]
    \le
    \varepsilon,
    \qquad
    \Pp\text{-a.s.}
\end{equation}
Indeed, the first inequality follows directly from the comparison
theorem. For the second one, let
\((Y^2,Z^2,K^2)\) be the solution defining
\(G_{t,t+\delta}^{t,x;v}[\xi-\varepsilon]\), and set
\[
    \widetilde Y_s=Y_s^2+\varepsilon,
    \qquad
    \widetilde Z_s=Z_s^2,
    \qquad
    \widetilde K_s=K_s^2.
\]
Then \(\widetilde Y_{t+\delta}=\xi\), and
\((\widetilde Y,\widetilde Z,\widetilde K)\) satisfies a BSDE whose
generator is
\(
    f\bigl(s,X_s^{t,x;v},y-\varepsilon,z,k,v_s\bigr).
\)
Since \(f\) is nonincreasing in \(y\), the comparison theorem yields
\[
    G_{t,t+\delta}^{t,x;v}[\xi]
    \le
    \widetilde Y_t
    =
    G_{t,t+\delta}^{t,x;v}[\xi-\varepsilon]
    +\varepsilon,
\]
which proves \eqref{eq:dpp-terminal-shift}.
Now fix an arbitrary \(v\in\Ucal_{t,t+\delta}\).
By \Cref{thm:value-random-initial-data}, for any
\(\varepsilon>0\), there exists a control
\(
    w^{v,\varepsilon}\in\Ucal_{t+\delta,T}
\)
such that
\[
    Y_{t+\delta}^{t+\delta,X_{t+\delta}^{t,x;v};w^{v,\varepsilon}}
    \ge
    u\bigl(t+\delta,X_{t+\delta}^{t,x;v}\bigr)-\varepsilon,
    \qquad
    \Pp\text{-a.s.}
\]
Define 
\(\widehat v_s^{\,v,\varepsilon}:=v_s\,\1_{[t,t+\delta]}(s)
+w_s^{v,\varepsilon}\,\1_{(t+\delta,T]}(s)\).
By the definition of the value function,
\[
\begin{aligned}
    u(t,x)
    \ge J\bigl(t,x;\widehat v^{\,v,\varepsilon}\bigr)
    = G_{t,t+\delta}^{t,x;v}
    \bigl[Y_{t+\delta}^{t+\delta,X_{t+\delta}^{t,x;v};w^{v,\varepsilon}}\bigr]
    \ge G_{t,t+\delta}^{t,x;v}
    \bigl[ u\bigl(t+\delta,X_{t+\delta}^{t,x;v}\bigr)-\varepsilon \bigr]
    \ge G_{t,t+\delta}^{t,x;v}
    \bigl[ u\bigl(t+\delta,X_{t+\delta}^{t,x;v}\bigr) \bigr]-\varepsilon,
    \quad \Pp\text{-a.s.},
\end{aligned}
\]
where the last inequality is due to
\eqref{eq:dpp-terminal-shift}.

Taking the essential supremum over \(v\) and letting \(\varepsilon\downarrow0\), we obtain

\[
    u(t,x)
    \ge \esssup_{v(\cdot)\in\Ucal_{t,t+\delta}}
    G_{t,t+\delta}^{t,x;v}
    \bigl[u\bigl(t+\delta,X_{t+\delta}^{t,x;v}\bigr)\bigr].
\]
This completes the proof.

\end{proof}

Then we show the time continuity of the value function.
\begin{theorem}[Time continuity of the value function]\label{thm:value-time-continuity}
    Suppose that \Cref{ass:standing} holds. Then the value function \(u\) is continuous with respect to the time variable.
\end{theorem}
\begin{proof}
Take \((t,x)\in[0,T)\times\R^{d'}\) and
\(0<\delta\le T-t\). Fix \(\varepsilon>0\).
By the characterization of the deterministic value function in
\Cref{thm:value-function-deterministic}, there exists
\(\widehat v^\varepsilon\in\Ucal^t\) such that
\(
    u(t,x)
    \le
    J\bigl(t,x;\widehat v^\varepsilon\bigr)+\varepsilon.
\)
Set
\(
    v^\varepsilon
    :=
    \widehat v^\varepsilon|_{[t,t+\delta]}
    \in\Ucal_{t,t+\delta},
    \ 
    w^\varepsilon
    :=
    \widehat v^\varepsilon|_{[t+\delta,T]}
    \in\Ucal_{t+\delta,T}.
\)
Thus
\(
    J\bigl(t,x;\widehat v^\varepsilon\bigr)
    =
    G_{t,t+\delta}^{t,x;v^\varepsilon}
    \bigl[
        Y_{t+\delta}^{
            t+\delta,
            X_{t+\delta}^{t,x;v^\varepsilon};
            w^\varepsilon
        }
    \bigr].
\)
By \Cref{thm:value-random-initial-data},
\[
    Y_{t+\delta}^{
        t+\delta,
        X_{t+\delta}^{t,x;v^\varepsilon};
        w^\varepsilon
    }
    \le
    u\bigl(
        t+\delta,
        X_{t+\delta}^{t,x;v^\varepsilon}
    \bigr),
    \qquad
    \Pp\text{-a.s.}
\]
Hence, by the comparison theorem,
\(
    J\bigl(t,x;\widehat v^\varepsilon\bigr)
    \le
    G_{t,t+\delta}^{t,x;v^\varepsilon}
    \bigl[
        u\bigl(
            t+\delta,
            X_{t+\delta}^{t,x;v^\varepsilon}
        \bigr)
    \bigr].
\)
On the other hand, \Cref{thm:dpp} implies
\[
    G_{t,t+\delta}^{t,x;v^\varepsilon}
    \bigl[
        u\bigl(
            t+\delta,
            X_{t+\delta}^{t,x;v^\varepsilon}
        \bigr)
    \bigr]
    \le u(t,x),
    \qquad
    \Pp\text{-a.s.}
\]
Therefore,
\[
\begin{aligned}
    G_{t,t+\delta}^{t,x;v^\varepsilon}
    \bigl[
        u\bigl(
            t+\delta,
            X_{t+\delta}^{t,x;v^\varepsilon}
        \bigr)
    \bigr]
    \le u(t,x)
    \le
    G_{t,t+\delta}^{t,x;v^\varepsilon}
    \bigl[
        u\bigl(
            t+\delta,
            X_{t+\delta}^{t,x;v^\varepsilon}
        \bigr)
    \bigr]
    +\varepsilon,
    \qquad
    \Pp\text{-a.s.}
\end{aligned}
\]
For simplicity, in what follows we write
\(v=v^\varepsilon\). Hence
\begin{equation}\label{eq:time-continuity}
    |u(t,x)-u(t+\delta,x)|
    \le
    |I_\delta^1|+|I_\delta^2|+\varepsilon,
    \qquad
    \Pp\text{-a.s.},
\end{equation}
where
    \(I^1_\delta = G^{t,x;v}_{t,t+\delta}\bigl[u\bigl(t+\delta,X_{t+\delta}^{t,x;v}\bigr)\bigr] - G^{t,x;v}_{t,t+\delta}[u(t+\delta,x)]\), 
    \(I^2_\delta = G^{t,x;v}_{t,t+\delta}[u(t+\delta,x)] - u(t+\delta,x)\).

    For \(I^1_\delta\), by a similar argument as in the proof for \Cref{thm:bsde-uniform-initial-stability}\,\ref{thm:bsde-uniform-initial-stability-ii}, we have
    \begin{align}\label{eq:I1-1}
        \E|I^1_\delta|^2 &\le C\varpi^2 \biggl( \frac{A}{m_1-A} \biggr) + C(m_1)\E|u(t+\delta,X_{t+\delta}^{t,x;v}) - u(t+\delta,x)|^2,
    \end{align}
    for any \(m_1>A\). Also from the proof of \Cref{thm:value-growth-spatial-continuity}, we know that for any \(\varepsilon >0\), there exist \(v_1^\varepsilon (\cdot), v_2^\varepsilon (\cdot) \in \Ucal_{t+\delta,T}\), such that
    \begin{equation}\label{eq:I1-2}
    \begin{aligned}
        &\quad \E|u(t+\delta,X_{t+\delta}^{t,x;v}) - u(t+\delta,x)|^2\\
        &\le C| Y_s^{t+\delta,X^{t,x;v}_{t+\delta};v_1^\varepsilon} - Y_s^{t+\delta,x;v_1^\varepsilon} |^2_{\mathcal{S}^2} + C| Y_s^{t+\delta,X^{t,x;v}_{t+\delta};v_2^\varepsilon} - Y^{t+\delta,x,v_2^\varepsilon}_s |^2_{\mathcal{S}^2} + C\varepsilon^2 \\
        &\le C\varpi^2 \biggl( \frac{A}{m_2-A} \biggr) + C(m_2)\bigl(\E|X_{t+\delta}^{t,x;v} - x|^2+\varpi_\ast\bigl(C\E|X_{t+\delta}^{t,x;v} - x|^2\bigr)\bigr)+ C\varepsilon^2,
    \end{aligned}
    \end{equation}
    
where the last inequality follows from
\eqref{eq:continuous-dependence-final}.
By the a priori estimate for SDEs with jumps, we have
\(
    \E|
        X_{t+\delta}^{t,x;v}-x
    |^2
    \le
    C(1+|x|^2)\delta.
\)
Consequently, by \eqref{eq:I1-1} and \eqref{eq:I1-2}, we have
\[
\begin{aligned}
    \E|I_\delta^1|^2
    \le
    C\varpi^2\biggl(\frac{A}{m_1-A}\biggr)
    +
    C(m_1)\varpi^2\biggl(\frac{A}{m_2-A}\biggr)
    +
    C(m_1)C(m_2)
    [
        (1+|x|^2)\delta
        +
        \varpi_\ast(
            C(1+|x|^2)\delta
        )
    ]
    +
    C(m_1)\varepsilon^2.
\end{aligned}
\]
Take the limit successively in the order: \(\delta\downarrow0\), \(\varepsilon\downarrow0\), \(m_2\uparrow\infty\), \(m_1\uparrow\infty\), and then \(\E|I^1_\delta|^2 \rightarrow 0\) follows.

    The term \(I^2_\delta\) can be estimated directly.
    Applying It\^{o}'s formula to \(\exp(\alpha r)|\tilde{Y}^{t,x;v}_r - u(t+\delta,x)|^2\) on \([t, t+\delta]\), we have
    \begin{equation}\label{eq:time-continuity-2}
    \begin{aligned}
	    &\quad \exp(\alpha t)|\tilde{Y}_{t}^{t,x;v}-u( t+\delta ,x ) |^2+\int_t^{t+\delta}{\alpha \exp(\alpha r)|\tilde{Y}_{r}^{t,x;v}-u( t+\delta ,x ) |^2\dr}\\
        &\quad +\int_t^{t+\delta}{\exp(\alpha r)|\tilde{Z}_{r}^{t,x;v}|^2\dr}+\int_t^{t+\delta}{\int_E{\exp(\alpha r)|\tilde{K}_{r}^{t,x;v}( e ) |^2\mu ( \dr,\de )}}\\
	    &=\int_t^{t+\delta}{2\exp(\alpha r)\bigl( \tilde{Y}_{r}^{t,x;v}-u( t+\delta ,x ) \bigr) f\biggl( r,X_{r}^{t,x;v},\tilde{Y}_{r}^{t,x;v},\tilde{Z}_{r}^{t,x;v},\int_E{\ell_r( e ) \tilde{K}_{r}^{t,x;v}( e )}\lambda(\de) ,v_r \biggr) \dr}\\
	    &\quad -\int_t^{t+\delta}{2\exp(\alpha r)\bigl( \tilde{Y}_{r-}^{t,x;v}-u( t+\delta ,x ) \bigr) \ip[\big]{\tilde{Z}_{r}^{t,x;v}}{\dW_r}}\\
        &\quad -\int_t^{t+\delta}{\int_E{2\exp(\alpha r)\bigl( \tilde{Y}_{r-}^{t,x;v}-u( t+\delta ,x ) \bigr) \tilde{K}_{r}^{t,x;v}( e ) \tilmu( \dr,\de )}},
    \end{aligned}
    \end{equation}
among which we have
    \begin{align*}
        &\quad\bigl( \tilde{Y}_{r}^{t,x;v}-u( t+\delta ,x ) \bigr) f\biggl( r,X_{r}^{t,x;v},\tilde{Y}_{r}^{t,x;v},\tilde{Z}_{r}^{t,x;v},\int_E{\ell_r( e ) \tilde{K}_{r}^{t,x;v}( e )}\lambda(\de) ,v_r \biggr) \\
        &\le \bigl( \tilde{Y}_{r}^{t,x;v}-u( t+\delta ,x ) \bigr) f\biggl( r,X_{r}^{t,x;v},u( t+\delta, x ) ,\tilde{Z}_{r}^{t,x;v},\int_E{\ell_r( e ) \tilde{K}_{r}^{t,x;v}( e ) \lambda(\de)},v_r \biggr) \\
        &\le C| \tilde{Y}_{r}^{t,x;v}-u( t+\delta ,x ) |^2+C(1+|x|^{2p_0})+| X_{r}^{t,x;v} |^2+\frac{1}{2}| \tilde{Z}_{r}^{t,x;v} |^2+\frac{1}{2}\int_E{| \tilde{K}_{r}^{t,x;v}( e ) |^2\lambda(\de)}.
    \end{align*}
    Choosing \(\alpha\) sufficiently large and taking expectations on both sides of \eqref{eq:time-continuity-2}, we have
    \begin{equation*}
    \begin{aligned}
        \E|I^2_\delta|^2 &\le C\delta(1+|x|^{2p_0}) + C\E\int_t^{t+\delta}{|X_{r}^{t,x;v}|^2\dr}\le C\delta (1+|x|^{2p_0}).
    \end{aligned}
    \end{equation*}
    Since \(|u(t,x)-u(t+\delta,x)|\) is deterministic, taking expectations on both sides of \eqref{eq:time-continuity}, we have
    \begin{equation*}
    \begin{aligned}
        |u(t,x)-u(t+\delta,x)| &\le \E|I^1_\delta| + \E|I^2_\delta| + \varepsilon \le \sqrt{\E|I^1_\delta|^2} + \sqrt{\E|I^2_\delta|^2} + \varepsilon.
    \end{aligned}
    \end{equation*}
Taking \(\delta \downarrow 0\) and then \(\varepsilon \downarrow 0\), we obtain the right-continuity of \(u\).
    The left-continuity can be obtained similarly by applying \Cref{thm:dpp} on \([t-\delta,t]\).
\end{proof}

\section{Hamilton-Jacobi-Bellman integro-differential equation}\label{sec:hjb}

In this section, we study the HJB integro-differential equation associated with the stochastic recursive control problem \eqref{eq:value-function}. For this
we define the operators \(\Aop\varphi,\Bop\varphi,\Bop^\delta(\varphi,u),\Cop\varphi\) and \(\Cop^\delta(\varphi,u)\) as follows:
\begin{align*}
    \Aop\varphi(t,x,v)&=\frac{1}{2}\tr\bigl( \sigma \sigma^\top(t,x,v)D^2\varphi (t,x) \bigr) + \langle b(t,x,v),D\varphi (t,x) \rangle, \\
    \Bop\varphi (t,x,v) &= \int_E{( \varphi (t,x+\gamma(t,x,v,e)) - \varphi (t,x) - \ip{D\varphi (t,x)}{\gamma(t,x,v,e)} ) \lambda(\de)}, \\
    \Bop^\delta(\varphi,u) (t,x,v) &= \int_{E_\delta}{( \varphi (t,x+\gamma(t,x,v,e)) - \varphi (t,x) - \ip{D\varphi (t,x)}{\gamma(t,x,v,e)} ) \lambda(\de)} \\
    &\quad + \int_{E_\delta^c}{( u (t,x+\gamma(t,x,v,e)) - u (t,x) - \ip{D\varphi (t,x)}{\gamma(t,x,v,e)} ) \lambda(\de)},\\
    \Cop\varphi (t,x,v) &= \int_E{( \varphi (t,x+\gamma(t,x,v,e)) - \varphi (t,x) )}\ell(t,e)\lambda(\de), \\
    \Cop^\delta(\varphi,u) (t,x,v) &= \int_{E_\delta}{( \varphi (t,x+\gamma(t,x,v,e)) - \varphi (t,x) )}\ell(t,e)\lambda(\de) \\
    &\quad + \int_{E_\delta^c}{( u (t,x+\gamma(t,x,v,e)) - u (t,x) )}\ell(t,e)\lambda(\de),
\end{align*}
where \(E_\delta = \{ e\in E: |e|<\delta \}\) and \(E_\delta^c = E\setminus E_\delta\), \(\varphi\in C^{1,2}([0,T]\times \R^{d'})\). Actually, the HJB equation in our case is a partial integro-differential equation as below:
\begin{equation}\label{eq:hjb}
\left\{
\begin{aligned}
&\frac{\partial}{\partial t}u(t,x)
   \begin{aligned}[t]
+ \underset{v\in U}{\sup}\Bigl\{&\Aop u(t,x,v)+\Bop u(t,x,v) \\
      &+ f(t,x,u(t,x),\sigma^\top(t,x,v)Du(t,x),\Cop u(t,x,v),v)\Bigr\} = 0,\ (t,x)\in [0,T)\times \mathbb{R}^{d'},
   \end{aligned} \\
&u(T,x) = h(x).
\end{aligned}
\right.
\end{equation}

We shall give two equivalent definitions of viscosity solutions to \eqref{eq:hjb}.

\begin{definition}[Global definition]\label{def:viscosity-global}
    A real-valued continuous function \(u\in C([0,T]\times \R^{d'})\) is called a viscosity subsolution (respectively, supersolution) of equation~\eqref{eq:hjb}
    if \(u(T,x)\le h(x)\) (respectively, \(u(T,x)\ge h(x)\)) for any \(x\in \R^{d'}\),
    and if for any test function \(\varphi \in C^{1,2}([0,T]\times \R^{d'})\) and \((t,x)\in [0,T)\times \R^{d'}\)
    such that \(u-\varphi\) attains a global maximum (respectively, minimum) at \((t,x)\), we have
    \begin{align*}
        \frac{\partial}{\partial t}\varphi (t,x) &+ \underset{v\in U}{\sup}\Bigl\{\Aop\varphi (t,x,v)+ \Bop\varphi(t,x,v) \\
        &+ f(t,x,u(t,x),\sigma^\top(t,x,v)D\varphi(t,x),\Cop\varphi(t,x,v),v)\Bigr\} \ge 0\ (\text{respectively,}\le 0).
    \end{align*}
\end{definition}

Next, we introduce the parabolic neighborhood of \((t,x)\) as \(B((t,x),R) := (t-R^2,t+R^2)\times \{y\in \R^{d'}: |y-x|<R\}\).
For any fixed pair \((\delta,\Lambda)\) with \(\delta>0\) and \(\Lambda>L_2\), we give the following local definition of viscosity solutions.
\begin{definition}[Local definition]\label{def:viscosity-local-i}
    A real-valued continuous function \(u\in C([0,T]\times \R^{d'})\) is called a viscosity subsolution (respectively, supersolution) of equation~\eqref{eq:hjb}
    if \(u(T,x)\le h(x)\) (respectively, \(u(T,x)\ge h(x)\)) for any \(x\in \R^{d'}\),
    and if for any test function \(\varphi \in C^{1,2}([0,T]\times \R^{d'})\) and \((t,x)\in [0,T)\times \R^{d'}\)
    such that \(u-\varphi\) attains a maximum (respectively, minimum) in \(B((t,x),R_{\delta,x})\) at \((t,x)\), we have
    \begin{align*}
        \frac{\partial}{\partial t}\varphi (t,x) &+ \underset{v\in U}{\sup}\Bigl\{\Aop\varphi (t,x,v) + \Bop^\delta (\varphi,u)(t,x,v) \\
        &+ f(t,x,u(t,x),\sigma^\top(t,x,v)D\varphi(t,x),\Cop^\delta (\varphi,u)(t,x,v),v)\Bigr\} \ge 0\ (\text{respectively,}\le 0),
    \end{align*}
    where \(R_{\delta,x} := \Lambda(1+|x|)\delta\).
\end{definition}

A viscosity solution is a function which is both a viscosity subsolution and a viscosity supersolution.

We next introduce the equivalence between the global and local definitions. 
The argument in \cite[Lemma~4.1]{LiPeng2009} applies with slight modification, and hence we omit the details here.

\begin{lemma}[Equivalence of the viscosity formulations]
\label{lem:viscosity-definitions-equivalent}
The global formulation in \Cref{def:viscosity-global} and the localized
formulation in \Cref{def:viscosity-local-i} are equivalent.
\end{lemma}

Next, we introduce some notation. For \(l>0\), define the following truncated functions:
\begin{align*}
    \Pi_l(y)=(|y|\wedge l)\sgn(y),\ y\in \R,\ \ \ \text{and}\ \ \
    f^l(t,x,y,z,k,v) = f(t,x,\Pi_l(y),z,k,v).
\end{align*}

For \(m>A\), define
\begin{align*}
f^{m,l}(t,x,y,z,k,v):=
\inf_{\substack{(u_0,u_1,u_2,u_3,u) \in\\
\R^{d'} \times \R \times \R^d \times \R \times U}}\Bigl\{&f^l(t,u_0,u_1,u_2,u_3,u)\\
&+m(|x-u_0|+|y-u_1|+|z-u_2|+|k-u_3|+|v-u|)\Bigr\}.
\end{align*}

We introduce the following lemma, the proof of which is similar to 
that of \cite[Lemma~1]{LepeltierSanMartin1997} and is omitted here.

\begin{lemma}\label{lem:full-inf-convolution}
Fix \(l>0\). For \(m>A\), let \(f^{m,l}\) be defined as above.
Then the following properties hold.

\begin{enumerate}[label=\textup{(\roman*)},ref=\textup{(\roman*)},wide=0pt]

\item\label{lem:full-inf-convolution-i}
The function \(f^{m,l}\) is finite and continuous on
\(
[0,T]\times\R^{d'}\times\R\times\R^d\times\R\times U.
\)
Moreover, for any
\(
q=(x,y,z,k,v),\ 
q'=(x',y',z',k',v')
\)
in
\(\R^{d'}\times\R\times\R^d\times\R\times U\),
\[
\begin{aligned}
    |f^{m,l}(t,x,y,z,k,v)-f^{m,l}(t,x',y',z',k',v')|&\le
    m(|x-x'|+|y-y'|+|z-z'|+|k-k'|+|v-v'|),\\
    |f^{m,l}(t,x,y,z,k,v)|&\le C(1+|x|+l^{p_0}+|z|+|k|).
\end{aligned}
\]

\item
For any \((t,x,z,k,v)\in[0,T]\times\R^{d'}\times\R^d\times\R\times U\), the $\R\to\R$ mapping
\(
y\longmapsto f^{m,l}(t,x,y,z,k,v)
\)
is nonincreasing. For any \((t,x,y,z,v)\in[0,T]\times\R^{d'}\times\R\times\R^d\times U\), the $\R\to\R$ mapping
\(
k\longmapsto f^{m,l}(t,x,y,z,k,v)
\)
is nondecreasing.

\item
For any
\((t,x,y,z,k,v)\in[0,T]\times\R^{d'}\times\R\times\R^d\times\R\times U\), the sequence
\(\{f^{m,l}(t,x,y,z,k,v)\}_{m>A}\) is nondecreasing in \(m\),
and
\[
f^{m,l}(t,x,y,z,k,v)
\le
f^l(t,x,y,z,k,v).
\]

\item
As \(m\to\infty\),
\(
f^{m,l}\uparrow f^l
\)
uniformly on every compact subset of
\(
[0,T]\times\R^{d'}\times\R\times\R^d\times\R\times U.
\)

\end{enumerate}
\end{lemma}

Then we can define two families of backward systems:
\begin{equation*}
\begin{aligned}
    Y^{m,l;t,x;v}_s &= h\bigl(X_{T}^{t,x;v}\bigr) + \int_s^T{f^{m,l}\biggl(r,X_{r}^{t,x;v},Y_{r}^{m,l;t,x;v},Z_{r}^{m,l;t,x;v},\int_E{\ell_r(e)K_{r}^{m,l;t,x;v}(e)}\lambda(\de),v_r\biggr)\dr} \\
    &\quad - \int_s^T{\ip[\big]{Z_{r}^{m,l;t,x;v}}{\dW_r}} - \int_s^T{\int_E{K_{r}^{m,l;t,x;v}(e)\tilmu(\dr,\de)}},\ s\in [t,T],
\end{aligned}
\end{equation*}
and
\begin{equation*}
\begin{aligned}
    Y^{l;t,x;v}_s &= h\bigl(X_{T}^{t,x;v}\bigr) + \int_s^T{f^l\biggl(r,X_{r}^{t,x;v},Y_{r}^{l;t,x;v},Z_{r}^{l;t,x;v},\int_E{\ell_r(e)K_{r}^{l;t,x;v}(e)}\lambda(\de),v_r\biggr)\dr} \\
    &\quad - \int_s^T{\ip[\big]{Z_{r}^{l;t,x;v}}{\dW_r}} - \int_s^T{\int_E{K_{r}^{l;t,x;v}(e)\tilmu(\dr,\de)}},\ s\in [t,T].
\end{aligned}
\end{equation*}

For \((t,x)\in[0,T]\times\R^{d'}\), define the corresponding value functions by
\[
u^{m,l}(t,x):=\esssup_{v\in\Ucal_{t,T}}Y_t^{m,l;t,x;v},
\qquad
u^l(t,x):=\esssup_{v\in\Ucal_{t,T}}Y_t^{l;t,x;v}.
\]

\begin{theorem}[Monotone convergence of truncated BSDEs]\label{thm:truncated-bsde-monotone-convergence}
    Suppose that \Cref{ass:standing} holds. For any \((t,x) \in [0,T]\times \R^{d'}\), \(v \in \Ucal_{t,T}\) and \(l>0\), we have
    \(Y^{m,l;t,x;v} \uparrow Y^{l;t,x;v}\) in \(\mathcal{S}^2\), \(Z^{m,l;t,x;v} \rightarrow Z^{l;t,x;v}\) in \(\mathcal{L}^2\) and \(K^{m,l;t,x;v} \rightarrow K^{l;t,x;v}\) in \(\mathcal{L}^2_\lambda\), as \(m \rightarrow \infty\).
\end{theorem}

\begin{proof}

The truncated generators \(f^l\) inherit the same modulus \(\varpi\) in
both the state and control variables, as well as the same uniform-continuity
constants in \((z,k)\), independently of \(l\). In particular, the modulus
\(\Omega_m\) in the uniform initial-stability estimate can be chosen
independently of \(l\).

    Since \(f^{m,l} \uparrow f^l\), from \Cref{lem:bsde-existence-comparison}\,\ref{lem:bsde-existence-comparison-i}, we can deduce that for \(m>m_0:=\lfloor2A+1\rfloor\), \(Y^{m_0,l;t,x;v}_s \le Y^{m,l;t,x;v}_s \le Y^{l;t,x;v}_s\), \(\Pp\)-a.s.
    and \(\lVert\sup_m |Y^{m,l;t,x;v}_s| \rVert_{\mathcal{S}^2}^2 \le \lVert Y^{l;t,x;v}_s \rVert_{\mathcal{S}^2}^2 + \lVert Y^{m_0,l;t,x;v}_s \rVert_{\mathcal{S}^2}^2 < C(1+|x|^2+l^{2p_0})\).
Applying It\^{o}'s formula to \(|Y^{m,l;t,x;v}_r|^2\) on \([s,T]\), we have
    \begin{equation*}
    \begin{aligned}
        &\quad |Y^{m,l;t,x;v}_s|^2 + \int_s^T{|Z_{r}^{m,l;t,x;v}|^2\dr} + \int_s^T{\int_E{|K_{r}^{m,l;t,x;v}(e)|^2\mu(\dr,\de)}} \\
        &= 2\int_s^T{Y_{r}^{m,l;t,x;v}f^{m,l}\biggl(r,X_{r}^{t,x;v},Y_{r}^{m,l;t,x;v},Z_{r}^{m,l;t,x;v},\int_E{\ell_r(e)K_{r}^{m,l;t,x;v}(e)}\lambda(\de),v_r\biggr)\dr} \\
        &\quad +|h(X_{T}^{t,x;v})|^2 - 2\int_s^T{Y_{r-}^{m,l;t,x;v}\ip[\big]{Z_{r}^{m,l;t,x;v}}{\dW_r}} - 2\int_s^T{\int_E{Y_{r-}^{m,l;t,x;v}K_{r}^{m,l;t,x;v}(e)\tilmu(\dr,\de)}}.
    \end{aligned}
    \end{equation*}
Taking expectations on both sides, by \Cref{lem:full-inf-convolution}\,\ref{lem:full-inf-convolution-i}, we have
    \begin{equation*}
    \begin{aligned}
        &\quad \E|Y^{m,l;t,x;v}_s|^2 + \E\int_s^T{|Z_{r}^{m,l;t,x;v}|^2\dr} + \E\int_s^T{\int_E{|K_{r}^{m,l;t,x;v}(e)|^2\lambda(\de)\dr}} \\
        &\le \E|h( X_{T}^{t,x;v} ) |^2+C\E\int_s^T{|Y_{r}^{m,l;t,x;v}|}
        \biggl( 1+|X_{r}^{t,x;v}|+|l|^{p_0}+|Z_{r}^{m,l;t,x;v}|+\biggl| \int_E{\ell_r( e ) K_{r}^{m,l;t,x;v}( e )}\lambda(\de) \biggr| \biggr) \dr \\
        &\le C\bigl(1+l^{2p_0}+\lVert X^{t,x;v}_s \rVert_{\mathcal{S}^2}^{2}+\lVert Y^{m,l;t,x;v}_s \rVert_{\mathcal{S}^2}^{2}\bigr) + \frac{1}{2}\E\int_s^T{|Z_{r}^{m,l;t,x;v}|^2\dr} 
        + \frac{1}{2}\E\int_s^T{\int_E{|K_{r}^{m,l;t,x;v}(e)|^2\lambda(\de)\dr}},
    \end{aligned}
    \end{equation*}
    which, together with the prior uniform estimate of \(Y^{m,l;t,x;v}_s\), leads to
    \begin{equation*}
    \begin{aligned}
        \lVert Y^{m,l;t,x;v}_s \rVert_{\mathcal{S}^2}^{2} + \lVert Z^{m,l;t,x;v}_s \rVert_{\mathcal{L}^2}^{2} + \lVert K^{m,l;t,x;v}_s \rVert_{\mathcal{L}^2_\lambda}^{2}
        \le C\bigl(1+l^{2p_0}+\lVert X^{t,x;v}_s \rVert_{\mathcal{S}^2}^{2}+\lVert Y^{m,l;t,x;v}_s \rVert_{\mathcal{S}^2}^{2}\bigr) < C(1+l^{2p_0}+|x|^2).
    \end{aligned}
    \end{equation*}

Since \(Y^{m,l;t,x;v}_s\) is increasing and has an \(\mathcal{S}^2\) upper bound \(Y^{l;t,x;v}_s\), we know that \(Y^{m,l;t,x;v}_s\) has a pathwise limit, denoted by \(\tilde{Y}^{l;t,x;v}_s\).
    By the dominated convergence theorem, we have \(Y^{m,l;t,x;v}_s \uparrow \tilde{Y}^{l;t,x;v}_s\) in \(\mathcal{L}^2\).

    We set
    \begin{align*}
        \check{Y}^{m,n}_s &= Y^{m,l;t,x;v}_s - Y^{n,l;t,x;v}_s, \
        \check{Z}^{m,n}_s = Z^{m,l;t,x;v}_s - Z^{n,l;t,x;v}_s, \
        \check{K}^{m,n}_s(e) = K^{m,l;t,x;v}_s(e) - K^{n,l;t,x;v}_s(e),
    \end{align*}
    and
    \begin{align*}
        \check{f}^{m,n}_s &= f^{m,l}\biggl(s,X_{s}^{t,x;v},Y_{s}^{m,l;t,x;v},Z_{s}^{m,l;t,x;v},\int_E{\ell_s(e)K_{s}^{m,l;t,x;v}(e)}\lambda(\de),v_s\biggr) \\
        &\quad - f^{n,l}\biggl(s,X_{s}^{t,x;v},Y_{s}^{n,l;t,x;v},Z_{s}^{n,l;t,x;v},\int_E{\ell_s(e)K_{s}^{n,l;t,x;v}(e)}\lambda(\de),v_s\biggr).
    \end{align*}
Applying It\^{o}'s formula to \(|\check{Y}^{m,n}_r|^2\) on \([t,T]\)
    and taking expectations on both sides, we have
    \begin{equation*}
    \begin{aligned}
	    &\quad \E\int_t^T{|\check{Z}^{m,n}_r|^2\dr}+\E\int_t^T{\int_E{|\check{K}^{m,n}_r(e)|^2\lambda(\de) \dr}}\\
	    &\le 2\E\int_t^T{ \check{Y}^{m,n}_r  \check{f}^{m,n}_r \dr}\le 2\biggl( \E\int_t^T{|\check{Y}^{m,n}_r|^2\dr} \biggr) ^{1/2}\biggl( \E\int_t^T{| \check{f}^{m,n}_r |^2\dr} \biggr) ^{1/2}.
    \end{aligned}
    \end{equation*}
Moreover,
    \begin{equation*}
    \begin{aligned}
        &\quad \biggl\| f^{m,l}\biggl(r,X_{r}^{t,x;v},Y_{r}^{m,l;t,x;v},Z_{r}^{m,l;t,x;v},\int_E{\ell_r( e ) K_{r}^{m,l;t,x;v}( e )}\lambda(\de) ,v_r \biggr) \biggr\|^2_{\mathcal{L}^2} \\
        &\le C\bigl( 1+\lVert X^{t,x;v}_r \rVert_{\mathcal{S}^2}^{2}+\lvert l \rvert^{2p_0}+\lVert Z^{m,l;t,x;v}_r \rVert_{\mathcal{L}^2}^{2}+\lVert K^{m,l;t,x;v}_r \rVert_{\mathcal{L}^2_\lambda}^{2} \bigr)\le C(1+|x|^2+l^{2p_0}).
    \end{aligned}
    \end{equation*}
    Hence we have
    \begin{align*}
        \lVert \check{Z}^{m,n}_s \rVert_{\mathcal{L}^2}^{2} + \lVert \check{K}^{m,n}_s \rVert_{\mathcal{L}^2_\lambda}^{2} &\le C(1+|x|^2+l^{2p_0})\lVert \check{Y}^{m,n}_s \rVert_{\mathcal{L}^2},
    \end{align*}
    which means that \(Z^{m,l;t,x;v}_s\) and \(K^{m,l;t,x;v}_s(e)\) are Cauchy sequences in \(\mathcal{L}^2\) and \(\mathcal{L}^2_\lambda\), respectively, with their limits denoted by \(\tilde{Z}^{l;t,x;v}_s\) and \(\tilde{K}^{l;t,x;v}_s(e)\).

    Finally, we need to show that \(\tilde{Y}^{l;t,x;v}_s = Y^{l;t,x;v}_s\), \(\tilde{Z}^{l;t,x;v}_s = Z^{l;t,x;v}_s\), \(\tilde{K}^{l;t,x;v}_s = K^{l;t,x;v}_s\) and \(Y^{m,l;t,x;v}_s \uparrow Y^{l;t,x;v}_s\) in \(\mathcal{S}^2\).
    Without loss of generality, we have
    \[\biggl(Y^{m,l;t,x;v}_s,\ Z^{m,l;t,x;v}_s,\ \int_{E}{\ell_s(e)K^{m,l;t,x;v}_s(e)}\lambda(\de)\biggr) \rightarrow \biggl(\tilde{Y}^{l;t,x;v}_s,\ \tilde{Z}^{l;t,x;v}_s,\ \int_{E}{\ell_s(e)\tilde{K}^{l;t,x;v}_s(e)}\lambda(\de)\biggr),\]
    \((\ds\otimes\Pp)\)-a.e. and in \(\mathcal{L}^2\).
    Taking a subsequence if necessary, we have
    \[\sup_m |Y^{m,l;t,x;v}_s|,\  \sup_m |Z^{m,l;t,x;v}_s|,\  \sup_m \biggl|\int_{E}{\ell_s(e)K^{m,l;t,x;v}_s(e)}\lambda(\de)\biggr| \in L^1([t,T]\times\Omega,\ds\otimes\Pp).\]
    On one hand, we can deduce the following estimate:
    \begin{equation*}
    \begin{aligned}
        &\quad \biggl| f^{m,l}\biggl(r,X_{r}^{t,x;v},Y_{r}^{m,l;t,x;v},Z_{r}^{m,l;t,x;v},\int_E{\ell_r(e)K_{r}^{m,l;t,x;v}(e)}\lambda(\de),v_r\biggr)\biggr|\\
        &\le C\biggl( 1+|X_{r}^{t,x;v}|+|l|^{p_0}+\sup_m |Z_{r}^{m,l;t,x;v}|+\sup_m \biggl|\int_E{\ell_r(e)K_{r}^{m,l;t,x;v}(e)}\lambda(\de)\biggr| \biggr),
    \end{aligned}
    \end{equation*}
    which belongs to \(L^1([t,T]\times\Omega,\ds\otimes\Pp)\). On the other hand, we have the following pathwise convergence:
    \begin{equation*}
    \begin{aligned}
        &\quad f^{m,l}\biggl(r,X_{r}^{t,x;v},Y_{r}^{m,l;t,x;v},Z_{r}^{m,l;t,x;v},\int_E{\ell_r(e)K_{r}^{m,l;t,x;v}(e)}\lambda(\de),v_r\biggr) \\
        &\rightarrow f^l\biggl(r,X_{r}^{t,x;v},\tilde{Y}_{r}^{l;t,x;v},\tilde{Z}_{r}^{l;t,x;v},\int_E{\ell_r(e)\tilde{K}_{r}^{l;t,x;v}(e)}\lambda(\de),v_r\biggr),\ \ (\dr\otimes\Pp)\text{-a.e.}
    \end{aligned}
    \end{equation*}
    By the dominated convergence theorem, the above convergence also holds in \(L^1([t,T]\times\Omega,\ds\otimes\Pp)\). Passing to a subsequence, we have the following pathwise convergence:
    \begin{equation*}
    \begin{aligned}
        &\quad \int_s^T{f^{m,l}\biggl(r,X_{r}^{t,x;v},Y_{r}^{m,l;t,x;v},Z_{r}^{m,l;t,x;v},\int_E{\ell_r(e)K_{r}^{m,l;t,x;v}(e)}\lambda(\de),v_r\biggr)\dr} \\
        &\rightarrow \int_s^T{f^l\biggl(r,X_{r}^{t,x;v},\tilde{Y}_{r}^{l;t,x;v},\tilde{Z}_{r}^{l;t,x;v},\int_E{\ell_r(e)\tilde{K}_{r}^{l;t,x;v}(e)}\lambda(\de),v_r\biggr)\dr},\quad\text{for all }s \in [t,T],\ \Pp\text{-a.s.}
    \end{aligned}
    \end{equation*}
    On the other hand, by the BDG inequality, we have
    \begin{align*}
        &\quad \E\Biggl[ \sup_{t\le s\le T}\biggl| \int_s^T{\ip[\big]{Z_{r}^{m,l;t,x;v}}{\dW_r}} - \int_s^T{\ip[\big]{\tilde{Z}_{r}^{l;t,x;v}}{\dW_r}} \biggr|^2 \Biggr]
        \le C\E\int_t^T{|Z_{r}^{m,l;t,x;v} - \tilde{Z}_{r}^{l;t,x;v}|^2\dr},\\
        &\quad \E\Biggl[ \sup_{t\le s\le T}\biggl| \int_s^T{\int_E{K_{r}^{m,l;t,x;v}(e)\tilmu(\dr,\de)}} - \int_s^T{\int_E{\tilde{K}_{r}^{l;t,x;v}(e)\tilmu(\dr,\de)}} \biggr|^2 \Biggr] \\
        &\le C\E\int_t^T{\int_E{|K_{r}^{m,l;t,x;v}(e) - \tilde{K}_{r}^{l;t,x;v}(e)|^2\lambda(\de)\dr}}.
    \end{align*}

Consequently, we have the desired pathwise convergence after passing to a subsequence:
    \begin{equation*}
    \begin{aligned}
        &\int_s^T{\ip[\big]{Z_{r}^{m,l;t,x;v}}{\dW_r}} \rightarrow \int_s^T{\ip[\big]{\tilde{Z}_{r}^{l;t,x;v}}{\dW_r}},\\
        &\int_s^T{\int_E{K_{r}^{m,l;t,x;v}(e)\tilmu(\dr,\de)}} \rightarrow \int_s^T{\int_E{\tilde{K}_{r}^{l;t,x;v}(e)\tilmu(\dr,\de)}},\quad\text{for all }s\in [t,T],\ \Pp\text{-a.s.}
    \end{aligned}
    \end{equation*}
    Passing to the limit in the BSDE that \(Y^{m,l;t,x;v}_s\) satisfies in the sense of \(\Pp\text{-a.s.}\),
    we obtain the conclusion from the uniqueness of the solution of the BSDE with jumps.

    Finally, for the \(\mathcal{S}^2\) convergence of \(Y^{m,l;t,x;v}_s\), we have
    \begin{equation*}
    \begin{aligned}
	    &\quad \sup_{t\le s\le T}|Y_{s}^{m,l;t,x;v}-Y_{s}^{l;t,x;v}|\\
        &\le \int_t^T{\biggl| f^{m,l}\biggl( r,X_{r}^{t,x;v},Y_{r}^{m,l;t,x;v},Z_{r}^{m,l;t,x;v},\int_E{\ell_r( e ) K_{r}^{m,l;t,x;v}( e )}\lambda(\de) ,v_r \biggr) }\\
	    &\quad \quad \quad  -f^l\biggl( r,X_{r}^{t,x;v},Y_{r}^{l;t,x;v},Z_{r}^{l;t,x;v},\int_E{\ell_r( e ) K_{r}^{l;t,x;v}( e )}\lambda(\de) ,v_r \biggr) \biggr|\dr\\
	    &\quad +\sup_{t\le s\le T}\biggl| \int_s^T{\ip[\big]{\bigl( Z_{r}^{m,l;t,x;v}-Z_{r}^{l;t,x;v} \bigr)}{\dW_r}} \biggr|+\sup_{t\le s\le T}\biggl| \int_s^T{\int_E{\bigl( K_{r}^{m,l;t,x;v}( e ) -K_{r}^{l;t,x;v}( e ) \bigr) \tilmu( \dr,\de )}} \biggr|\\
        &\rightarrow 0,\ \Pp\text{-a.s.}
    \end{aligned}
    \end{equation*}
Note that we also have \(\sup_m |Y_{s}^{m,l;t,x;v}| \in \mathcal{S}^2\). By the dominated convergence theorem, we have \(Y^{m,l;t,x;v}_s \rightarrow Y^{l;t,x;v}_s\) in \(\mathcal{S}^2\).

Since \(Y^{m,l;t,x;v}\uparrow Y^{l;t,x;v}\), the \(\mathcal{S}^2\)
convergence obtained along the subsequence immediately extends
to the whole sequence. Thus, the proof is complete.

\end{proof}

\begin{lemma}[Uniform convergence of truncated BSDEs]\label{lem:truncated-bsde-uniform-convergence}
    Suppose that \Cref{ass:standing} holds. For any compact subset \(K\) in \([0,T]\times \R^{d'}\), we have
    \begin{align*}
        \underset{(t,x)\in K}{\sup}\esssup_{v\in \Ucal_{t,T}}\E|Y^{l;t,x;v}_t - Y^{t,x;v}_t| \rightarrow 0,\quad\text{as }l \rightarrow \infty.
    \end{align*}
\end{lemma}

\begin{proof}
    Since \(K\) is compact, its projection in \(\R^{d'}\) can be covered by a finite \(\delta\)-net \(\{x_1,x_2,\cdots,x_N\}\).
    Therefore, for any \((t,x)\in K\), there exists \(x_i\), \(i\in\{1,2,\cdots,N\}\), such that \(|x-x_i|<\delta\). So, we have
    \allowdisplaybreaks[2]
    \begin{align*}
        &\quad \sup_{0\le t\le T}\esssup_{v\in \Ucal_{t,T}}\E| Y_{t}^{l;t,x;v}-Y_{t}^{t,x;v} |\\
        &\le \sup_{0\le t\le T}\esssup_{v\in \Ucal_{t,T}}\bigl( \E| Y_{t}^{l;t,x;v}-Y_{t}^{l;t,x_i;v} |+\E| Y_{t}^{t,x_i;v}-Y_{t}^{t,x;v} | \bigr) +\sup_{0\le t\le T}\esssup_{v\in \Ucal_{t,T}}\E| Y_{t}^{l;t,x_i;v}-Y_{t}^{t,x_i;v} |\\
        &\le \sup_{0\le t\le T}\Bigl( \sup_{| x-y |\le \delta}\esssup_{v\in \Ucal_{t,T}}\E| Y_{t}^{l;t,x;v}-Y_{t}^{l;t,y;v} |+\sup_{| x-y |\le \delta}\esssup_{v\in \Ucal_{t,T}}\E| Y_{t}^{t,x;v}-Y_{t}^{t,y;v} | \Bigr)
        +\sup_{0\le t\le T}\esssup_{v\in \Ucal_{t,T}}\E| Y_{t}^{l;t,x_i;v}-Y_{t}^{t,x_i;v} |\\
        &\le \sup_{| x-y |\le \delta}\Bigl( \sup_{0\le t\le T}\esssup_{v\in \Ucal_{t,T}}\E| Y_{t}^{l;t,x;v}-Y_{t}^{l;t,y;v} |+\sup_{0\le t\le T}\esssup_{v\in \Ucal_{t,T}}\E| Y_{t}^{t,x;v}-Y_{t}^{t,y;v} | \Bigr) \\
        &\quad +N\sup_{1\le i\le N}\sup_{0\le t\le T}\esssup_{v\in \Ucal_{t,T}}\E| Y_{t}^{l;t,x_i;v}-Y_{t}^{t,x_i;v} |.
    \end{align*}
Denote the first term on the right-hand side of the last inequality by \(I\) and set
\[
II:=\sup_{1\le i\le N}\sup_{0\le t\le T}\esssup_{v\in\Ucal_{t,T}}\E|Y_t^{l;t,x_i;v}-Y_t^{t,x_i;v}|.
\]
Thus the right-hand side is \(I+N\cdot II\).
    By \Cref{thm:bsde-uniform-initial-stability}\,\ref{thm:bsde-uniform-initial-stability-i}, we have that for any \(\varepsilon >0\), there exists \(\delta >0\) independent of \(l\) such that \(I<\varepsilon \).

    To estimate \(II\), set
    \begin{align*}
        \bar{Y}^l_s := Y^{l;t,x;v}_s - Y^{t,x;v}_s,\ \bar{Z}^l_s := Z^{l;t,x;v}_s - Z^{t,x;v}_s,\ \bar{K}^l_s(e) := K^{l;t,x;v}_s(e) - K^{t,x;v}_s(e).
    \end{align*}
    Then applying It\^{o}'s formula to \(\exp(\alpha r)|\bar{Y}^l_r|^2\) on \([s,T]\), we have
    \begin{equation*}
    \begin{aligned}
        &\quad \exp(\alpha s)|\bar{Y}^l_s|^2 +\int_s^T{\alpha \exp(\alpha r)|\bar{Y}^l_r|^2\dr} +\int_s^T{\exp(\alpha r)|\bar{Z}^l_r|^2\dr} + \int_s^T{\int_E{\exp(\alpha r)|\bar{K}^l_r(e)|^2\mu(\dr,\de)}} \\
        &=2\int_s^T{\exp(\alpha r)\bar{Y}^l_r\Biggl( f^l\biggl(r,X_{r}^{t,x;v},Y_{r}^{l;t,x;v},Z_{r}^{l;t,x;v},\int_E{\ell_r(e)K_{r}^{l;t,x;v}(e)}\lambda(\de),v_r\biggr) }\\
        &\quad  - f\biggl(r,X_{r}^{t,x;v},Y_{r}^{t,x;v},Z_{r}^{t,x;v},\int_E{\ell_r(e)K_{r}^{t,x;v}(e)}\lambda(\de),v_r\biggr) \Biggr) \dr - 2\int_s^T{\exp(\alpha r)\bar{Y}^l_{r-}\ip[\big]{\bar{Z}^l_r}{\dW_r}} \\
        &\quad - 2\int_s^T{\int_E{\exp(\alpha r)\bar{Y}^l_{r-}\bar{K}^l_r(e)\tilmu(\dr,\de)}}.
    \end{aligned}
    \end{equation*}
    Meanwhile, since \(f\) is monotone in \(y\) and uniformly continuous with respect to \(z,k\), we have the following estimate:
    \begin{equation*}
    \begin{aligned}
        &\quad \bar{Y}^l_r\bigl(
        f^l\bigl(r,X_{r}^{t,x;v},Y_{r}^{l;t,x;v},Z_{r}^{l;t,x;v},\widehat{K}_{r}^{l;t,x;v},v_r\bigr)
        -f\bigl(r,X_{r}^{t,x;v},Y_{r}^{t,x;v},Z_{r}^{t,x;v},\widehat{K}_{r}^{t,x;v},v_r\bigr)
        \bigr) \\
        &\le \bar{Y}^l_r\bigl(
        f\bigl(r,X_{r}^{t,x;v},\Pi_l\bigl(Y_{r}^{l;t,x;v}\bigr),Z_{r}^{l;t,x;v},\widehat{K}_{r}^{l;t,x;v},v_r\bigr)
        -f\bigl(r,X_{r}^{t,x;v},\Pi_l\bigl(Y_{r}^{l;t,x;v}\bigr),Z_{r}^{t,x;v},\widehat{K}_{r}^{t,x;v},v_r\bigr)
        \bigr) \\
        &\quad + \bigl[\bigl(Y_{r}^{l;t,x;v} - \Pi_l\bigl(Y_{r}^{l;t,x;v}\bigr)\bigr)
        +\bigl(\Pi_l\bigl(Y_{r}^{l;t,x;v}\bigr) - Y_{r}^{t,x;v}\bigr)\bigr] \\
        &\qquad\times \bigl(
        f\bigl(r,X_{r}^{t,x;v},\Pi_l\bigl(Y_{r}^{l;t,x;v}\bigr),Z_{r}^{t,x;v},\widehat{K}_{r}^{t,x;v},v_r\bigr)
        -f\bigl(r,X_{r}^{t,x;v},Y_{r}^{t,x;v},Z_{r}^{t,x;v},\widehat{K}_{r}^{t,x;v},v_r\bigr)
        \bigr) \\
        &\le \rho|\bar{Y}^l_r||\bar{Z}^l_r|
        +\rho \1_{\{\bar{Y}^l_r \int_E{\ell_r(e)\bar{K}^l_r(e)}\lambda(\de)>0\}}
        \bar{Y}^l_r \int_E{\ell_r(e)\bar{K}^l_r(e)}\lambda(\de)  + \Delta_l(r)+2\varpi\Bigl(\frac{2A}{\rho}\Bigr)|\bar{Y}^l_r|,
    \end{aligned}
    \end{equation*}
    where
    \[
        \widehat{K}_{r}^{l;t,x;v}
        :=\int_E \ell_r(e)K_{r}^{l;t,x;v}(e)\lambda(\de),
        \quad
        \widehat{K}_{r}^{t,x;v}
        :=\int_E \ell_r(e)K_{r}^{t,x;v}(e)\lambda(\de),
    \]
    and
    \begin{align*}
        \Delta_l(r) &:= \bigl[Y_r^{l;t,x;v} - \Pi_l\bigl(Y_r^{l;t,x;v}\bigr)\bigr] \Biggl[f\biggl(r,X_{r}^{t,x;v},\Pi_l\bigl(Y_{r}^{l;t,x;v}\bigr),Z_{r}^{t,x;v},\int_E{\ell_r(e)K_{r}^{t,x;v}(e)}\lambda(\de),v_r\biggr)  \\
        &\qquad  - f\biggl(r,X_{r}^{t,x;v},Y_{r}^{t,x;v},Z_{r}^{t,x;v},\int_E{\ell_r(e)K_{r}^{t,x;v}(e)}\lambda(\de),v_r\biggr)\Biggr].
    \end{align*}
Taking \(s=t\), we have
    \begin{equation*}
    \begin{aligned}
        &\quad \exp(\alpha t)|\bar{Y}^l_t|^2 +\int_t^T{\alpha \exp(\alpha r)|\bar{Y}^l_r|^2\dr} +\int_t^T{\exp(\alpha r)|\bar{Z}^l_r|^2\dr} + \int_t^T{\int_E{\exp(\alpha r)|\bar{K}^l_r(e)|^2\mu(\dr,\de)}} \\
        &\le C\int_t^T{\exp(\alpha r)\Bigl(\varpi^2\Bigl(\frac{2A}{\rho}\Bigr)+|\bar{Y}^l_r|^2 \Bigr)\dr} +2\int_t^T{\exp(\alpha r)\Delta_l(r)\dr} \\
        &\quad - 2\int_t^T{\exp(\alpha r)\bar{Y}^l_{r-}\ip[\bigg]{\bar{Z}^l_r}{\dW_r - \rho \1_{\{\bar{Y}^l_{r-} \bar{Z}^l_r \neq 0 \}} \frac{\bar{Y}^l_{r-} \bar{Z}^l_r}{|\bar{Y}^l_{r-}||\bar{Z}^l_{r}|}\dr}} \\
        &\quad - 2\int_t^T{\int_E{\exp(\alpha r)\bar{Y}^l_{r-}\bar{K}^l_r(e)\bigl(\tilmu(\dr,\de) - \rho \1_{\{\bar{Y}^l_{r-} \int_E{\ell_r(e)\bar{K}^l_r(e)}\lambda(\de)>0\}} \ell_r(e)\lambda(\de)\dr\bigr)} }.
    \end{aligned}
    \end{equation*}
    Choose \(\alpha = C+1\). Using a similar Girsanov transformation argument as in \Cref{lem:bsde-existence-comparison}\,\ref{lem:bsde-existence-comparison-i}, taking the conditional expectation under the new probability measure \(\Pp^\rho\) at time \(t\) first, and then taking the expectation under \(\Pp\), we have
    \begin{equation}\label{eq:truncated-bsde-convergence-1}
    \begin{aligned}
        \E|\bar{Y}^l_t|^2 \le C\varpi^2\Bigl(\frac{2A}{\rho}\Bigr) + C\E\biggl(\frac{\Theta_T}{\Theta_t}\int_t^T{|\Delta_l(r)|\dr}\biggr),
    \end{aligned}
    \end{equation}
where
\[
\begin{gathered}
\Theta_s=\mathcal E\biggl(
\int_t^\cdot \rho\1_{\{\bar Y^l_{r-}\bar Z^l_r\ne0\}}
\ip[\bigg]{\frac{\bar Y^l_{r-}\bar Z^l_r}{|\bar Y^l_{r-}||\bar Z^l_r|}}{\dW_r}
+\int_t^\cdot\int_E
\rho\1_{\{\bar Y^l_{r-}\int_E\ell_r(e)\bar K^l_r(e)\lambda(\de)>0\}}
\ell_r(e)\tilmu(\dr,\de)
\biggr)_s,\quad s\in[t,T],\\
\Theta_t=1,\quad \dP^\rho=\Theta_T\,\dP.
\end{gathered}
\]
By \Cref{lem:bsde-existence-comparison}\,\ref{lem:bsde-existence-comparison-iii}, we have
    \begin{equation*}
    \begin{aligned}
        &\quad \E\int_t^T{|\Delta_l(r)|^2\dr} \\
        &\le C\E\int_t^T{|Y_r^{l;t,x;v}|^2\Biggl(1+|X_{r}^{t,x;v}|^2 + |Y^{t,x;v}_r|^{2p_0} + |Y^{l;t,x;v}_r|^{2p_0} + |Z^{t,x;v}_r|^2 + \biggl|\int_E{\ell_r(e)K_{r}^{t,x;v}(e)}\lambda(\de)\biggr|^2\Biggr)\dr} \\
        &\le C\bigl(1+\lVert X^{t,x;v}_s \rVert_{\mathcal{S}^4}^{4} + \lVert Y^{t,x;v}_s \rVert_{\mathcal{S}^{4p_0}}^{4p_0} + \lVert Y^{l;t,x;v}_s \rVert_{\mathcal{S}^{4p_0}}^{4p_0} + \lVert Z^{t,x;v}_s \rVert_{\mathcal{L}^4}^{4} + \lVert K^{t,x;v}_s \rVert_{\mathcal{L}^4_\lambda}^{4}\bigr) \\
        &\le C(p_0)(1+|x|^{4p_0}).
    \end{aligned}
    \end{equation*}
    This estimate holds uniformly with respect to \(l\).

    Define \(A_l := \bigl\{ (s,\omega) \in [t,T]\times\Omega: |Y_s^{l;t,x;v}(\omega)| > l \bigr\}\) and \(B_l := \bigl\{ \omega \in \Omega: \text{there exists }s\in[t,T]\text{ such that }|Y_s^{l;t,x;v}(\omega)| > l \bigr\}\). Obviously, we have \(A_l \subseteq [t,T]\times B_l\). Therefore
    \begin{align*}
        \E\int_t^T{\1_{A_l}(r,\omega)\dr} &\le T\Pp(B_l) \le T\Pp\Bigl( \sup_{t\le s\le T}|Y_s^{l;t,x;v}| > l \Bigr) \\
        &\le \frac{T}{l^2}\E\Bigl[ \sup_{t\le s\le T}|Y_s^{l;t,x;v}|^2 \Bigr] \le \frac{C(1+|x|^2)}{l^2}\rightarrow 0,\quad\text{as }l \rightarrow \infty.
    \end{align*}
Due to the fact that \(\Delta_l(\omega,r)=0\) if \((r,\omega) \notin A_l\),
    we can derive from \eqref{eq:truncated-bsde-convergence-1} that
    \begin{equation}\label{eq:truncated-bsde-convergence-2}
    \begin{aligned}
        \E|\bar{Y}^l_t|^2 &\le C\varpi^2\Bigl(\frac{2A}{\rho}\Bigr) + C\E\biggl(\frac{\Theta_T}{\Theta_t}\int_t^T{|\Delta_l(r)|\1_{A_l}(r,\omega)\dr}\biggr)\\
        &\le C\varpi^2\Bigl(\frac{2A}{\rho}\Bigr)+C \E\biggl( \frac{\Theta_T^4}{\Theta_t^4} \biggr) ^{1/4} \biggl( \E\int_t^T{|\Delta_l(r)|^2\dr} \biggr) ^{1/2}\biggl( \E\int_t^T{\1_{A_l}(r,\omega)\dr} \biggr) ^{1/4} \\
        &\le C\varpi^2\Bigl(\frac{2A}{\rho}\Bigr)+C(\rho,K)l^{-1/2},
    \end{aligned}
    \end{equation}
in which, by a similar argument as in \Cref{lem:bsde-existence-comparison}\,\ref{lem:bsde-existence-comparison-i}, we can show that \(\Theta_T/\Theta_t\) has a uniform \(L^4\) bound independent of \(t\) and \(v\).
    Fix \(\rho\) first, and let \(l \rightarrow \infty\). Then we have \(C(\rho)l^{-1/2} \rightarrow 0\).
    Finally, taking \(\rho \rightarrow \infty\), we have \(\varpi^2(2A/\rho) \rightarrow 0\). Hence \(\E|\bar{Y}^l_t|^2 \rightarrow 0\) as \(l \rightarrow \infty\).

    From the proof, we see that the final constant \(C\) in \eqref{eq:truncated-bsde-convergence-2} is independent of the choices of both \(t\) and \(v\), which means that the inequality holds uniformly for any \(t\in [0,T]\) and \(v\in \Ucal_{t,T}\).
    Therefore, for any \(\varepsilon >0\) and fixed \(x_i\), \(i\in\{1,2,\cdots,N\}\), there exists \(l(i)\) such that when \(l\ge l(i)\), we have
    \begin{align*}
        \sup_{0\le t\le T}\esssup_{v\in \Ucal_{t,T}}\E| Y_{t}^{l;t,x_i;v}-Y_{t}^{t,x_i;v} |
        <\varepsilon /N.
    \end{align*}
    Then set \(l^* = \max\limits_{i\in\{1,2,\cdots,N\}}\ l(i)\). When \(l\ge l^*\), we have \(II < \varepsilon /N\).
    Thus, for any \(\varepsilon >0\), choose \(\delta\) and \(l\), such that when \(l\ge l^*\), we have
    \begin{align*}
        \sup_{(t,x)\in K}\esssup_{v\in \Ucal_{t,T}}\E| Y_{t}^{l;t,x;v}-Y_{t}^{t,x;v} |
        < I + N\cdot II < 2\varepsilon.
    \end{align*}
The proof for \Cref{lem:truncated-bsde-uniform-convergence} is complete.
\end{proof}

\begin{theorem}[Uniform convergence of truncated value functions]\label{thm:truncated-value-uniform-convergence}
    Suppose that \Cref{ass:standing} holds. If \(K \subset [0,T] \times \R^{d'}\) is a compact set, then
 for any given \(l\), \(u^{m,l}\rightarrow u^l\) uniformly in \(K\) as \(m\rightarrow \infty\).
    Moreover, \(u^l \rightarrow u\) uniformly in \(K\) as \(l \rightarrow \infty\).
\end{theorem}
\begin{proof}
    For any given \(l>0\), since we already know that \(Y^{m,l;t,x;v}_t \uparrow Y^{l;t,x;v}_t\) pathwise for any \((t,x)\in [0,T]\times \R^{d'}\) and \(v\in \Ucal_{t,T}\),
    we have \(u^{m,l}(t,x) \uparrow u^l(t,x)\) as \(m\rightarrow \infty\).
Moreover, since both \(u^{m,l}\) and \(u^l\) are continuous with respect to \((t,x)\) in the compact set \(K\), Dini's theorem yields that \(u^{m,l}\rightarrow u^l\) uniformly in \(K\), as \(m\rightarrow \infty\).

    To prove the second claim, for any given \(\varepsilon >0\), from \Cref{thm:value-random-initial-data}, we know that for any \((t,x)\in [0,T]\times \R^{d'}\), there exist \(v^1_{t,x},v^2_{t,x}\in \Ucal_{t,T}\) such that
    \begin{align*}
        Y^{l;t,x;v^1_{t,x}}_t \le u^l(t,x) \le Y^{l;t,x;v^2_{t,x}}_t + \varepsilon\ \ \ \text{and}\ \ \
        Y^{t,x;v^2_{t,x}}_t \le u(t,x) \le Y^{t,x;v^1_{t,x}}_t + \varepsilon.
    \end{align*}
    Since both \(u(t,x)\) and \(u^l(t,x)\) are deterministic, it follows that
    \begin{equation*}
    \begin{aligned}
        \sup_{(t,x)\in K}|u^l(t,x) - u(t,x)| &\le \sup_{(t,x)\in K}\bigl( \E| Y^{l;t,x;v^1_{t,x}}_t - Y^{t,x;v^1_{t,x}}_t | + \E| Y^{l;t,x;v^2_{t,x}}_t - Y^{t,x;v^2_{t,x}}_t | + 2\varepsilon \bigr) \\
        &\le 2\sup_{(t,x)\in K}\esssup_{v\in \Ucal_{t,T}}\E| Y^{l;t,x;v}_t - Y^{t,x;v}_t | + 2\varepsilon\le 6\varepsilon,
    \end{aligned}
    \end{equation*}
    where the last inequality comes from \Cref{lem:truncated-bsde-uniform-convergence}.
    Hence we have the uniform convergence of \(u^l\) to \(u\) on \(K\).
\end{proof}

\begin{theorem}[Stability of viscosity solutions]\label{thm:viscosity-stability}
Suppose that \Cref{ass:standing} holds. Let
\[
f_n:
[0,T]\times\R^{d'}\times\R\times\R^d\times\R\times U
\longrightarrow\R,
\qquad n\in\mathbb N,
\]
be a sequence of generators satisfying \Cref{ass:A1,ass:A5,ass:A6,ass:A7}, where the corresponding constants and moduli of
continuity may depend on \(n\). Suppose further that

    \begin{enumerate}[label=\textup{(\roman*)},ref=\textup{(\roman*)},wide=0pt]
    \item \(u_n\) converges to \(u\) uniformly in every compact subset of \([0,T]\times \R^{d'}\) as \(n\to\infty\),

    \item \(f_n\) converges to \(f\) uniformly in every compact subset of \([0,T]\times \R^{d'} \times \R \times \R^d \times \R \times U\) as \(n\to\infty\),

    \item \(u_n\) is a viscosity solution of the HJB equation \eqref{eq:hjb} associated with \(f_n\) for any \(n\in\mathbb{N}\).
    \end{enumerate}
    Then \(u\) is a viscosity solution of the HJB equation \eqref{eq:hjb} associated with \(f\).
\end{theorem}
\begin{proof}
 We only verify the viscosity subsolution property using \Cref{def:viscosity-local-i}. The supersolution property can be verified similarly.

    We choose \(\delta=1/(2L_2+2)\) and \(\Lambda=L_2+1\).
    Let \(\varphi \in C^{1,2}([0,T]\times \R^{d'})\) such that \(u-\varphi\) attains its strict local maximum at \((t,x)\) in \(B((t,x),R_{\delta,x})\), and let \(\tilde{\varphi}\) be the continuation of \(\varphi\) on \(B((t,x),3R_{\delta,x})\) such that \(\tilde{\varphi} = \varphi\) on \(B((t,x),L_2(1+|x|)\delta)\),
    \(\tilde{\varphi} \in C^{1,2}([0,T]\times \R^{d'})\) and \(u-\tilde{\varphi}\) attains its strict local maximum at \((t,x)\) in \(B((t,x),3R_{\delta,x})\). Since as \(n\to\infty\), \(u_n\rightarrow u\) uniformly in any compact subset of \([0,T]\times \R^{d'}\), there exists a sequence \((t_n,x_n)\) such that \((t_n,x_n)\rightarrow (t,x)\) and \(u_n-\tilde{\varphi}\) attains its local maximum at \((t_n,x_n)\) in \(B((t,x),3R_{\delta,x})\).
Without loss of generality, choose a subsequence such that \((t_n,x_n) \in B((t,x),R_{\delta,x})\) for any \(n\in\mathbb{N}\).
 Obviously, \(B((t_n,x_n),R_{\delta,x_n}) \subset B((t,x),3R_{\delta,x})\).
    So from the viscosity subsolution property of \(u_n\), we know
    \begin{equation*}
    \begin{aligned}
        \frac{\partial \tilde{\varphi}}{\partial t}(t_n,x_n) + \underset{v\in U}{\sup}\Bigl\{&\Aop\tilde{\varphi}(t_n,x_n,v) + \Bop^\delta(\tilde{\varphi},u_n)(t_n,x_n,v)\\
        &+f_n( t_n,x_n,u_n(t_n,x_n),\sigma^\top(t_n,x_n,v)D\tilde{\varphi}(t_n,x_n),\Cop^\delta(\tilde{\varphi},u_n)(t_n,x_n,v),v ) \Bigr\} \ge 0.
    \end{aligned}
    \end{equation*}
    For \(v\in U\), define
    \begin{align*}
        H_n(v):=&\Aop\tilde{\varphi}(t_n,x_n,v)+\Bop^\delta(\tilde{\varphi},u_n)(t_n,x_n,v) \\
        &+f_n( t_n,x_n,u_n(t_n,x_n),\sigma^\top(t_n,x_n,v)D\tilde{\varphi}(t_n,x_n),\Cop^\delta(\tilde{\varphi},u_n)(t_n,x_n,v),v ),
    \end{align*}
    and
    \begin{align*}
        H(v):=&\Aop\tilde{\varphi}(t,x,v)+\Bop^\delta(\tilde{\varphi},u)(t,x,v) 
        +f( t,x,u(t,x),\sigma^\top(t,x,v)D\tilde{\varphi}(t,x),\Cop^\delta(\tilde{\varphi},u)(t,x,v),v ).
    \end{align*}
For any \(n\in\mathbb{N}\), choose \(v_n\in U\) such that \(H_n(v_n)\ge \sup_{v\in U}H_n(v)-\frac1n\).
Then
    \begin{align}\label{eq:stability-viscosity-inequality-n}
        \frac{\partial \tilde{\varphi}}{\partial t}(t_n,x_n)+H_n(v_n)\ge -\frac1n.
    \end{align}
    Since \(U\) is compact, taking a subsequence if necessary, we have
    \(v_n\to \bar v\in U\).

Then we show that as \(n\to\infty\),
    \begin{align}\label{eq:stability-hamiltonian-convergence}
        H_n(v_n)\longrightarrow H(\bar v).
    \end{align}
First note that
    \begin{align*}
        &\Aop\tilde{\varphi}(t_n,x_n,v_n)\longrightarrow \Aop\tilde{\varphi}(t,x,\bar v),\\
        &\sigma^\top(t_n,x_n,v_n)D\tilde{\varphi}(t_n,x_n)\longrightarrow \sigma^\top(t,x,\bar v)D\tilde{\varphi}(t,x),\\
        &u_n(t_n,x_n)\longrightarrow u(t,x).
    \end{align*}
On the other hand, by the dominated convergence theorem we obtain
    \begin{align*}
        \Bop^\delta(\tilde{\varphi},u_n)(t_n,x_n,v_n)\longrightarrow \Bop^\delta(\tilde{\varphi},u)(t,x,\bar v)\ \ \ \text{and}\ \ \
        \Cop^\delta(\tilde{\varphi},u_n)(t_n,x_n,v_n)\longrightarrow \Cop^\delta(\tilde{\varphi},u)(t,x,\bar v).
    \end{align*}
Since \(f_n\to f\) uniformly in any compact subset as \(n\to\infty\), we conclude that
    \begin{align*}
        &f_n( t_n,x_n,u_n(t_n,x_n),\sigma^\top(t_n,x_n,v_n)D\tilde{\varphi}(t_n,x_n),\Cop^\delta(\tilde{\varphi},u_n)(t_n,x_n,v_n),v_n ) \notag\\
        &\qquad\longrightarrow
        f( t,x,u(t,x),\sigma^\top(t,x,\bar v)D\tilde{\varphi}(t,x),\Cop^\delta(\tilde{\varphi},u)(t,x,\bar v),\bar v ),
    \end{align*}
which implies \eqref{eq:stability-hamiltonian-convergence}.

Taking \(n\to\infty\) in \eqref{eq:stability-viscosity-inequality-n}, we get
    \begin{align*}
        \frac{\partial \tilde{\varphi}}{\partial t}(t,x)+H(\bar v)\ge 0.
    \end{align*}
    Since \(H(\bar v)\le \sup_{v\in U}H(v)\), it follows that
    \begin{align*}
        \frac{\partial \tilde{\varphi}}{\partial t}(t,x) + \underset{v\in U}{\sup}\Bigl\{&\Aop\tilde{\varphi}(t,x,v)
        + \Bop^\delta(\tilde{\varphi},u)(t,x,v)
        +f( t,x,u(t,x),\sigma^\top(t,x,v)D\tilde{\varphi}(t,x),\Cop^\delta(\tilde{\varphi},u)(t,x,v),v ) \Bigr\} \ge 0.
    \end{align*}
    It is vital to see that \(\tilde{\varphi}(t,x+\gamma(t,x,v,e)) = \varphi(t,x+\gamma(t,x,v,e))\) for any \(e\in E_\delta\) and \(v\in U\).
    Hence we have
    \begin{align*}
        \frac{\partial \varphi}{\partial t}(t,x) + \underset{v\in U}{\sup}\Bigl\{&\Aop\varphi(t,x,v)
        + \Bop^\delta(\varphi,u)(t,x,v)
        +f( t,x,u(t,x),\sigma^\top(t,x,v)D\varphi(t,x),\Cop^\delta(\varphi,u)(t,x,v),v ) \Bigr\} \ge 0.
    \end{align*}

Thus \(u\) satisfies the subsolution inequality in
\Cref{def:viscosity-local-i}. By
\Cref{lem:viscosity-definitions-equivalent}, \(u\) is a viscosity
subsolution of \eqref{eq:hjb}.
The supersolution case is analogous,
except that one fixes an arbitrary \(v\in U\), passes to the limit in
the corresponding inequality, and then takes the supremum over \(v\in U\).
\end{proof}

\begin{theorem}[Viscosity solution property]\label{thm:value-viscosity-solution}
Suppose that \Cref{ass:standing} holds. Then the value function
\(u(t,x)\) is a viscosity solution of the HJB equation~\eqref{eq:hjb}.
\end{theorem}
\begin{proof}
For each \(l>0\) and \(m>A\), \(f^{m,l}\) is Lipschitz continuous with respect to \((x,y,z,k,v)\).
By \cite[Theorem~4.1]{LiPeng2009}
(the proof is very similar when the kernel \(\ell(x,e)\) is replaced by the kernel \(\ell(t,e)\) used in our situation),
it follows immediately that \(u^{m,l}\) is a viscosity solution
of the HJB equation associated with \(f^{m,l}\).

Fix \(l>0\). By \Cref{thm:truncated-value-uniform-convergence}, \(u^{m,l}\to u^{l}\) uniformly on every compact subset of
\([0,T]\times \R^{d'}\) as \(m\to\infty\). Moreover,
\(f^{m,l}\to f^{l}\) uniformly in any compact subset of
\([0,T]\times \R^{d'}\times \R\times \R^{d}\times \R\times U\).
Hence, by \Cref{thm:viscosity-stability}, \(u^{l}\) is a viscosity solution of  the HJB
equation associated with \(f^{l}\).
Furthermore, by \Cref{thm:truncated-value-uniform-convergence} again, \(u^{l}\to u\) uniformly in every compact subset of
\([0,T]\times \R^{d'}\) as \(l\to\infty\). On the other hand,
we know \(f^{l}\to f\) uniformly in every compact
subset of \([0,T]\times \R^{d'}\times \R\times \R^{d}\times \R\times U\).
Therefore, by \Cref{thm:viscosity-stability}, we know that \(u\) is a
viscosity solution of the HJB equation \eqref{eq:hjb}.
\end{proof}

\section{Uniqueness of viscosity solutions}\label{sec:uniqueness}

In this section, we may assume without loss of generality that the
monotonicity constant \(\kappa\) in \assitemref{ass:A5} is negative, by an
exponential change of variables if necessary.

We denote by \(\Pcal^{2,+}u(t,x)\) and \(\Pcal^{2,-}u(t,x)\) the parabolic superjet and subjet of \(u\) at \((t,x)\), respectively, defined as follows:
\begin{align*}
    \Pcal^{2,+}u(t,x) &= \Bigl\{ (a,p,X)\in \R\times \R^{d'} \times S^{d'}: u(s,y) \le u(t,x) + a(s-t) + \langle p,y-x \rangle\\
    & \qquad +\frac{1}{2}\langle X(y-x),y-x \rangle + o(|s-t|+|y-x|^2),\ \text{as }(s,y) \to (t,x) \Bigr\}; \\
    \Pcal^{2,-}u(t,x) &= \Bigl\{ (a,p,X)\in \R\times \R^{d'} \times S^{d'}: u(s,y) \ge u(t,x) + a(s-t) + \langle p,y-x \rangle \\
    & \qquad +\frac{1}{2}\langle X(y-x),y-x \rangle + o(|s-t|+|y-x|^2),\ \text{as }(s,y) \to (t,x) \Bigr\}.
\end{align*}

The limiting parabolic semijets \(\overline{\Pcal}^{2,\pm}u(t,x)\) consist of all triples \((a,p,X)\) for which there are points \((t_j,x_j)\) in the domain of \(u\) and triples
\(
(a_j,p_j,X_j)\in\Pcal^{2,\pm}u(t_j,x_j)
\)
such that
\[
(t_j,x_j,u(t_j,x_j),a_j,p_j,X_j)
\longrightarrow (t,x,u(t,x),a,p,X).
\]

The next lemma follows from the nonlocal Jensen--Ishii lemma of Barles and Imbert \cite[Lemma~1 and Corollary~2]{BarlesImbert2008}; see also \cite[Remarks~1 and 2 in Section~1]{BarlesImbert2008}.
Related parabolic formulations for fully nonlinear PIDEs can also be found in 
Jakobsen and Karlsen~\cite[Theorem~2.2]{JakobsenKarlsen2005} and Seydel~\cite[Corollary~5.13]{Seydel2009}.

\begin{lemma}[Nonlocal parabolic Jensen--Ishii lemma]\label{lem:jensen-ishii}
Suppose that \Cref{ass:standing} holds. Let
\(u\in C((0,T)\times\mathbb R^{d'})\) be a viscosity subsolution of
\eqref{eq:hjb}, bounded from above, and let
\(w\in C((0,T)\times\mathbb R^{d'})\) be a viscosity supersolution of
\eqref{eq:hjb}, bounded from below. Let
\[
 \Phi\in C^2((0,T)\times\mathbb R^{d'}\times\mathbb R^{d'}),
\]
and suppose that \((\bar t,\bar x,\bar y)\) is a global maximum point of
\[
 (t,x,y)\longmapsto u(t,x)-w(t,y)-\Phi(t,x,y).
\]
Set
\[
 p_1:=D_x\Phi(\bar t,\bar x,\bar y),
 \qquad
 p_2:=-D_y\Phi(\bar t,\bar x,\bar y).
\]
For \(\alpha>0\), define the sup-convolution
\begin{equation*}
\begin{aligned}
 \Phi_\alpha(x,y)
 :=\sup_{|(\xi,\eta)-(x,y)|\le 1}\biggl\{
 \Phi(\bar t,\xi,\eta)
 -p_1\cdot(\xi-x)+p_2\cdot(\eta-y)
 -\frac{|\xi-x|^2+|\eta-y|^2}{2\alpha}
 \biggr\}.
\end{aligned}
\end{equation*}
For any \(\delta>0\), there exists \(\bar\alpha=\bar\alpha(\delta)>0\)
such that, for \(0<\alpha<\bar\alpha\), there exist
\[
 (a_1,p_1,X_1)\in\overline{\mathcal P}^{2,+}u(\bar t,\bar x),
 \qquad
 (a_2,p_2,X_2)\in\overline{\mathcal P}^{2,-}w(\bar t,\bar y),
\]
satisfying
\begin{equation*}
 a_1-a_2=\partial_t\Phi(\bar t,\bar x,\bar y),
\end{equation*}
and
\begin{equation*}
 -\frac1\alpha I_{2d'}
 \le
 \begin{pmatrix}X_1&0\\0&-X_2\end{pmatrix}
 \le D^2_{x,y}\Phi_\alpha(\bar x,\bar y)
 =D^2_{x,y}\Phi(\bar t,\bar x,\bar y)+o_\alpha(1).
\end{equation*}
Let
\[
 \Phi_\alpha^u(t,x):=\Phi_\alpha(x,\bar y),
 \qquad
 \Phi_\alpha^w(t,y):=-\Phi_\alpha(\bar x,y).
\]
Then
\begin{equation*}
\begin{aligned}
0\le\;&a_1+\sup_{\nu\in U}\biggl\{
 \frac12\operatorname{tr}\bigl[
   \sigma\sigma^\top(\bar t,\bar x,\nu)X_1
 \bigr]
 +\langle b(\bar t,\bar x,\nu),p_1\rangle
 +\Bop^\delta\bigl(\Phi_\alpha^u,u\bigr)(\bar t,\bar x,\nu)\\
&\qquad
 +f\bigl(\bar t,\bar x,u(\bar t,\bar x),
          \sigma^\top(\bar t,\bar x,\nu)p_1,
          \Cop^\delta\bigl(\Phi_\alpha^u,u\bigr)(\bar t,\bar x,\nu),\nu\bigr)
\biggr\},
\end{aligned}
\end{equation*}
and
\begin{equation*}
\begin{aligned}
0\ge\;&a_2+\sup_{\nu\in U}\biggl\{
 \frac12\operatorname{tr}\bigl[
   \sigma\sigma^\top(\bar t,\bar y,\nu)X_2
 \bigr]
 +\langle b(\bar t,\bar y,\nu),p_2\rangle
 +\Bop^\delta\bigl(\Phi_\alpha^w,w\bigr)(\bar t,\bar y,\nu)\\
&\qquad
 +f\bigl(\bar t,\bar y,w(\bar t,\bar y),
          \sigma^\top(\bar t,\bar y,\nu)p_2,
          \Cop^\delta\bigl(\Phi_\alpha^w,w\bigr)(\bar t,\bar y,\nu),\nu\bigr)
\biggr\}.
\end{aligned}
\end{equation*}
\end{lemma}

\begin{theorem}[Comparison theorem]\label{thm:comparison}
Suppose that \Cref{ass:standing} holds.
Let \(0\le a<b\le T\), and let
\(u_1,u_2\in C([a,b]\times\R^{d'})\)
be a viscosity subsolution and a viscosity
supersolution of the HJB equation \eqref{eq:hjb}
on \((a,b)\times\R^{d'}\), respectively.
Suppose that \(u_1\) is bounded from above and \(u_2\)
is bounded from below. If
\[
    u_1(b,x)\le u_2(b,x), \qquad x\in\R^{d'},
\]
then
\[
    u_1(t,x)\le u_2(t,x), \qquad (t,x)\in[a,b]\times\R^{d'}.
\]
\end{theorem}
\begin{proof}
For simplicity, we only prove the case \(a=0\) and \(b=T\).
Argue by contradiction. Suppose that
\(\sup_{(0,T]\times\R^{d'}}(u_1-u_2)>0\).
Choose \(\beta>0\) small enough such that
\(N:=\sup_{(0,T]\times\mathbb R^{d'}}\{u_1(t,x)-u_2(t,x)-\beta/t\}>0\).
For \(\eta>0\), define
\[
N^\eta:=\sup_{(0,T]\times\mathbb R^{d'}}\Bigl\{ u_1(t,x)-u_2(t,x)-\frac{\beta}{t}-\eta|x|^2 \Bigr\}.
\]

Since \(u_1\) is bounded from above and \(u_2\) is bounded from
below, the supremum \(N^\eta\) is attained. Moreover,
\(
    N^\eta\longrightarrow N
    \text{ as }
    \eta\downarrow0.
\)
For a given \(\eta>0\) and any \(\varepsilon>0\), define
\[
N^{\eta,\varepsilon}:=\sup_{(0,T]\times\mathbb R^{d'}\times\mathbb R^{d'}}\biggl\{u_1(t,x)-u_2(t,y)-\frac{\beta}{t}-\eta|x|^2-\frac{|x-y|^2}{\varepsilon^2}\biggr\}.
\]
Let \((t_\varepsilon,x_\varepsilon,y_\varepsilon)\) be a maximum point
of \(N^{\eta,\varepsilon}\). By a standard argument, 
after passing to a subsequence if necessary, there exists
a maximum point \((t_\eta,x_\eta)\) of \(N^\eta\) such that
\[
(t_\varepsilon,x_\varepsilon,y_\varepsilon)
\to
(t_\eta,x_\eta,x_\eta),
\ \ \ \text{and}\ \ \
\frac{|x_\varepsilon-y_\varepsilon|^2}{\varepsilon^2}\to0,
 \ \ \
\text{as }\varepsilon\downarrow0.
\]
In addition,
\(
    \eta|x_\eta|^2\longrightarrow0
    \text{ as }
    \eta\downarrow0.
\)

Fix \(\eta>0\) sufficiently small so that \(N^\eta>0\).
We claim that \(t_\varepsilon<T\) for all sufficiently small
\(\varepsilon\).
Indeed, for fixed \(\eta\), the maximizing points
\(x_\varepsilon,y_\varepsilon\) remain in a compact subset
\(K_\eta\) of \(\R^{d'}\).
Thus, there exists a continuity modulus \(\omega_\eta\) of 
\(u_2(T,\cdot)\) on \(K_\eta\), such that
\(
    \omega_\eta(r)\rightarrow0 \text{ as }r\downarrow0.
\)
If \(t_\varepsilon=T\) along a subsequence, then we have
\[
    0<N^\eta
    \le N^{\eta,\varepsilon}
    \le
    u_2(T,x_\varepsilon)-u_2(T,y_\varepsilon)-\frac{\beta}{T}
    \le
    \omega_\eta(|x_\varepsilon-y_\varepsilon|)-\frac{\beta}{T}.
\]
Since \(|x_\varepsilon-y_\varepsilon|\to0\) as
\(\varepsilon\downarrow0\), the right-hand side converges
to \(-\beta/T<0\), which is a contradiction.
Hence \(t_\varepsilon<T\) for all sufficiently small
\(\varepsilon\).

Set \(\Phi(t,x,y)=\beta/t+\eta|x|^2+|x-y|^2/\varepsilon^2\).
Then we have
\[
p_1=D_x\Phi(t_\varepsilon,x_\varepsilon,y_\varepsilon)
=2\eta x_\varepsilon+2\frac{x_\varepsilon-y_\varepsilon}{\varepsilon^2},
\quad p_2=-D_y\Phi(t_\varepsilon,x_\varepsilon,y_\varepsilon)
=2\frac{x_\varepsilon-y_\varepsilon}{\varepsilon^2}.
\]
\[
\partial_t\Phi(t_\varepsilon,x_\varepsilon,y_\varepsilon)
=-\frac{\beta}{t_\varepsilon^2}, \quad
D^2_{x,y}\Phi(t_\varepsilon,x_\varepsilon,y_\varepsilon)=\frac{2}{\varepsilon^2}
\begin{pmatrix}
(\eta \varepsilon^2+1) I & -I \\
-I & I
\end{pmatrix}.
\]
Applying \Cref{lem:jensen-ishii}, we obtain
\[
(a_1,p_1,X_1)\in \bar{\mathcal P}^{2,+}u_1(t_\varepsilon,x_\varepsilon)\ \ \ \text{and}\ \ \
(a_2,p_2,X_2)\in \bar{\mathcal P}^{2,-}u_2(t_\varepsilon,y_\varepsilon)
\]
with
\[
a_1-a_2=-\frac{\beta}{t_\varepsilon^2}\ \ \ \text{and}\ \ \
-\frac1\alpha I_{2d'}
\le
\begin{pmatrix}
X_1&0\\
0&-X_2
\end{pmatrix}
\le
D^2_{x,y}\Phi(t_\varepsilon,x_\varepsilon,y_\varepsilon)+o_\alpha(1)I_{2d'}.
\]
Hence we have the viscosity inequalities for \(u_1\) and \(u_2\) at \((t_\varepsilon,x_\varepsilon)\) and \((t_\varepsilon,y_\varepsilon)\), respectively,
\begin{align*}
    a_1+\sup_{v\in U}&\biggl\{ \frac{1}{2}\tr\bigl( \sigma \sigma^\top( t_{\varepsilon},x_{\varepsilon},v ) X_1 \bigr) +b^\top( t_{\varepsilon},x_{\varepsilon},v ) p_1 + \Bop^\delta(\Phi_\alpha(\cdot,y_\varepsilon),u_1)(t_{\varepsilon},x_{\varepsilon},v) \\
        &\quad +  f(t_{\varepsilon},x_{\varepsilon},u_1(t_{\varepsilon},x_{\varepsilon}),\sigma^\top(t_{\varepsilon},x_{\varepsilon},v)p_1,\Cop^\delta(\Phi_\alpha(\cdot,y_\varepsilon),u_1)(t_{\varepsilon},x_{\varepsilon},v),v) \biggr\} \ge 0,\\
    a_2+\sup_{v\in U}&\biggl\{ \frac{1}{2}\tr\bigl( \sigma \sigma^\top( t_{\varepsilon},y_{\varepsilon},v ) X_2 \bigr) +b^\top( t_{\varepsilon},y_{\varepsilon},v ) p_2 + \Bop^\delta(-\Phi_\alpha(x_\varepsilon,\cdot),u_2)(t_{\varepsilon},y_{\varepsilon},v) \\
        &\quad +  f(t_{\varepsilon},y_{\varepsilon},u_2(t_{\varepsilon},y_{\varepsilon}),\sigma^\top(t_{\varepsilon},y_{\varepsilon},v)p_2,\Cop^\delta(-\Phi_\alpha(x_\varepsilon,\cdot),u_2)(t_{\varepsilon},y_{\varepsilon},v),v) \biggr\} \le 0,
\end{align*}
which gives
\begin{equation}\label{eq:viscosity-inequality}
\begin{aligned}
-\frac{\beta}{t_\varepsilon^2}+\sup_{v\in U}\biggl\{ &\frac{1}{2}\tr\bigl( \sigma \sigma^\top( t_{\varepsilon},x_{\varepsilon},v ) X_1 -\sigma \sigma^\top( t_{\varepsilon},y_{\varepsilon},v ) X_2 \bigr) +b^\top( t_{\varepsilon},x_{\varepsilon},v ) p_1 - b^\top( t_{\varepsilon},y_{\varepsilon},v ) p_2\\
&+ \Bop^\delta(\Phi_\alpha(\cdot,y_\varepsilon),u_1)(t_{\varepsilon},x_{\varepsilon},v) - \Bop^\delta(-\Phi_\alpha(x_\varepsilon,\cdot),u_2)(t_{\varepsilon},y_{\varepsilon},v) \\
&+ f(t_{\varepsilon},x_{\varepsilon},u_1(t_{\varepsilon},x_{\varepsilon}),\sigma^\top(t_{\varepsilon},x_{\varepsilon},v)p_1,\Cop^\delta(\Phi_\alpha(\cdot,y_\varepsilon),u_1)(t_{\varepsilon},x_{\varepsilon},v),v) \\
&- f(t_{\varepsilon},y_{\varepsilon},u_2(t_{\varepsilon},y_{\varepsilon}),\sigma^\top(t_{\varepsilon},y_{\varepsilon},v)p_2,\Cop^\delta(-\Phi_\alpha(x_\varepsilon,\cdot),u_2)(t_{\varepsilon},y_{\varepsilon},v),v) \biggr\} \ge 0.
\end{aligned}
\end{equation}

Then we estimate the terms in \eqref{eq:viscosity-inequality}. First, by the Lipschitz continuity of \(\sigma\) and \(b\), we have
\begin{equation*}
\begin{aligned}
\frac12\operatorname{tr}\bigl[
\sigma\sigma^\top(t_\varepsilon,x_\varepsilon,v)X_1
-\sigma\sigma^\top(t_\varepsilon,y_\varepsilon,v)X_2\bigr]
\le C\frac{|x_\varepsilon-y_\varepsilon|^2}{\varepsilon^2}
+C\eta(1+|x_\varepsilon|^2)
+o_\alpha(1)(1+|x_\varepsilon|^2+|y_\varepsilon|^2),
\end{aligned}
\end{equation*}
\begin{equation*}
\begin{aligned}
|\sigma^\top(t_\varepsilon,x_\varepsilon,v)p_1-\sigma^\top(t_\varepsilon,y_\varepsilon,v)p_2|
\le C\frac{|x_\varepsilon-y_\varepsilon|^2}{\varepsilon^2}
+C\eta(1+|x_\varepsilon|^2),
\end{aligned}
\end{equation*}
and
\begin{align*}
b(t_\varepsilon,x_\varepsilon,v)\cdot p_1
-b(t_\varepsilon,y_\varepsilon,v)\cdot p_2
\le C\frac{|x_\varepsilon-y_\varepsilon|^2}{\varepsilon^2}
+C\eta(1+|x_\varepsilon|^2).
\end{align*}
For the small-jump part in the nonlocal terms, we have
\begin{equation}\label{eq:small-jump-B}
\begin{aligned}
&\int_{E_\delta}
\Bigl\{
\Phi(t_\varepsilon,x_\varepsilon+\gamma(t_\varepsilon,x_\varepsilon,v,e),y_\varepsilon)
-\Phi(t_\varepsilon,x_\varepsilon,y_\varepsilon)
-p_1\cdot\gamma(t_\varepsilon,x_\varepsilon,v,e)
\\
&\qquad \ -
\bigl[
-\Phi(t_\varepsilon,x_\varepsilon,y_\varepsilon+\gamma(t_\varepsilon,y_\varepsilon,v,e))
+\Phi(t_\varepsilon,x_\varepsilon,y_\varepsilon)
-p_2\cdot\gamma(t_\varepsilon,y_\varepsilon,v,e)
\bigr]
\Bigr\}\lambda(\de)\\
\le &C\biggl(\eta+\frac{1}{\varepsilon^2}\biggr)(1+|x_\varepsilon|^2+|y_\varepsilon|^2)\int_{E_\delta}{(1\wedge |e|^2)\lambda(\de)},
\end{aligned}
\end{equation}
and
\begin{equation*}
\begin{aligned}
&\int_{E_\delta}
\Bigl\{\Phi(t_\varepsilon,x_\varepsilon+\gamma(t_\varepsilon,x_\varepsilon,v,e),y_\varepsilon)
-\Phi(t_\varepsilon,x_\varepsilon,y_\varepsilon)
\\&\qquad \ -[
-\Phi(t_\varepsilon,x_\varepsilon,y_\varepsilon+\gamma(t_\varepsilon,y_\varepsilon,v,e))
+\Phi(t_\varepsilon,x_\varepsilon,y_\varepsilon)
]
\Bigr\}
\ell(t_\varepsilon,e)\,\lambda(\de)\\
\le&C\biggl\{\biggl(\eta+\frac{1}{\varepsilon^2}\biggr)(1+|x_\varepsilon|^2+|y_\varepsilon|^2)
+\frac{|x_\varepsilon-y_\varepsilon|^2}{\varepsilon^2}\biggr\}\int_{E_\delta}(1\wedge |e|^2)\lambda(\de).
\end{aligned}
\end{equation*}
For the large-jump part in the nonlocal terms, the maximality of
\(u_1(t,x)-u_2(t,y)-\Phi(t,x,y)\)
at \((t_\varepsilon,x_\varepsilon,y_\varepsilon)\) implies
\[
\begin{aligned}
&u_1(t_\varepsilon,x_\varepsilon+\gamma(t_\varepsilon,x_\varepsilon,v,e))
-u_1(t_\varepsilon,x_\varepsilon)
-u_2(t_\varepsilon,y_\varepsilon+\gamma(t_\varepsilon,y_\varepsilon,v,e))
+u_2(t_\varepsilon,y_\varepsilon)
\\\le&
\Phi(t_\varepsilon,x_\varepsilon+\gamma(t_\varepsilon,x_\varepsilon,v,e),y_\varepsilon+\gamma(t_\varepsilon,y_\varepsilon,v,e))
-\Phi(t_\varepsilon,x_\varepsilon,y_\varepsilon)\\
=&\frac{2}{\varepsilon^2}\ip{x_\varepsilon-y_\varepsilon}{\gamma(t_\varepsilon,x_\varepsilon,v,e)-\gamma(t_\varepsilon,y_\varepsilon,v,e)}
+\frac{1}{\varepsilon^2}|\gamma(t_\varepsilon,x_\varepsilon,v,e)-\gamma(t_\varepsilon,y_\varepsilon,v,e)|^2
\\&+2\eta \ip{x_\varepsilon}{\gamma(t_\varepsilon,x_\varepsilon,v,e)}
+\eta|\gamma(t_\varepsilon,x_\varepsilon,v,e)|^2.
\end{aligned}
\]
Hence we have
\begin{equation}
\begin{aligned}
&\int_{E_\delta^c}
\Bigl\{u_1(t_\varepsilon,x_\varepsilon+\gamma(t_\varepsilon,x_\varepsilon,v,e))
-u_1(t_\varepsilon,x_\varepsilon)-p_1\cdot\gamma(t_\varepsilon,x_\varepsilon,v,e)\\
&\qquad \  -[u_2(t_\varepsilon,y_\varepsilon+\gamma(t_\varepsilon,y_\varepsilon,v,e))
-u_2(t_\varepsilon,y_\varepsilon)-p_2\cdot\gamma(t_\varepsilon,y_\varepsilon,v,e)]\Bigr\}\lambda(\de)\\
\le&C\eta(1+|x_\varepsilon|^2)
+C\frac{|x_\varepsilon-y_\varepsilon|^2}{\varepsilon^2},
\end{aligned}
\end{equation}
and
\begin{equation}\label{eq:large-jump-C}
\begin{aligned}
&\int_{E_\delta^c}\Bigl[
u_1(t_\varepsilon,x_\varepsilon+\gamma(t_\varepsilon,x_\varepsilon,v,e))
-u_1(t_\varepsilon,x_\varepsilon)
-
u_2(t_\varepsilon,y_\varepsilon+\gamma(t_\varepsilon,y_\varepsilon,v,e))
+u_2(t_\varepsilon,y_\varepsilon)\Bigr]
\ell(t_\varepsilon,e)\,\lambda(\de)\\
\le&C\eta(1+|x_\varepsilon|^2)
+C\frac{|x_\varepsilon-y_\varepsilon|^2}{\varepsilon^2}.
\end{aligned}
\end{equation}
Finally, for fixed \(\delta,\varepsilon,\eta\), by the local \(C^2\)-
convergence \(\Phi_\alpha\to\Phi\) and \eqref{eq:small-jump-B}-\eqref{eq:large-jump-C}, we have
\begin{equation*}
\begin{aligned}
&\quad \Bop^\delta(\Phi_\alpha(\cdot,y_\varepsilon),u_1)(t_\varepsilon,x_\varepsilon,v)
-\Bop^\delta(-\Phi_\alpha(x_\varepsilon,\cdot),u_2)(t_\varepsilon,y_\varepsilon,v)
\\&\le \Bop^\delta(\Phi(\cdot,y_\varepsilon),u_1)(t_\varepsilon,x_\varepsilon,v)
-\Bop^\delta(-\Phi(x_\varepsilon,\cdot),u_2)(t_\varepsilon,y_\varepsilon,v)
+o_\alpha(1)
\\&\le C\biggl(\eta+\frac{1}{\varepsilon^2}\biggr)(1+|x_\varepsilon|^2+|y_\varepsilon|^2)o_\delta(1)
+C\eta(1+|x_\varepsilon|^2)+C\frac{|x_\varepsilon-y_\varepsilon|^2}{\varepsilon^2}
+o_\alpha(1),
\end{aligned}
\end{equation*}
and
\begin{equation*}
\begin{aligned}
&\quad \Cop^\delta(\Phi_\alpha(\cdot,y_\varepsilon),u_1)(t_\varepsilon,x_\varepsilon,v)
-\Cop^\delta(-\Phi_\alpha(x_\varepsilon,\cdot),u_2)(t_\varepsilon,y_\varepsilon,v)
\\&\le \Cop^\delta(\Phi(\cdot,y_\varepsilon),u_1)(t_\varepsilon,x_\varepsilon,v)
-\Cop^\delta(-\Phi(x_\varepsilon,\cdot),u_2)(t_\varepsilon,y_\varepsilon,v)
+o_\alpha(1)
\\&\le C\biggl(\eta+\frac{1}{\varepsilon^2}\biggr)(1+|x_\varepsilon|^2+|y_\varepsilon|^2)o_\delta(1)
+C\eta(1+|x_\varepsilon|^2)+C\frac{|x_\varepsilon-y_\varepsilon|^2}{\varepsilon^2}+o_\alpha(1).
\end{aligned}
\end{equation*}
Therefore, for any \(v\in U\), we have
\begin{equation}\label{eq:estimate-hamiltonian}
\begin{aligned}
&\frac{1}{2}\tr\bigl( \sigma \sigma^\top( t_{\varepsilon},x_{\varepsilon},v ) X_1 -\sigma \sigma^\top( t_{\varepsilon},y_{\varepsilon},v ) X_2 \bigr) +b^\top( t_{\varepsilon},x_{\varepsilon},v ) p_1 - b^\top( t_{\varepsilon},y_{\varepsilon},v ) p_2
\\&+\Bop^\delta(\Phi_\alpha(\cdot,y_\varepsilon),u_1)(t_\varepsilon,x_\varepsilon,v)
-\Bop^\delta(-\Phi_\alpha(x_\varepsilon,\cdot),u_2)(t_\varepsilon,y_\varepsilon,v)
\\
&+f(
t_\varepsilon,x_\varepsilon,u_1(t_\varepsilon,x_\varepsilon),
\sigma^\top(t_\varepsilon,x_\varepsilon,v)p_1,
\Cop^\delta(\Phi_\alpha(\cdot,y_\varepsilon),u_1)(t_\varepsilon,x_\varepsilon,v),v
)
\\
&-f(
t_\varepsilon,y_\varepsilon,u_2(t_\varepsilon,y_\varepsilon),
\sigma^\top(t_\varepsilon,y_\varepsilon,v)p_2,
\Cop^\delta(-\Phi_\alpha(x_\varepsilon,\cdot),u_2)(t_\varepsilon,y_\varepsilon,v),v
)\\
\le&CI_{\eta,\varepsilon}+C\biggl(\eta+\frac{1}{\varepsilon^2}\biggr)(1+|x_\varepsilon|^2+|y_\varepsilon|^2)o_\delta(1)+o_\alpha(1)+\kappa(u_1(t_\varepsilon,x_\varepsilon)-u_2(t_\varepsilon,y_\varepsilon))\\
\quad &+\varpi(|x_\varepsilon-y_\varepsilon|)
+\varpi(CI_{\eta,\varepsilon})
+\varpi\Biggl(CI_{\eta,\varepsilon}+C\biggl(\eta+\frac{1}{\varepsilon^2}\biggr)(1+|x_\varepsilon|^2+|y_\varepsilon|^2)o_\delta(1)+o_\alpha(1)\Biggr),
\end{aligned}
\end{equation}
where \(I_{\eta,\varepsilon}:=\eta(1+|x_\varepsilon|^2)
+|x_\varepsilon-y_\varepsilon|^2/\varepsilon^2\).

Combining \eqref{eq:estimate-hamiltonian} with \eqref{eq:viscosity-inequality}, we have
\begin{equation*}
\begin{aligned}
    -\frac{\beta}{t_\varepsilon^2} &+C\biggl(\eta+\frac{1}{\varepsilon^2}\biggr)(1+|x_\varepsilon|^2+|y_\varepsilon|^2)o_\delta(1)+o_\alpha(1)
    +\varpi(|x_\varepsilon-y_\varepsilon|)+\kappa(u_1(t_\varepsilon,x_\varepsilon)-u_2(t_\varepsilon,y_\varepsilon))\\
    &+CI_{\eta,\varepsilon}+\varpi(CI_{\eta,\varepsilon})+\varpi\Biggl(CI_{\eta,\varepsilon}+C\biggl(\eta+\frac{1}{\varepsilon^2}\biggr)(1+|x_\varepsilon|^2+|y_\varepsilon|^2)o_\delta(1)+o_\alpha(1)\Biggr) \ge 0.
\end{aligned}
\end{equation*}
Taking the lower limits in the order
\[
\alpha\downarrow0,\qquad
\delta\downarrow0,\qquad
\varepsilon\downarrow0,\qquad
\eta\downarrow0,
\]
we have
\[\liminf_{\eta\downarrow0}\kappa(u_1(t_\eta,x_\eta)-u_2(t_\eta,x_\eta))\ge\frac{\beta}{T^2}>0,\]
which contradicts our assumption that \(N>0\). Therefore,
\[
u_1(t,x)\le u_2(t,x),
\qquad (t,x)\in(0,T]\times\mathbb R^{d'}.
\]
\end{proof}

To prove uniqueness, we need to impose an additional uniform bound on \(\gamma\).

\begin{enumerate}[label=\textup{(A\arabic*)},ref=\textup{(A\arabic*)},start=8,wide=0pt]
\item\label{ass:A8} There exists a constant \(M_\gamma>0\) such that \(|\gamma(t,x,v,e)| \le M_\gamma(1\wedge |e|)\) for any \((t,x,v,e)\in[0,T]\times\R^{d'}\times U\times E\).
\end{enumerate}

The next step is to construct a reference function.
This is similar to the standard reference-function argument of \cite[Lemma~3.8]{BarlesBuckdahnPardoux1997}.
Hence we omit the proof.
\begin{lemma}[Reference function]\label{lem:reference-function}
Suppose that \Cref{ass:standing} and \ref{ass:A8} hold. Then for \(\widetilde A>0\), \(S\in(0,T]\) and any \(\rho\ge 2A\), there exists \(C_\rho>0\),
depending only on \(\rho,\widetilde A,T,M_\gamma\) and the constants in \Cref{ass:standing}, such that the reference function
\begin{equation*}
\chi(t,x)
:=
\exp\{
(C_\rho(S-t)+\widetilde A)[
\log((1+|x|^2)^{1/2})+1
]^2
\},
\qquad (t,x)\in[0,S]\times\mathbb R^{d'},
\end{equation*}
satisfies
\[
\chi_t(t,x)+\sup_{v\in U}\bigl\{\Aop\chi(t,x,v)+\Bop\chi(t,x,v)+\kappa\chi(t,x)+\rho |\sigma^\top(t,x,v)D\chi(t,x)|
+\rho(\Cop\chi(t,x,v))^+\bigr\}\le 0
\]
on 
\(\bigl[S-\widetilde A/C_\rho,S\bigr]\cap[0,T] \times \mathbb R^{d'}\).
\end{lemma}

Next we introduce a function class necessary for the uniqueness result.
\[
\Gamma:=\bigl\{u\in C([0,T]\times\mathbb R^{d'}):\exists A_\Gamma>0\text{ such that }|u(t,x)|\le A_\Gamma \exp\bigl\{A_\Gamma[\log(1+|x|^2)+1]^2\bigr\}\bigr\}.
\]

\begin{theorem}[Uniqueness in the growth class]\label{thm:uniqueness-growth}
Suppose that \Cref{ass:standing} and \ref{ass:A8} hold.
If \(u_1,u_2\in\Gamma\) are two viscosity solutions of the HJB equation \eqref{eq:hjb}, then
\[
u_1(t,x)=u_2(t,x),
\qquad
(t,x)\in[0,T]\times\mathbb R^{d'}.
\]
\end{theorem}
\begin{proof}

Let \(A_1,A_2>0\) be such that
\(
    |u_i(t,x)|
    \le
    A_i
    \exp\{
        A_i[\log(1+|x|^2)+1]^2
    \},
    \ i=1,2.
\)
Then fix 
\(
    \widetilde A>4\max\{A_1,A_2\}.
\)
It is immediate that for any \(\varepsilon>0\),
\(u_1-\varepsilon\chi\) is bounded from above and
\(u_2+\varepsilon\chi\) is bounded from below.

It is enough to prove \(u_1\le u_2\). Fix \(\rho\ge 2A\).
We first prove a perturbation claim. For \(S\in(0,T]\), let \(\chi\) be the reference
function in \Cref{lem:reference-function} on the strip
\(I_S\times\mathbb R^{d'}\), where \(I_S:=\bigl[S-\widetilde A/C_\rho,S\bigr]\cap[0,T]\).
For \(\varepsilon>0\), define
\[
u_1^{\varepsilon,\rho}(t,x):=u_1(t,x)-\varepsilon\chi(t,x)-2\varpi\Bigl(\frac{2A}{\rho}\Bigr)(S-t).
\]
We claim that \(u_1^{\varepsilon,\rho}\) is a viscosity subsolution on
\(I_S\times\mathbb R^{d'}\).
To see it, take \(\phi\in C^{1,2}(I_S\times\mathbb R^{d'})\) such that
\(u_1^{\varepsilon,\rho}-\phi\)
attains a global maximum at \((t,x)\in I_S^o\times\mathbb R^{d'}\), where \(I_S^o\) is the interior of \(I_S\). Define
\[
\widetilde\phi(s,y):=\phi(s,y)+\varepsilon\chi(s,y)+2\varpi\Bigl(\frac{2A}{\rho}\Bigr)(S-s).
\]

Then \(u_1-\widetilde\phi\) attains its maximum over
\(I_S\times\R^{d'}\) at \((t,x)\).
We may choose \(\delta>0\) sufficiently small such that
\(B((t,x),R_{\delta,x})\subset I_S\times\R^{d'}.\)
Applying \Cref{def:viscosity-local-i}, we obtain
\[
\widetilde\phi_t(t,x)
+\sup_{v\in U}\bigl\{
\Aop\widetilde\phi(t,x,v)+\Bop^\delta(\widetilde\phi,u_1)(t,x,v)
+f\bigl(t,x,u_1(t,x),\sigma^\top(t,x,v)D\widetilde\phi(t,x),
\Cop^\delta(\widetilde\phi,u_1)(t,x,v),v\bigr)
\bigr\}\ge 0.
\]
By the maximality of
\(u_1-\widetilde\phi\), it follows that
\[
\widetilde\phi_t(t,x)
+\sup_{v\in U}\bigl\{
\Aop\widetilde\phi(t,x,v)+\Bop\widetilde\phi(t,x,v)
+f\bigl(t,x,u_1(t,x),\sigma^\top(t,x,v)D\widetilde\phi(t,x),\Cop\widetilde\phi(t,x,v),v\bigr)
\bigr\} \ge 0.
\]

On the other hand,
\[
\begin{aligned}
& f(t,x,u_1(t,x),\sigma^\top(t,x,v)D\widetilde\phi(t,x),\Cop\widetilde\phi(t,x,v),v)-f\bigl(t,x,u_1^{\varepsilon,\rho}(t,x),\sigma^\top(t,x,v)D\phi(t,x),\Cop\phi(t,x,v),v\bigr)
\\ \le&
\kappa\Bigl(\varepsilon\chi(t,x)+2\varpi\Bigl(\frac{2A}{\rho}\Bigr)(S-t)\Bigr)
+\varpi(\varepsilon |\sigma^\top(t,x,v)D\chi(t,x)|)
+\varpi(\varepsilon(\Cop\chi(t,x,v))^+)
\\ \le&\varepsilon \kappa \chi(t,x)+\varepsilon \rho|\sigma^\top(t,x,v)D\chi(t,x)|
+\varepsilon \rho(\Cop\chi(t,x,v))^+ +2\varpi\Bigl(\frac{2A}{\rho}\Bigr).
\end{aligned}
\]
Hence
\[
\begin{aligned}
&\widetilde\phi_t+\Aop\widetilde\phi+\Bop\widetilde\phi
+f(t,x,u_1,\sigma^\top D\widetilde\phi,\Cop\widetilde\phi,v)
\\
\le&\phi_t+\Aop\phi+\Bop\phi
+f\bigl(t,x,u_1^{\varepsilon,\rho},\sigma^\top D\phi,\Cop\phi,v\bigr)+\varepsilon
\bigl[\chi_t+\Aop\chi+\Bop\chi+\kappa\chi+\rho |\sigma^\top D\chi|+\rho(\Cop\chi)^+\bigr].
\end{aligned}
\]
By \Cref{lem:reference-function}, the bracketed term is nonpositive.
Taking the supremum over \(v\in U\), we obtain
\[
\phi_t(t,x)
+
\sup_{v\in U}
\bigl\{
\Aop\phi(t,x,v)
+
\Bop\phi(t,x,v)
+
f\bigl(t,x,u_1^{\varepsilon,\rho}(t,x),\sigma^\top(t,x,v)D\phi(t,x),\Cop\phi(t,x,v),v\bigr)
\bigr\}
\ge0.
\]
Thus \(u_1^{\varepsilon,\rho}\) is a viscosity subsolution.

Similarly, for any constant \(r\ge0\), define
\[
u_{2,r}^{\varepsilon,\rho}(t,x)
:=
u_2(t,x)+\varepsilon\chi(t,x)+2\varpi\Bigl(\frac{2A}{\rho}\Bigr)(S-t)+r,
\]
and we can prove that \(u_{2,r}^{\varepsilon,\rho}\) is a viscosity supersolution on
\(I_S\times\mathbb R^{d'}\).

We now apply \Cref{thm:comparison} on successive time strips. Set
\[
S_0:=T,
\qquad
S_{j+1}:=\biggl(S_j-\frac{\widetilde A}{C_\rho}\biggr)^+,
\qquad j\ge0.
\]
Since \(C_\rho<\infty\), this sequence reaches \(0\) after finitely many steps.
At \(S_0=T\), the terminal condition implies
\(u_1(T,x)\le u_2(T,x)\),  \(x\in\mathbb R^{d'}\).
On the strip \([S_1,S_0]\), the functions
\(u_1^{\varepsilon,\rho}\) and \(u_{2,0}^{\varepsilon,\rho}\)
are a bounded-above subsolution and a bounded-below supersolution, respectively. Moreover,
\[
u_1^{\varepsilon,\rho}(S_0,x)
=
u_1(S_0,x)-\varepsilon\chi(S_0,x)
\le
u_2(S_0,x)+\varepsilon\chi(S_0,x)
=
u_{2,0}^{\varepsilon,\rho}(S_0,x).
\]
Applying \Cref{thm:comparison} on \([S_1,S_0]\), we obtain
\[
u_1(t,x)-u_2(t,x)
\le
2\varepsilon\chi(t,x)
+
4\varpi\Bigl(\frac{2A}{\rho}\Bigr)(S_0-t),
\qquad
(t,x)\in[S_1,S_0]\times\mathbb R^{d'}.
\]
Taking \(\varepsilon\downarrow0\), we get
\[
u_1(t,x)-u_2(t,x)
\le
4\varpi\Bigl(\frac{2A}{\rho}\Bigr)(S_0-t),
\qquad
(t,x)\in[S_1,S_0]\times\mathbb R^{d'}.
\]
In particular,
\[
u_1(S_1,x)-u_2(S_1,x)
\le
4\varpi\Bigl(\frac{2A}{\rho}\Bigr)(S_0-S_1).
\]
The proof is by induction. Suppose that for some \(j\ge0\),
\[
u_1(S_j,x)-u_2(S_j,x)
\le
r_j,
\qquad x\in\mathbb R^{d'},
\]
where \(r_j:=4\varpi(2A/\rho)(T-S_j)\).
On the strip \([S_{j+1},S_j]\), applying \Cref{thm:comparison} again to
\[
u_1-\varepsilon\chi-2\varpi\Bigl(\frac{2A}{\rho}\Bigr)(S_j-t)
\ \ \ \text{and}\ \ \ 
u_2+\varepsilon\chi+2\varpi\Bigl(\frac{2A}{\rho}\Bigr)(S_j-t)+r_j,
\]
and noticing that the terminal inequality at \(t=S_j\) is satisfied, we have
\[
u_1(t,x)-u_2(t,x)
\le
r_j
+
2\varepsilon\chi(t,x)
+
4\varpi\Bigl(\frac{2A}{\rho}\Bigr)(S_j-t).
\]
Taking \(\varepsilon\downarrow0\), we get
\[
u_1(t,x)-u_2(t,x)
\le
r_j+4\varpi\Bigl(\frac{2A}{\rho}\Bigr)(S_j-t)
\le
4\varpi\Bigl(\frac{2A}{\rho}\Bigr)(T-t),
\qquad
(t,x)\in[S_{j+1},S_j]\times\mathbb R^{d'}.
\]
After finitely many steps, the whole interval \([0,T]\) is covered. Hence
\[
u_1(t,x)-u_2(t,x)
\le
4\varpi\Bigl(\frac{2A}{\rho}\Bigr) T,
\qquad
(t,x)\in[0,T]\times\mathbb R^{d'}.
\]
Finally, taking \(\rho\to\infty\), we obtain
\[
u_1(t,x)\le u_2(t,x),
\qquad
(t,x)\in[0,T]\times\mathbb R^{d'}.
\]
The reverse inequality can be obtained by interchanging \(u_1\) and \(u_2\). Therefore 
\(u_1=u_2\) on \([0,T]\times\mathbb R^{d'}\).
\end{proof}

\section*{Data availability}
Data sharing is not applicable to this article as no datasets were generated or analysed during this theoretical study.

\end{document}